%% file: eulergothic_arxiv.tex
\documentclass{amsart}
\pdfoutput=1

\usepackage{amsthm,amsmath,amssymb,latexsym,amsxtra,amscd,amssymb,color,mathtools,xspace}
\usepackage{mathtools}
\usepackage{units}
\usepackage{enumitem}
\usepackage{array}
\usepackage{booktabs}
\usepackage{calc}
\usepackage{pgf}
\usepackage{scalefnt}
\usepackage{listings}
\usepackage{multicol}
\usepackage{longtable}
\allowdisplaybreaks

\usepackage{bm}

\usepackage{placeins}
\usepackage[obeyspaces,hyphens,spaces]{url}
\usepackage{letltxmacro}

\usepackage{graphicx,color}
\graphicspath{{img/}}

\usepackage{tikz}
\usetikzlibrary{matrix,arrows,calc}
\usetikzlibrary[patterns]
\usetikzlibrary{decorations.markings}
\usetikzlibrary{decorations.pathmorphing,patterns}

\usepackage{csquotes}

\usepackage{enumitem}

\usepackage{amsthm}
\usepackage{hyperref}

\usepackage{aliascnt}

\input{macros}

\begin{document}

\title[Gothic Teichm\"{u}ller curves] {Euler characteristics of \\
  Gothic Teichm\"{u}ller curves}

\author{Martin M\"oller}
\address{
Institut f\"ur Mathematik, Goethe--Universit\"at Frankfurt,
Robert-Mayer-Str. 6--8,
60325 Frankfurt am Main, Germany
}
\email{moeller@math.uni-frankfurt.de}

\author{David Torres-Teigell}
\address{
Institut f\"ur Mathematik, Goethe--Universit\"at Frankfurt,
Robert-Mayer-Str. 6--8,
60325 Frankfurt am Main, Germany
}
\email{torres@math.uni-frankfurt.de}
\begin{abstract}
We compute the Euler characteristics of the recently discovered
series of Gothic Teichm\"{u}ller curves. The main tool is the
construction of `Gothic' Hilbert modular forms vanishing at the
images of these Teichm\"{u}ller curves.
\par
Contrary to all previously known examples, the Euler characteristics
is not proportional to the Euler characteristic of the ambient
Hilbert modular surfaces. This results in interesting `varying'
phenomena for Lyapunov exponents.
\end{abstract}

\maketitle
\tableofcontents
\noindent
\SaveTocDepth{1}

\input{sec_intro} 

\input{sec_gothic_def}


\input{sec_complementary}

\input{sec_HMVAR}
\input{sec_linebundles}


\input{sec_thetagothic}

\input{sec_modularcurves}

\input{sec_cusps}

\input{sec_examples}


\input{sec_asymptotics}

\input{sec_lyapunov}

\newpage
\begin{appendix}
\input{sec_tables}
 \end{appendix}


\bibliographystyle{alpha}

\input{sources.bbl}

\end{document}

%% file: macros.tex
\lstdefinelanguage{Sage}[]{Python}
{morekeywords={False,sage,True},sensitive=true}
\def\equationautorefname~#1\null{(#1)\null}

\theoremstyle{plain}
\newtheorem{theorem}{Theorem}[section]
\newaliascnt{lemma}{theorem}
\newtheorem{lemma}[lemma]{Lemma}
\aliascntresetthe{lemma}
\newaliascnt{prop}{theorem}
\newtheorem{prop}[prop]{Proposition}
\aliascntresetthe{prop}
\newaliascnt{cor}{theorem}
\newtheorem{cor}[cor]{Corollary}
\aliascntresetthe{cor}
\newaliascnt{conj}{theorem}

\aliascntresetthe{conj}
\newaliascnt{dfn}{theorem}

\aliascntresetthe{dfn}
\newaliascnt{ex}{theorem}

\aliascntresetthe{ex}

\makeatletter%
\protected\def\W{\@ifnextchar[{\W@arg}{W_D}}
\def\W@arg[#1]{%
  \def\@D{D}%
  \def\c@W##1,##2,##3\@nil{\def\c@D{##1}\def\n@D{##2}}\c@W#1,,\@nil%
  \def\@W##1=##2=##3\@nil{\def\@DD{##1}\def\@DDD{##2}}\expandafter\@W\c@D==\@nil%
  \ifx\@D\@DD W_{\@DDD}%
    \expandafter\ifx\n@D\relax\relax\else(\n@D)\fi%
  \else W_D(#1)\fi%
}
\protected\def\parexp#1{\@ifnextchar^{(#1)}{#1}}
\protected\def\parind#1{\@ifnextchar_{(#1)}{#1}}

\def\defbb#1{\expandafter\def\csname b#1\endcsname{\mathbb{#1}}}
\def\defcal#1{\expandafter\def\csname c#1\endcsname{\mathcal{#1}}}
\def\deffrak#1{\expandafter\def\csname frak#1\endcsname{\mathfrak{#1}}}
\def\defop#1{\expandafter\def\csname#1\endcsname{\operatorname{#1}}}
\def\defbold#1{\expandafter\def\csname bb#1\endcsname{\mathbf{#1}}}

\def\defcals#1{\@defcals#1\@nil}
\def\@defcals#1{\ifx#1\@nil\else\defcal{#1}\expandafter\@defcals\fi}
\def\deffraks#1{\@deffraks#1\@nil}
\def\@deffraks#1{\ifx#1\@nil\else\deffrak{#1}\expandafter\@deffraks\fi}
\def\defbbs#1{\@defbbs#1\@nil}
\def\@defbbs#1{\ifx#1\@nil\else\defbb{#1}\expandafter\@defbbs\fi}
\def\defops#1{\@defops#1,\@nil}
\def\@defops#1,#2\@nil{\if\relax#1\relax\else\defop{#1}\fi\if\relax#2\relax\else\expandafter\@defops#2\@nil\fi}
\def\defbolds#1{\@defbolds#1\@nil}
\def\@defbolds#1{\ifx#1\@nil\else\defbold{#1}\expandafter\@defbolds\fi}
\makeatother%
\def\mat#1#2#3#4{\begin{pmatrix}#1&#2\\#3&#4\\ \end{pmatrix}} 
\def\sm#1#2#3#4{\bigl(\smallmatrix#1&#2\\#3&#4\endsmallmatrix\bigr)} 

\defbbs{ZHQCNPARDV}
\defcals{OCFGHRMXPYZEVULA}
\deffraks{abcdefopqgmB}
\defops{SL,mod,Spec,Re,Gal,Tr,tr,End,GL,Hom,Kern,Bild,PGL,Aut,id,Pic,div,Supp,Fix,Jac,im,Im,PT,Pic,Id,Prym,diag,mult,tr,Sp,pr,Red}
\defbolds{uzvxeabcd}

\def\i{\mathrm{i}}

\def\A{A}
\def\B{B}
\def\piA{\pi_{A}}
\def\piB{\pi_{B}}
\def\p{p}
\def\abel{\overline{\varphi}}
\def\pabel{\varphi}
\def\={\;=\;}  \def\+{\,+\,} 
\def\ve{\epsilon}
\def\be{\begin{equation}}   \def\ee{\end{equation}}     
\def\bes{\begin{equation*}}    \def\ees{\end{equation*}}
\def\ba{\be\begin{aligned}} \def\ea{\end{aligned}\ee}   
\def\bas{\bes\begin{aligned}}  \def\eas{\end{aligned}\ees}
\newcommand{\legendre}[2]{\left(\frac{#1}{#2}\right)}

\newcommand{\bftau}{{\boldsymbol{\tau}}}
\newcommand{\Teichmuller}{Teich\-m\"uller~}
\newcommand{\RED}{{\Red_{23}}}
\newcommand{\tRED}{{\widetilde{\Red}_{23}}}
\newcommand{\REDIRR}{{\Red_{16}}}
\newcommand{\MB}{{\mathcal{B}}}
\newcommand{\avdg}{M}

\DeclarePairedDelimiterX{\inbrack}[1]{[}{]}{%
  #1}
\newcommand{\sizebrackets}[2][0]{%
  \ifcase#1\relax
    \inbrack{#2}\or         
    \inbrack[\big]{#2}\or   
    \inbrack[\Big]{#2}\or   
    \inbrack[\bigg]{#2}\or  
    \inbrack[\Bigg]{#2}     
  \else
    \inbrack*{#2}
  \fi}

\renewcommand{\th}[3]{
\pgfmathsetmacro \sp {{2*(#3-2)}}
\vartheta\sizebrackets[#3]{\begin{matrix} #1 \\[\sp pt] #2 \end{matrix}}
}

\newcounter{savedtocdepth}
\newcommand*{\SaveTocDepth}[1]{%
  \addtocontents{toc}{%
    \protect\setcounter{savedtocdepth}{\protect\value{tocdepth}}%
    \protect\setcounter{tocdepth}{#1}%
  }%
}

\newcommand{\moduli}[1][g]{{\mathcal M}_{#1}}
\newcommand{\omoduli}[1][g]{{\Omega \mathcal M}_{#1}}

\definecolor{lgray}{RGB}{214, 214, 214}
\definecolor{dgray}{RGB}{115, 115, 115}

\newcommand*{\x}{1}%
\newcommand*{\y}{1}%

\renewcommand*{\x}{4/3}%
\renewcommand*{\y}{4/sqrt(3)}%

\newcommand*{\ShiftP}{($({\x*(3/4*sqrt(3) + 3/4)}, {(1/4*sqrt(3) + 3/4)*\y}) - ({\x*((-3)/4)}, {(1/4*sqrt(3))*\y})$)}%
\newcommand*{\ShiftQ}{($({\x*((-1)/4*sqrt(3) - 3/4)}, {(3/4*sqrt(3) + 3/4)*\y}) -  ({\x*(1/2*sqrt(3))}, {(0)*\y})$)}%

\pgfdeclarepatternformonly{wide horizontal lines}{\pgfpointorigin}{\pgfqpoint{100pt}{1pt}}
{\pgfqpoint{100pt}{6pt}}
{
	\pgfsetlinewidth{0.4pt}
	\pgfpathmoveto{\pgfqpoint{0pt}{0.5pt}}
	\pgfpathlineto{\pgfqpoint{100pt}{0.5pt}}
	\pgfusepath{stroke}
}

\pgfdeclarepatternformonly{closely horizontal lines}{\pgfpointorigin}{\pgfqpoint{100pt}{1pt}}
{\pgfqpoint{100pt}{2pt}}
{
	\pgfsetlinewidth{0.4pt}
	\pgfpathmoveto{\pgfqpoint{0pt}{0.5pt}}
	\pgfpathlineto{\pgfqpoint{100pt}{0.5pt}}
	\pgfusepath{stroke}
}

\pgfdeclarepatternformonly{wide vertical lines}{\pgfpointorigin}{\pgfqpoint{1pt}{100pt}}
{\pgfqpoint{6pt}{100pt}}
{
	\pgfsetlinewidth{0.4pt}
	\pgfpathmoveto{\pgfqpoint{0.5pt}{0pt}}
	\pgfpathlineto{\pgfqpoint{0.5pt}{100pt}}
	\pgfusepath{stroke}
}

\pgfdeclarepatternformonly{closely vertical lines}{\pgfpointorigin}{\pgfqpoint{1pt}{100pt}}
{\pgfqpoint{3pt}{100pt}}
{
	\pgfsetlinewidth{0.4pt}
	\pgfpathmoveto{\pgfqpoint{0.5pt}{0pt}}
	\pgfpathlineto{\pgfqpoint{0.5pt}{100pt}}
	\pgfusepath{stroke}
}

%% file: sec_intro.tex
\section{Introduction}

\Teichmuller curves are complex geodesics in the moduli space
of curves~$\moduli[g]$. They arise as the $\SL_2(\bR)$-orbits of flat surfaces
with optimal dynamics, called Veech surfaces. If the Veech surface
is not obtained by a covering construction from a lower genus surface,
it is called primitive and the resulting \Teichmuller curve is
called primitive too. There are very few constructions of primitive
Teichm\"uller curves (see \cite{moellerICM} for a list of known examples).
Each infinite collection of primitive \Teichmuller curves in a
fixed genus stems from an invariant submanifold `like the minimal
stratum $\omoduli[2](2)$' in genus two (see Section~\ref{sec:semicath}), by the finiteness results from Eskin-Filip-Wright~\cite{EFWfinite}. While the geometry of $\omoduli[2](2)$ and of the Prym
loci is well-understood now, the geometry of the two invariant submanifolds
`like $\omoduli[2](2)$' recently discovered by
Eskin-McMullen-Mukamel-Wright in~\cite{EMMW} is basically unexplored. Here we
focus on the Gothic locus~$\Omega G \subset \omoduli[4](2,2,2)$ of
flat genus four surfaces, introduced already in~\cite{MMW}.
\par
While interest in \Teichmuller curves originates from dynamics,
their geometry is strongly determined by modular forms. \Teichmuller
curves in an infinite series of fixed genus always map via the
Torelli map to the locus of real multiplication, i.e.\ to a
Hilbert modular surface (\cite{moeller06} together with \cite{EFWfinite}).
Conversely, the intersection of~$\Omega G$ with the locus of real
multiplication by the order~$\cO_D$ is a union of \Teichmuller
curves~$G_D$. These  \Teichmuller curves are primitive if and only if $D$
is not a square, which we assume in the rest of this paper. The
modular forms in question are thus Hilbert modular forms, supposed
to cut out the \Teichmuller curves~$G_D$ inside the Hilbert modular surface.
\par
Contrary to the expectation from the situation in genus two and in the
Prym loci of genus three and four, there is no Hilbert modular form whose
vanishing locus is precisely equal to the Gothic \Teichmuller curves~$G_D$!
Yet, there is a `Gothic' Hilbert modular form~$\cG_D$ whose vanishing locus
is only slightly larger than $G_D$, the difference being a collection of
modular curves, whose parameters can all be computed. 
\par
\medskip
In order to state the results, we roughly recall the definition
of~$\Omega G$, see Section~\ref{sec:defGothic} for more details.
A flat surface $(X,\omega)$ in the stratum~$\omoduli[4](2,2,2)$
is Gothic if it admits an involution~$J$ leaving~$\omega$
anti-invariant and fixing the zeros of~$\omega$ and an `odd' degree
three map $X \to B$ to an elliptic curve~$B$ mapping all the zeros
to a single point. The involution~$J$ induces a degree two map
$X \to A$ to another elliptic curve~$A$. The complement of both~$A$
and~$B$ in the Jacobian of~$X$ inherits a polarization of type~$(1,6)$
as we show in Section~\ref{sec:complementary}. Consequently, the
number~$6$ plays
a prominent role in the paper: the Gothic \Teichmuller curves~$G_D$
naturally live on Hilbert modular surfaces $X_D(\frakb)$ where
$\frakb$ is an $\cO_D$-ideal of norm~$6$.
\par
Our first goal is to give a natural decomposition of~$G_D$ into (perhaps still reducible)
components and to compute explicitly their Euler characteristics. They can be written in terms of Euler characteristics of those
Hilbert modular surfaces and of the reducible locus $\RED$,
parametrizing $(2,3)$-polarized products of elliptic curves with real
multiplication by $\cO_{D}$.
\par
\begin{theorem} \label{thm:introVol} Let $D$ be a non-square discriminant. The Gothic \Teichmuller curve $G_D$ is non-empty if and only if $D \equiv 0,1,4,9,12,16\bmod{24}$. \\ 
In this case, $G_{D}$ consists of different sub-curves $G_{D}(\frakb)$ corresponding to different $\cO_{D}$-ideals $\frakb$ of norm 6. The Euler characteristics of all these sub-curves agree and are equal to
\bes \label{eq:volcomp}
-\chi(G_D(\frakb)) =  \frac32 \, \chi(X_D(\frakb)) + 2\, \chi(\RED(\frakb))\,.
\ees
\end{theorem}
\par
We give a completely explicit formula in Theorem~\ref{thm:volumes}, and a
table for small discriminants can be found at the end of the paper.
The Euler characteristic of the Hilbert modular surface~$X_D(\frakb)$
is equal to the Euler characteristic of a standard Hilbert modular surface
if~$D$ is fundamental, and differs by a simple factor in general,
see Proposition~\ref{prop:vol16} for the complete formula. In any case, $\chi(X_{D}(\frakb))$ is independent of the choice of the ideal of norm 6. We strongly suspect the sub-curves~$G_D(\frakb)$ defined in the
theorem to be irreducible but we do not attempt to prove this here.
\par
\medskip
The presence of modular curves in the vanishing locus of the
Hilbert modular form~$\cG_D$ has another consequence that makes
characteristic invariants of the Gothic \Teichmuller curves
behave differently than all the examples known so far. We phrase this in
terms of Lyapunov exponents in Section~\ref{sec:volandlyap} and restate it
geometrically here.
\par
Teichm\"uller curves are Kobayashi geodesic algebraic curves~$C$ in Hilbert
modular surfaces. If $z \mapsto (z,\varphi(z))$ is the universal covering map
of a Kobayashi geodesic, then for any $M \in \GL_2(\bQ(\sqrt{D}))$ the map
$z \mapsto (Mz,M^\sigma\varphi(z))$ descends to another Kobayashi geodesic.
All modular curves arise by this twisting procedure from the diagonal
and the twists of Teichm\"uller curves are interesting special curves
on Hilbert modular surfaces. However, this twisting does not change
the most basic algebraic invariant, 
\bes
\lambda_2(C) \= (C \cdot [\omega_2])\,/\,(C \cdot [\omega_1])\,,
\ees
where $[\omega_i]$ are the foliation classes on the Hilbert modular surface.
For modular curves $\lambda_2=1$ and, in general, the list of known $\lambda_2(C)$ of Kobayashi geodesics $C$ was a rather short (and finite) list (see the summary in \cite[Section~1]{MZ}). As a consequence of the
decomposition of $\{\cG_D = 0\}$ into several components
we obtain:
\par
\begin{cor}
The sequence of invariants $\lambda_2(G_D)$ is infinite and tends
to $3/13$ for $D \to \infty$.
\end{cor}
\par 
This corollary is proved in the equivalent formulation of Theorem~\ref{thm:LyapOmegaG}, see also Proposition~\ref{prop:LyapGD}. It is an open question whether for a {\em fixed} Hilbert modular surface
the set of $\lambda_2(C)$ for all its Kobayashi geodesics~$C$ is finite
or infinite.
\par
\medskip
We next summarize the main steps in the proof of Theorem~\ref{thm:introVol}
and explain the origin of the Gothic modular form~$\cG_D$. Analyzing the
definition of the Gothic locus (Section~\ref{sec:ExGothic}), we obtain that the image of
a Gothic Veech surface in its $(1,6)$-polarized Prym abelian surface is a curve with a triple point at the origin and horizontal tangents at three
non-zero two-torsion points (Section~\ref{sec:GothMF}). To construct these images as the
vanishing locus of a theta function, we need to impose $5$ conditions,
two stemming from the multiplicity at the origin and the rest from the behaviour at the two-torsion points.
The (odd) theta functions vary in a $3$-dimensional
projective space, so we can impose the first three conditions and, by restricting to a divisor $\{\cG_D =0\}$
in the Hilbert modular surface, we can also satisfy the last two conditions. Teichm\"uller curves exist due to dimension miracles. From our point of view
this is manifested by the last two conditions holding simultaneously along $\{\cG_D =0\}$,
due to theta value relations at $2$-torsion points (Section~\ref{sec:linbd}).
\par
Contrary to the previous known cases in $\omoduli[2](2)$ and the Prym locus, the vanishing locus of the Gothic modular form $\cG_{D}$ contains some `spurious' components apart from the Gothic \Teichmuller curves. These components form the $(2,3)$-reducible locus, points in the Hilbert modular surface corresponding to products of elliptic curves with the natural $(2,3)$-product polarization (Section~\ref{sec:modularcurves}). By studying the vanishing order of the modular form along both the Gothic \Teichmuller curve and the reducible locus (Section~\ref{sec:divGD}), we can finally relate their Euler characteristics with the Euler characteristic of the Hilbert modular surface in which they live (Section~\ref{sec:volandlyap}). This also allows us to give a formula for the Lyapunov exponents of the individual Gothic \Teichmuller curves and to compute those of the Gothic locus.
\par
{\bf Acknowledgements:} The authors thank Ronen Mukamel for sharing
insights, in particular his program to compute Veech groups
(\cite{mukamelfundamental})
that provided valuable cross-checks, recorded in the table in the appendix. 
The authors also thank Don Zagier for useful conversations, in particular
in connection with Section~\ref{sec:AsyDiv}, and the referee for many helpful suggestions.\\
The authors acknowledge support from the LOEWE-Schwerpunkt ``Uniformisierte
Strukturen in Arithmetik und Geometrie'' and from the DFG-Projekt ``Classification of Teichm\"{u}ller curves MO 1884/2--1''.

%
%

%% file: sec_gothic_def.tex
\section{Examples of Gothic \Teichmuller curves} \label{sec:ExGothic}

In this section we introduce the Gothic locus and the 
Gothic \Teichmuller curves, following \cite{MMW}. Not all
Gothic \Teichmuller curves can be presented in the shape of
a Gothic cathedral. In fact, the simplest example of a Gothic
\Teichmuller curve already appeared in work of Ward~\cite{ward}
on triangular billiards.

\subsection{The Gothic locus} \label{sec:defGothic}

Given a Riemann surface $X$ with an involution~$J$ we say that
a map $\piB:X \to B$ is {\em odd}, if there exists an involution $j:\B\to \B$ such that 
$\piB\circ J=j\circ\piB$. Following \cite{MMW} we define the
{\em Gothic locus}~$\Omega G$ to be the set of Riemann surfaces
$(X,\omega) \in \Omega\cM_{4}(2^{3},0^{3})$ such that
\begin{itemize}
\item[i)] there exists an involution $J \in \Aut(X)$ whose
fixed points are the six marked points, the zeros $Z = \cZ(\omega) 
= \{z_1,z_2,z_3\}$ and the marked regular points  $P = \{p_1,p_2,p_3\}$,
\item[ii)] the one-form $\omega$ is $J$-antiinvariant, that is $J^{*}\omega=-\omega$, and
\item[iii)] there exists a genus one curve~$B$ and an odd map $\pi_{B}:X\to B$ of 
degree~$3$ such that $|\piB(Z)|=1$.
\end{itemize}
\par
Every flat surface $(X,\omega)\in\Omega G$ in the  Gothic locus thus
comes with  maps
\be \label{eq:Gothic}
\begin{tikzpicture}[ampersand replacement=\&]
\matrix (m) [matrix of math nodes, row sep=1.75em, column sep=3.5em, text height=1.5ex, text depth=0.25ex]
{ \& X \& \\
\A \& \& \B \\
\& \bP^{1} \& \\
};
\path[->,font=\scriptsize]
(m-1-2) edge [above] node {$\piA$} (m-2-1)
(m-1-2) edge [above] node {$\piB$} (m-2-3)
(m-2-1) edge [above] node {$\p$} (m-3-2)
(m-2-3) edge [above] node {$r$} (m-3-2)
(m-1-2) edge [right] node {$h$} (m-3-2)
;
\end{tikzpicture}\ee
where
\begin{itemize}
	\item $\piA:X\to X/J\cong \A$ is of degree 2;
	\item $\piB:X\to \B$ is an odd, degree 3 ramified covering such
 that $|\piB(Z)|=1$, 
	\item $r:\B\to \B/j\cong\bP^{1}$ is the quotient map; and
	\item $\p:\A\to \bP^{1}$ is the degree 3 ramified covering 
that makes the diagram commutative.
\end{itemize}
\par
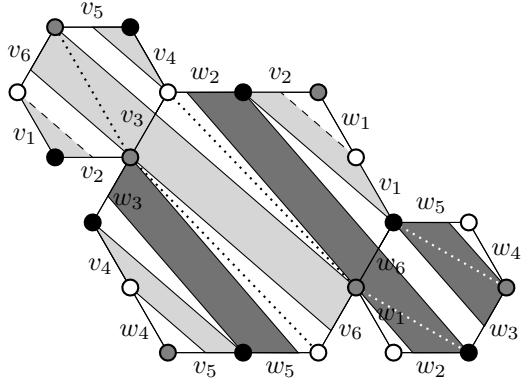
\begin{figure}[htb]
\centering
\input{cylinder_preperation_rotated}
\caption{The hexagon form in the Gothic locus (from~\cite{MMW})}
\label{cap:Gothichex}
\end{figure}
These maps can be illustrated on the hexagon form in
Figure~\ref{cap:Gothichex}. It admits an automorphism~$R$
of order~$6$ with $R^* \omega = \zeta_6 \omega$. Then $J = R^3$
and $\pi_A$ and $\pi_B$ are the quotients by $R^3$ and $R^2$
respectively. Note, however, that the map~$\pi_B$ will not be Galois in general.
The reason for the definition is that~$\Omega G$ turns out to be an unexpected
$\SL_2(\bR)$-orbit closure.
\par
\begin{theorem}[\cite{MMW}]
The Gothic locus is a closed irreducible variety of dimension~$4$, 
locally defined by linear equation in period coordinates.
\end{theorem}
\par
In fact, $v_1,\ldots,v_6, w_1,\ldots,w_6$ are periods on the $10$-dimensional
space $\omoduli[4](2,2,2)$, and $\sum_{i=1}^6 v_i = 0 = \sum_{i=1}^6 w_i$ by
construction. In fact, $v_1,\ldots,v_5, w_1,\ldots,w_5$ form a coordinate
system. In this coordinates $\Omega G$ is cut out by the conditions
\ba \label{eq:lineareq}
v_{i+3} \= -v_i, \quad w_{i+3} \= -w_i, \quad 
v_1+v_3+v_5 \= 0, \quad w_1+w_3+w_5 \= 0
\ea
for $i=1,2,3$.
\par
\medskip
The branch points of the maps in the diagram~\eqref{eq:Gothic} give a
collection of special points. We introduce notation for later use.
Given a point $x\in X$ we will denote the other points in the same
$\pi_B$-fiber by $\piB^{-1}(\piB(x))=\{x\eqqcolon x^{(1)},x^{(2)}, x^{(3)} \}$
and $h=r\circ\piB=\p\circ\piA$.
The preimages of the ramification points of~$h$ and their behavior under 
the maps $\piA$ and $\piB$ can be described in the following way.
\begin{itemize}
	\item The image point $e_{4}'=\piB(Z)$ is fixed by $j$, since each
$z_i$ is fixed by~$J$. This point is therefore sent by $r$ to a 
ramification point $e_{4}=r(e_{4}')$ of~$h$. In particular we can choose the
group law on~$B$ such that~$e_{4}'$ agrees with the origin $O$.
	\item The image points $e_{i}'=\piB(p_{i})$ for $i=1,2,3$ are also fixed by $j$, giving rise to the other three points of order 2 in $\B$. Their preimages under $\piB$ are given by $\piB^{-1}(e_{i}')=\{p_{i}, q_{i}, J(q_{i})\}$.
	\item There exist  three other ramification points of the map $h$, 
among the preimages of which there exist points $\{y_{i}, J(y_{i})\}$ for $i=1,2,3$ with ramification index~$2$ with respect to $h$ each. 
\end{itemize}
\par
\be\label{eq:hpreim}
\begin{tikzpicture}
\matrix (m) [matrix of math nodes, row sep=1.75em, column sep=2.0em, text height=1.5ex, text depth=0.25ex]
{ & \{p_i, q_i, J(q_i) \} & \\
\{\bar{p}_i, \bar{q}_i  \}  & & e_i' \\
& e_i & \\
};
\path[->,font=\scriptsize]
(m-1-2) edge [above] node {$\piA$} (m-2-1)
(m-1-2) edge [above] node {$\piB$} (m-2-3)
(m-2-1) edge [above] node {$\p$} (m-3-2)
(m-2-3) edge [above] node {$r$} (m-3-2)
(m-1-2) edge [right] node {$h$} (m-3-2)
;
\end{tikzpicture}
\begin{tikzpicture}
\matrix (m) [matrix of math nodes, row sep=1.75em, column sep=2.0em, text height=1.5ex, text depth=0.25ex]
{ & \{z_1, z_2, z_3  \} & \\
\{\bar{z}_1, \bar{z}_2, \bar{z}_3  \}  & & e_4' \\
& e_4 & \\
};
\path[->,font=\scriptsize]
(m-1-2) edge [above] node {$\piA$} (m-2-1)
(m-1-2) edge [above] node {$\piB$} (m-2-3)
(m-2-1) edge [above] node {$\p$} (m-3-2)
(m-2-3) edge [above] node {$r$} (m-3-2)
(m-1-2) edge [right] node {$h$} (m-3-2)
;

\end{tikzpicture}\ee
\par
Recall that the stratum $\omoduli[g](2,2,2)$ has two connected
components, distinguished by the parity of the spin structure. One can take a flat surface in $\Omega G$ (e.g. the hexagon form) and compute the winding numbers of a symplectic basis to prove that the Gothic locus lies in the component $\omoduli[g]^{even}(2,2,2)$ with even spin structure (see also the 
argument using $\theta_{\rm null}$  in \cite[Section~4]{MMW}). 
We will however not use this fact when cutting out in Section~\ref{sec:GothMF}
the image of the Veech surfaces in their Prym varieties with theta functions.
\par
The one-form~$\omega$ obviously belongs to the tangent space
to a three-dimensional subvariety of~$\Jac(X)$, the complement of~$A$,
since $\omega$ is $J$-invariant. We can reduce the considerations to
abelian surfaces, the complement of both~$A$ and~$B$ thanks to the 
following observation.
\par
\begin{lemma}[\cite{MMW}]
For $(X,\omega) \in \Omega G$ the $\pi_B$-pushforward is zero.
\end{lemma}
\par
\begin{proof}
The differential $(\pi_B)_*(\omega)$ vanishes at $e_4'$, since all
the $\pi_B$-preimages of that point are zeros of $\omega$, and 
this push-forward differential is holomorphic. 
On the elliptic curve~$B$ this implies $(\pi_B)_*(\omega) = 0$.
\end{proof}
\par

\subsection{Gothic \Teichmuller curves: Cathedrals and semiregular hexagons} \label{sec:semicath}

The Gothic locus~$\Omega G$ is `like $\omoduli[2](2)$' in a precise
sense: it is an affine invariant manifold of dimension four and
rank~two (in the sense of~\cite{wrightCyl}). In this situation, the intersection
with the locus where the Prym variety (as defined in detail
in Section~\ref{sec:complementary}) has real multiplication by a
quadratic field is a union of \Teichmuller curves. That is, if we let
\bas
\Omega G_D \= \{(X,\omega) \in \Omega G: \omega\,\,\, \text{is an eigenform for 
real multiplication}& \\
\text{by $\cO_D$ on $\Prym X$} &\}
\eas
the image $G_D \subset \moduli[4]$ is a finite union of Teich\-m\"uller
curves by  \cite[Theorem~1.7]{MMW}. We give flat pictures of
some of these \Teichmuller curves. 
\par
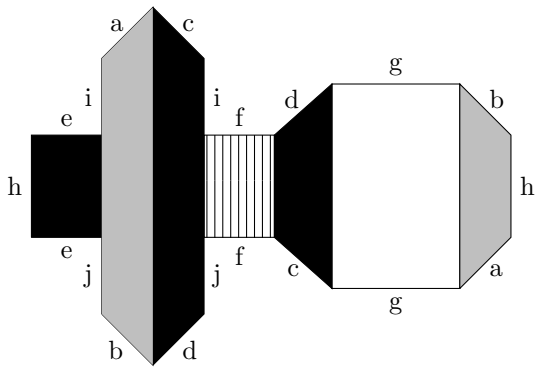
\begin{figure}[htb]
\centering
\input{pic_cathedral}
\caption{A cathedral-shaped surface in the Gothic locus}
\label{cap:Gothicpres}
\end{figure}

The first flat picture is the Gothic cathedral Figure~\ref{cap:Gothicpres}.
It was obtained in \cite{MMW} by shearing jointly the light gray cylinders
in Figure~\ref{cap:Gothichex} (which preserves membership in $\Omega G$)
and a cut and paste operation until these light gray cylinders become the ones
containing the sides $e$, $f$ and $g$ and the dark gray cylinders
transformed into the cylinders containing the sides~$i$ and $j$. After normalizing in the horizontal and vertical directions, one can furthermore assume that the periods have the following form
	\[
	\mathrm{a}=\mathrm{d}=\frac{1+i}{2} \,,\quad \mathrm{b}=\mathrm{c}=\frac{1-i}{2}\,,\quad \mathrm{h}=i\,, \quad
	\mathrm{e} = \mathrm{f} = \mathrm{g}/2 = \alpha\,,\quad \mathrm{i} = \mathrm{j} = \beta i\,,\]
for some $\alpha,\beta\in\bR$. For appropriate values of $\alpha$ and $\beta$
the ratios of the moduli of all vertical and of all horizontal cylinders are commensurable and, therefore, the Veech group of the cathedral contains parabolic elements fixing the vertical and horizontal direction. In fact, this happens whenever $\alpha = x + y \sqrt{d}$ and $\beta = -3x - 3/2 + 3y \sqrt{d}$ for $d>0$ and $x,y\in\bQ$. The product of such parabolic elements is then hyperbolic and has quadratic trace field $\bQ(\sqrt{d})$ and, consequently, Figure~\ref{cap:Gothicpres} generates a \Teichmuller
curve (\cite[Section~9]{MMW}). 
\par
A more precise computation shows that
e.g.\ for $x=0$, $y=1/2$ and $d=2$ the period matrix of $\Prym(X,\pi_A,\pi_B)$
(see Section~\ref{sec:complementary}) is equivalent to
	\[\Pi= \bigg(\begin{array}{rrrr}
-3 \sqrt{2} &-3 \sqrt{2} -1 + 3 \i & -3  + 3 \i & 3 \i \sqrt{2} - 3 \i \\
 3 \sqrt{2} & 3 \sqrt{2} -1 + 3 \i & -3 + 3 \i  &-3 \i \sqrt{2} - 3 \i
\end{array} \bigg)\,.\]
This abelian variety
admits real multiplication by $\cO_{288}$ as can by seen by the
analytic and rational representations 
	\[A_{\frac{\sqrt{288}}{2}}=\begin{pmatrix} \frac{\sqrt{288}}{2} & 0 \\ 0 & -\frac{\sqrt{288}}{2}\end{pmatrix}\quad\mbox{and}\quad
	R_{\frac{\sqrt{288}}{2}} =  \begin{pmatrix}  18&11&-3&-9\\ -18&-9&9&9\\ 18&15&-3&-3\\ 0&6&6&-6 \end{pmatrix}\,,\]
        i.e.\ the identity $A_{\frac{\sqrt{288}}{2}} \Pi = \Pi R_{\frac{\sqrt{288}}{2}}$
        holds.
\par
\medskip
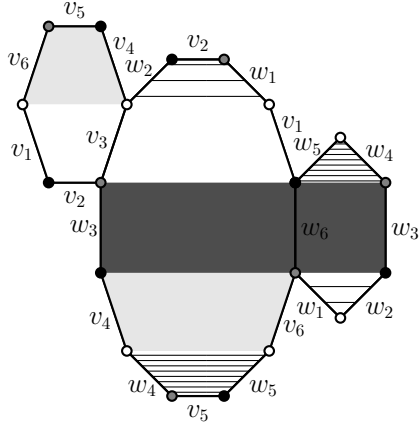
\begin{figure}[htb]
\centering
\input{ente_horizontal_schraffiert}
\caption{A duck-shaped surface in the Gothic locus}
\label{cap:ente}
\end{figure}
Alternatively, one can move the corners of the hexagon while maintaining
the relations~\eqref{eq:lineareq} and the surface becomes horizontally and vertically
periodic with cylinders as in Figure~\ref{cap:ente}. Concretely, we may take
\ba \label{eq:entecoords}
v_1 &= x+yi\,, \quad  &v_2 &= 2x\,, \quad &v_3 &= x-yi\,,  \\
w_1 &= 1-i\,,\quad   &w_2 &= 1+i\,, \quad &w_3 &= 2i
\ea 
for $x,y \in \bR$.
\par
\begin{prop} Let $x,y\in\mathbb{Q}(\sqrt{d})$ so that 
	\[\frac{(1+3x)}{y(1+x)}\,,\ \frac{x(y+3)}{(y+1)} \in \mathbb{Q}\,.\]
Then the flat surface in Figure~\ref{cap:ente} generates a Teichm\"{u}ller curve
in $G_{D}$ for some $D$ such that $\mathbb{Q}(\sqrt{D})=\mathbb{Q}(\sqrt{d})$.
\end{prop}
\par
\begin{proof} The coordinates in~\eqref{eq:entecoords} are chosen so that 
the flat surface  admits a horizontal and a vertical cylinder decomposition.
The  moduli  of the horizontal  cylinders are given
by $(m_{1},m_{2}) = ( 2(1+3x)/y , 2(1+x) )$ while the moduli of vertical
cylinders are $(m_{1}',m_{2}') = (2(y+1)/x , 2(y+3))$, the commensurability of
which is given by the above conditions on~$x$ and~$y$.
\end{proof}
\par
It is amusing to note that a curve in this series of Teichm\"{u}ller curves
in $G_D$ was in the literature long before the discovery of the whole series.
The (irreducible) curve $G_{12}$ will be our second running example. We recall the notation $\mathcal{T}(m,n)$ of Wright (\cite{Wri}) for the Veech-Ward-Bouw-M\"{o}ller curve generated by the
unfolding of the $(m,n,\infty)$-triangle (see also~\cite{ward}), and the `semiregular polygons' decomposition of the corresponding Veech surface $(Y_{m,n},\eta_{m,n})$ of Hooper (\cite{Hoo}).
\par
\begin{prop} \label{prop:G12summary}
The Teichm\"{u}ller curve $G_{12}$ agrees with the Veech-Ward-Bouw-M\"{o}ller curve
$\mathcal{T}(3,6)$. It is generated by the flat surface in Figure~\ref{cap:ente} with $x=\frac{\sqrt{3}}{3}$, $y=-\sqrt{3}$, which agrees with the semiregular polygon decomposition of $(Y_{3,6},\eta_{3,6})$ after scaling the axes by $3/4$ and $4/\sqrt{3}$. The Veech group of $G_{12}$ is the triangle group $\Delta(3,6,\infty)$, hence $\chi(G_{12}) =  -1/2$.
\end{prop}
\par
\begin{proof} The equivalence of the flat presentation is a straightforward 
check using the notation conventions given in the references. To see that
this example corresponds to discriminant $D=12$ in the Gothic series it is enough to check that
\[\Pi= \left(\begin{array}{rrrr}
\frac{(18\i-6)(\sqrt{3}+1)}{\sqrt{3}+3} & \frac{12\i (\sqrt{3}+1)}{\sqrt{3}+3} & 6\i-6 & 4 \\
\frac{(18\i-6)(\sqrt{3}-1)}{\sqrt{3}-3} & \frac{12\i (\sqrt{3}-1)}{\sqrt{3}-3} & 6\i-6 & 4
\end{array} \right)\]
gives the period matrix of the corresponding Prym variety $\Prym(X,\pi_{A},\pi_{B})$ (see Section~\ref{sec:complementary}) and it admits real multiplication by $\cO_{12}$ defined by the analytic and rational representation 
	\[A_{\frac{\sqrt{12}}{2}}=\begin{pmatrix} \frac{\sqrt{12}}{2} & 0 \\ 0 & -\frac{\sqrt{12}}{2}\end{pmatrix}\quad\mbox{and}\quad
	R_{\frac{\sqrt{12}}{2}} =  \begin{pmatrix}   0&0&3&-2\\  0&0&-3&3\\  3&2&0&0\\  3&3&0&0 \end{pmatrix}.\]

\end{proof}

%% file: cylinder_preperation_rotated.tex
\usetikzlibrary{calc}

\definecolor{lgray}{RGB}{214, 214, 214}
\definecolor{dgray}{RGB}{115, 115, 115}

\begin{tikzpicture} 
\begin{scope}[rotate=-60]
\newdimen\R
\R=2cm 
\newdimen\r
\r=1cm 
\coordinate (A) at (0,0);
\path[draw] (A) --++(0:\R) coordinate (B) --++(60:\R) coordinate (C) --++(120:\R) coordinate (D) --++(180:\R) coordinate (E) --++(240:\R) coordinate (F) --++(300:\R) -- cycle;
%
\path[draw] (F) --++ (180:\r) coordinate (G) --++(240:\r) coordinate (H) --++(300:\r) coordinate (I) %
--++(360:\r) coordinate (J) --++(60:\r) coordinate (K) -- cycle;
%
\path[draw] (D) --++ (60:\r) coordinate (L) --++(0:\r) coordinate (M) --++(300:\r) coordinate (N) 
--++(240:\r) coordinate (O) --++(180:\r) coordinate (P) -- cycle;
\coordinate (AB) at ($(A)!0.5!(B)$);
\coordinate (BC) at ($(B)!0.5!(C)$);
\coordinate (CD) at ($(C)!0.5!(D)$);
\coordinate (DE) at ($(D)!0.5!(E)$);
\coordinate (EF) at ($(E)!0.5!(F)$);
\coordinate (FA) at ($(F)!0.5!(A)$);
\coordinate (GH) at ($(G)!0.5!(H)$);
\coordinate (JK) at ($(J)!0.5!(K)$);
\coordinate (EEF) at ($(E)!0.5!(EF)$);
\coordinate (BBC) at ($(B)!0.5!(BC)$);
\coordinate (Q) at (-2.17,.3);
\coordinate (R) at (2.83,2.025);
\coordinate (S) at (3.63,2.82);
\coordinate (T) at (-.87,1.95);
\coordinate (U) at (-.2,.35);
\coordinate (V) at (2.87,1.5);
\coordinate (W) at (3.8,3.82);
\coordinate (X) at (2.37,4.1);

\path[draw=black,fill=dgray] (EF) -- (T) -- (S) -- (N) -- cycle;
\path[draw=black,fill=dgray] (U) -- (BC) -- (V) -- (FA) -- cycle;
\path[draw=black,fill=dgray] (D) -- (W) -- (M) -- (X) -- cycle;
\path[draw=black,fill=lgray] (A) -- (AB) -- (BBC) -- (BC) -- cycle;
\path[fill=lgray] (I) -- (J) -- (JK) -- cycle;
\path[draw=black,fill=lgray] (F) -- (G) -- (GH) -- cycle;
\path[draw=black,fill=lgray] (H) -- (Q) -- (R) -- (CD) -- cycle;
\path[fill=lgray] (D) -- (DE) -- (EEF) -- (EF) -- cycle;


\draw (D) -- (EF);

\draw[thick,dotted] (FA) -- (H) (FA) -- (C) (F) -- (CD);
\draw[thick,dotted,color=white] (CD) -- (N) (D) -- (M); 
\draw[dashed] (DE) -- (EEF) (I) -- (JK);

%
\coordinate (A) at (0,0);
\path[draw] (A) --++(0:\R) coordinate (B) --++(60:\R) coordinate (C) --++(120:\R) coordinate (D) --++(180:\R) coordinate (E) --++(240:\R) coordinate (F) --++(300:\R) -- cycle;
%
\path[draw] (F) --++ (180:\r) coordinate (G) --++(240:\r) coordinate (H) --++(300:\r) coordinate (I) %
--++(360:\r) coordinate (J) --++(60:\r) coordinate (K) -- cycle;
%
\path[draw] (D) --++ (60:\r) coordinate (L) --++(0:\r) coordinate (M) --++(300:\r) coordinate (N) 
--++(240:\r) coordinate (O) --++(180:\r) coordinate (P) -- cycle;
%
\foreach \x in {B,CD,E,FA,H,M}{
	\filldraw[thick,draw=black,fill=gray] (\x) circle (3pt);
}
\foreach \x in {A,BC,D,EF,G,J,N}{
	\filldraw[thick,draw=black,fill=black] (\x) circle (3pt);
}
\foreach \x in {AB,C,DE,F,I,L,O}{
	\filldraw[thick,draw=black,fill=white] (\x) circle (3pt);
}

%
%

\foreach \from/\to/\name in {J/K/$v_{2}$, B/BC/$v_{5}$, BC/C/$w_{5}$, O/N/$w_{2}$
}
\path (\from) -- (\to) node[midway, auto=right] {\name};

\foreach \from/\to/\name in {C/P/$v_{6}$, N/M/$w_{3}$, A/K/$w_{3}$, P/D/$w_{6}$}
\path (\from) -- (\to) node[midway, auto=right, xshift=-3pt,yshift=2pt] {\name};

\foreach \from/\to/\name in {D/DE/$v_{1}$, DE/E/$w_{1}$, F/G/$v_{4}$, M/L/$w_{4}$}
\path (\from) -- (\to) node[midway, auto=right, xshift=-2pt,yshift=-5pt] {\name};

\foreach \from/\to/\name in {P/O/$w_{1}$}
\path (\from) -- (\to) node[midway, auto=left, xshift=-3pt,yshift=-5pt] {\name};

\foreach \from/\to/\name in {I/H/$v_{6}$, FA/F/$v_{3}$}
\path (\from) -- (\to) node[midway, auto=left, xshift=2pt,yshift=-4pt] {\name};

\foreach \from/\to/\name in {J/I/$v_{1}$, AB/A/$v_{4}$, B/AB/$w_{4}$}
\path (\from) -- (\to) node[midway, auto=left, xshift=4pt,yshift=2pt] {\name};

\foreach \from/\to/\name in {D/L/$w_{5}$, EF/E/$v_{2}$, F/EF/$w_{2}$, H/G/$v_{5}$
}
\path (\from) -- (\to) node[midway, auto=left] {\name};

\end{scope}
\end{tikzpicture}

%% file: pic_cathedral.tex

{\scalefont{1}

\begin{tikzpicture} [scale=1.7],transform shape]


 \filldraw[black] (0,0) -- (0,.8) node[midway,left] {h} -- (.55,.8) node[midway,above] {e} -- (.55,0) --(0,0) node[midway,below] {e};

 \draw[] (.55,-.6) -- (.55,1.4) -- (.95,1.8) node[midway,auto=left,xshift=2pt,yshift=-3pt] {a} -- (.95,-1) -- (.55,-.6) node[midway,auto=left,xshift=2pt,yshift=3pt] {b};
 \fill[lightgray] (.55,-.6) -- (.55,1.4) -- (.95,1.8) -- (.95,-1) -- (.55,-.6);
 \draw (.95,-1) -- (1.35,-.6) node[midway,auto=right,xshift=-2pt,yshift=3pt] {d} -- (1.35,1.4) -- (.95,1.8) node[midway, auto=right,xshift=-2pt,yshift=-3pt] {c};
 \fill[black] (.95,-1) -- (1.35,-.6) -- (1.35,1.4) -- (.95,1.8);
 \draw (1.35,0) -- (1.9,0) node[midway,below] {f} -- (1.9,.8) -- (1.35,.8) node[midway,above] {f}  --(1.35,0);
 \fill[pattern=vertical lines] (1.35,0) -- (1.9,0) -- (1.9,.8) -- (1.35,.8) --(1.35,0);
 \filldraw[black] (1.9,0) -- (2.35,-.4) node[midway,auto=right,xshift=2pt,yshift=3pt] {c} -- (2.35,1.2)  -- (1.9,.8) node[midway,auto=right,xshift=2pt,yshift=-3pt] {d} -- (1.9,0); 
 \draw (2.35,-.4) -- (3.35,-.4) node[midway,below] {g} -- (3.35,1.2) -- (2.35,1.2) node[midway,above] {g} -- (2.35,-.4);  
 \fill[lightgray] (3.35,-.4) -- (3.75,0) -- (3.75,.8) -- (3.35,1.2) -- (3.35,-.4);
 \draw (3.35,-.4) -- 
       (3.75,0) node[midway,auto=right,xshift=-2pt,yshift=3pt] {a} -- 
       (3.75,.8)  node[midway,right] {h} -- 
       (3.35,1.2) node[midway,auto=right,xshift=-2pt,yshift=-3pt] {b} -- 
       (3.35,-.4); 
 
 \coordinate (i1) at ($(.55,.8)!0.5!(.55,1.4)$); 
 \coordinate (i2) at ($(1.35,.8)!0.5!(1.35,1.4)$);     
 \coordinate (j1) at ($(.55,-.6)!0.5!(.55,0)$); 
 \coordinate (j2) at ($(1.35,-.6)!0.5!(1.35,0)$);    
 \draw[] (i1) node[xshift=-5pt] {i} ;
 \draw[] (i2) node[xshift=5pt] {i} ; 
 \draw[] (j1) node[xshift=-5pt] {j} ;
 \draw[] (j2) node[xshift=5pt] {j} ;

\end{tikzpicture}
}

%% file: ente_horizontal_schraffiert.tex
{\scalefont{1.6}
\begin{tikzpicture}[scale=.6, transform shape]

\draw[thick] 
({\x*(0)}, {(0)*\y})  coordinate (1) -- ({\x*(1/2*sqrt(3))}, {(0)*\y})  coordinate (2) -- ({\x*(1/2*sqrt(3) + 3/4)}, {(1/4*sqrt(3))*\y})  coordinate (3) -- ({\x*(3/4*sqrt(3) + 3/4)}, {(1/4*sqrt(3) + 3/4)*\y})  coordinate (4) -- ({\x*(3/4*sqrt(3) + 3/4)}, {(3/4*sqrt(3) + 3/4)*\y})  coordinate (5) -- ({\x*(1/2*sqrt(3) + 3/4)}, {(3/4*sqrt(3) + 3/2)*\y})  coordinate (6) -- ({\x*(1/2*sqrt(3))}, {(sqrt(3) + 3/2)*\y})  coordinate (7) -- ({\x*(0)}, {(sqrt(3) + 3/2)*\y})  coordinate (8) -- ({\x*((-3)/4)}, {(3/4*sqrt(3) + 3/2)*\y})  coordinate (9) -- ({\x*((-1)/4*sqrt(3) - 3/4)}, {(3/4*sqrt(3) + 3/4)*\y})  coordinate (10) -- ({\x*((-1)/4*sqrt(3) - 3/4)}, {(1/4*sqrt(3) + 3/4)*\y})  coordinate (11) -- ({\x*((-3)/4)}, {(1/4*sqrt(3))*\y})  coordinate (12) -- cycle;

\begin{scope}[shift={\ShiftP}]
\draw[thick]  
({\x*(0)}, {(0)*\y})  coordinate (13) -- ({\x*(3/4)}, {(1/4*sqrt(3))*\y})  coordinate (14) -- ({\x*(3/4)}, {(3/4*sqrt(3))*\y})  coordinate (15) -- ({\x*(0)}, {(sqrt(3))*\y})  coordinate (16) -- ({\x*((-3)/4)}, {(3/4*sqrt(3))*\y})  coordinate (17) -- ({\x*((-3)/4)}, {(1/4*sqrt(3))*\y})-- cycle;
\end{scope}

\begin{scope}[shift={\ShiftQ}]
\draw[thick]  
({\x*(0)}, {(0)*\y})  coordinate (18) -- ({\x*(1/2*sqrt(3))}, {(0)*\y})  coordinate (19)  -- ({\x*(3/4*sqrt(3))}, {(3/4)*\y})  coordinate (20) -- ({\x*(1/2*sqrt(3))}, {(3/2)*\y})  coordinate (21) -- ({\x*(0)}, {(3/2)*\y})  coordinate (22) -- ({\x*((-1)/4*sqrt(3))}, {(3/4)*\y})  coordinate (23) -- ({\x*(0)}, {(0)*\y}) -- cycle;
\end{scope}

\foreach \from/\to/\name in {23/18/$v_{1}$, 18/19/$v_{2}$, 19/20/$v_{3}$, 20/21/$v_{4}$, 21/22/$v_{5}$, 22/23/$v_{6}$, 4/13/$w_{1}$, 13/14/$w_{2}$, 14/15/$w_{3}$, 15/16/$w_{4}$, 16/17/$w_{5}$, 17/4/$w_{6}$, 2/1/$v_{5}$, 3/2/$w_{5}$, 4/3/$v_{6}$, 6/5/$v_{1}$, 7/6/$w_{1}$, 8/7/$v_{2}$, 9/8/$w_{2}$, 11/10/$w_{3}$, 12/11/$v_{4}$, 1/12/$w_{4}$ }
;

\fill[pattern=closely horizontal lines] (1) -- (2) -- (3) -- (12) -- (1) (17) -- (15) -- (16) -- (17);
\fill[color=white!80!gray] (23) -- (20) -- (21) -- (22) -- (23) (12) -- (3) -- (4) -- (11) -- (12);
\fill[pattern=wide horizontal lines] (9) -- (6) -- (7) -- (8) -- (9) (4) -- (13) -- (14) -- (4);
\fill[color=white] (18) -- (19) -- (20) -- (23) -- (18) (19) -- (5) -- (6) -- (20) -- (19);
\fill[color=black!40!gray] (11) -- (4) -- (5) -- (10) -- (11) (4) -- (14) -- (15) -- (5) -- (4);


\draw[thick] 
({\x*(0)}, {(0)*\y})  coordinate (1) -- ({\x*(1/2*sqrt(3))}, {(0)*\y})  coordinate (2) -- ({\x*(1/2*sqrt(3) + 3/4)}, {(1/4*sqrt(3))*\y})  coordinate (3) -- ({\x*(3/4*sqrt(3) + 3/4)}, {(1/4*sqrt(3) + 3/4)*\y})  coordinate (4) -- ({\x*(3/4*sqrt(3) + 3/4)}, {(3/4*sqrt(3) + 3/4)*\y})  coordinate (5) -- ({\x*(1/2*sqrt(3) + 3/4)}, {(3/4*sqrt(3) + 3/2)*\y})  coordinate (6) -- ({\x*(1/2*sqrt(3))}, {(sqrt(3) + 3/2)*\y})  coordinate (7) -- ({\x*(0)}, {(sqrt(3) + 3/2)*\y})  coordinate (8) -- ({\x*((-3)/4)}, {(3/4*sqrt(3) + 3/2)*\y})  coordinate (9) -- ({\x*((-1)/4*sqrt(3) - 3/4)}, {(3/4*sqrt(3) + 3/4)*\y})  coordinate (10) -- ({\x*((-1)/4*sqrt(3) - 3/4)}, {(1/4*sqrt(3) + 3/4)*\y})  coordinate (11) -- ({\x*((-3)/4)}, {(1/4*sqrt(3))*\y})  coordinate (12) -- cycle;

\begin{scope}[shift={\ShiftP}]
\draw[thick]  
({\x*(0)}, {(0)*\y})  coordinate (13) -- ({\x*(3/4)}, {(1/4*sqrt(3))*\y})  coordinate (14) -- ({\x*(3/4)}, {(3/4*sqrt(3))*\y})  coordinate (15) -- ({\x*(0)}, {(sqrt(3))*\y})  coordinate (16) -- ({\x*((-3)/4)}, {(3/4*sqrt(3))*\y})  coordinate (17) -- ({\x*((-3)/4)}, {(1/4*sqrt(3))*\y})-- cycle;
\end{scope}

\begin{scope}[shift={\ShiftQ}]
\draw[thick]  
({\x*(0)}, {(0)*\y})  coordinate (18) -- ({\x*(1/2*sqrt(3))}, {(0)*\y})  coordinate (19)  -- ({\x*(3/4*sqrt(3))}, {(3/4)*\y})  coordinate (20) -- ({\x*(1/2*sqrt(3))}, {(3/2)*\y})  coordinate (21) -- ({\x*(0)}, {(3/2)*\y})  coordinate (22) -- ({\x*((-1)/4*sqrt(3))}, {(3/4)*\y})  coordinate (23) -- ({\x*(0)}, {(0)*\y}) -- cycle;
\end{scope} 

\foreach \from/\to/\name in {2/1/$v_{5}$, 8/7/$v_{2}$}
\path (\from) -- (\to) node[midway, auto=left] {\name}
;

\foreach \from/\to/\name in {11/10/$w_{3}$}\path (\from) -- (\to) node[midway, auto=left,xshift=4pt] {\name};
\foreach \from/\to/\name in {12/11/$v_{4}$, 18/23/$v_{1}$}\path (\from) -- (\to) node[midway, auto=left,xshift=4pt,yshift=5pt] {\name};
\foreach \from/\to/\name in {3/2/$w_{5}$, 14/13/$w_{2}$}\path (\from) -- (\to) node[midway, auto=left,xshift=-4pt,yshift=1pt] {\name};
\foreach \from/\to/\name in {4/3/$v_{6}$}\path (\from) -- (\to) node[midway, auto=left,xshift=-5pt,yshift=0pt] {\name};
\foreach \from/\to/\name in {12/1/$w_{4}$, 4/13/$w_{1}$}\path (\from) -- (\to) node[midway, auto=right,xshift=10pt,yshift=1pt] {\name};
\foreach \from/\to/\name in {14/15/$w_{3}$, 4/17/$w_{6}$}\path (\from) -- (\to) node[midway, auto=right,xshift=-2pt,yshift=-2pt] {\name};
\foreach \from/\to/\name in {15/16/$w_{4}$, 6/7/$w_{1}$}\path (\from) -- (\to) node[midway, auto=right,xshift=-4pt,yshift=-6pt] {\name};
\foreach \from/\to/\name in {17/16/$w_{5}$, 9/8/$w_{2}$}\path (\from) -- (\to) node[midway, auto=left,xshift=8pt,yshift=-2pt] {\name};
\foreach \from/\to/\name in {5/6/$v_{1}$, 20/21/$v_{4}$}\path (\from) -- (\to) node[midway, auto=right,xshift=-6pt,yshift=2pt] {\name};
\foreach \from/\to/\name in {22/23/$v_{6}$,20/19/$v_{3}$}\path (\from) -- (\to) node[midway, auto=right,xshift=2pt,yshift=-8pt] {\name};
\foreach \from/\to/\name in 
{18/19/$v_{2}$, 21/22/$v_{5}$}
\path (\from) -- (\to) node[midway, auto=right] {\name}
;

\foreach \x in {1,4,7,10,15,22}
{
  \filldraw[thick,draw=black,fill=gray] (\x) circle (3pt); 
}

\foreach \x in {2,5,8,11,14,17,18,21}
{
  \filldraw[thick,draw=black,fill=black] (\x) circle (3pt); 
}

\foreach \x in {3,6,9,12,13,16,23}
{
  \filldraw[thick,draw=black,fill=white] (\x) circle (3pt); 
}
\end{tikzpicture}
}

%% file: sec_complementary.tex
\section{Prym varieties for two maps} 
\label{sec:complementary}

Given a finite collection of maps $\pi_i: X \to Y_i$ between curves,
the Prym variety $\Prym(X,\pi_1,\ldots,\pi_n)$ (in a generalized sense)
is the complementary abelian variety to the image of the maps $\pi_i^{*}:\Jac Y_i
\to \Jac X$, that is the perpendicular space to the tangent spaces
$\oplus \Omega_{Y_i}^\vee$ divided by its intersection with the period lattice.
The main goal of this section is to determine the signature of the
polarization on this Prym variety $\Prym(X) = \Prym(X,\pi_A,\pi_B)$ in the
case of a Gothic flat surface.
\par
\begin{prop}  \label{prop:Prymis16}
The restriction of the principal polarization on $\Jac(X)$
is a polarization on $\Prym X$ of type~$(1,6)$. Consequently,
the dual Prym variety $\Prym^\vee X$ has a natural polarization
of type~$(1,6)$, too.
\end{prop}
\par
We first recall some equivalent definition of complementary abelian subvarieties
in terms of endomorphisms. Let $(T,\cL)$ be an abelian variety, that is a complex
torus $T=V/\Lambda$ together with a positive definite line bundle~$\cL$. Given
an abelian subvariety $\iota: Y\to T$, one can define its exponent $e_{Y}$ as the
exponent $e(\cL)$ of the induced polarization $\iota^{*}\cL$, and its norm endomorphism
$N_{Y}\in\End(T)$ and symmetric idempotent $\varepsilon_{Y}\in\End_{\bQ}(T)$ as
\be \label{eq:NYdef}
N_{Y} \coloneqq \iota\,\psi_{\iota^{*}\cL}\,\check{\iota}\,\phi_{\cL}\qquad
\mbox{and}\qquad\varepsilon_{Y} \coloneqq \frac{1}{e_{Y}}N_{Y}\,,
\ee
respectively, where $\phi_{\cL}:T\to T^{\vee}$ is the isogeny associated to a line bundle~$\cL$ and $\psi_{\cL}=e(\cL)\phi^{-1}_{\cL}$ (see~\cite[Section 5.3]{birkenhakelange}). In the case of $\cL$ being a principal polarization, the
exponent of $Y$ is precisely $e_{Y}=\min \{n>0\,:\, n\varepsilon_{Y}\in\End(T)\}$ (see~\cite[Prop. 12.1.1]{birkenhakelange}).
\par
The assignment $Y\mapsto \varepsilon_{Y}$ and its inverse
$\varepsilon\mapsto X^{\varepsilon}\coloneqq\im(n\epsilon)$,
for any $n>0$ such that $n\varepsilon\in\End(T)$, induce a bijection between
the set of abelian subvarieties of~$T$ and the set of symmetric (with respect
to the Rosati involution $f \mapsto f' = \phi^{-1}_{\cL} \widehat{f} \phi_\cL$)
idempotents in $\End_{\bQ}(T)$. Accordingly, the canonical involution
$\varepsilon \mapsto 1-\varepsilon$ on the set of symmetric idempotents induces
an involution $Y\mapsto Z\coloneqq X^{1-\varepsilon_{Y}}$ on the set of abelian
subvarieties. The abelian subvariety $Z$ is called the \emph{complementary
abelian subvariety} of $Y$, and the exponent $e_{Z}$ agrees with $e_{Y}$ in
the case of $\cL$ being a principal polarization.
The map $(N_{Y},N_{Z}): X \to Y\times Z$ 
is an isogeny and the following identities 
	\bas
N_{Y}|_{Y} \= e_{Y}\Id \,, \quad N_{Y}|_{Z} \= 0 \= N_{Y}N_{Z} &\= 0 \,,
\quad e_{Y}N_{Z} + e_{Z}N_{Y} &\= e_{Y}e_{Z}\Id\,.
	\eas
hold (\cite[Section~5.3]{birkenhakelange}).
\par
Let now $\pi:X\to Y$ be a morphism between curves. The pullback map defines a
homomorphism $\pi^{*}:\Jac Y \to \Jac X$. This map is, moreover,  injective
whenever $\pi$ does not factor through a cyclic \'{e}tale cover of degree $\ge 2$.
Under these conditions, the Prym variety $\Prym(X,\pi)$ of the map~$\pi$ is defined
as the complementary abelian variety of $\pi^{*}(\Jac Y)$ (or, equivalently, as the
connected component of the identity of the kernel $\ker N_{\pi^{*}(\Jac Y)}$). The
Jacobian of $X$ decomposes, up to isogeny, as $\Jac X \sim \pi^{*}(\Jac Y)
\times \Prym(X,\pi)$. Note that, in general, $\Prym(X,\pi)$ is not a Prym variety
in the classical sense (see~\cite[Section 12]{birkenhakelange}), as the induced
polarization will not be a multiple of the principal polarization.
\par
\medskip
Consider now a pair of morphisms of curves $\pi_{1}:X\to Y_{1}$ and $\pi_{2}:X\to Y_{2}$, together with the corresponding homomorphisms $\pi_{1}^{*}:\Jac Y_{1} \to \Jac X$ and $\pi_{2}^{*}:\Jac Y_{2} \to \Jac X$. Assume moreover that there exist morphisms $g_{1}:Y_{1}\to Y$ and $g_{2}:Y_{2}\to Y$ to some curve $Y$ such that the diagram
\[\begin{tikzpicture}
\matrix (m) [matrix of math nodes, row sep=1.5em, column sep=3.5em, text height=1.5ex, text depth=0.25ex]
{ & X & \\
Y_{1} & & Y_{2} \\
& Y & \\
};
\path[->,font=\scriptsize]
(m-1-2) edge [above] node {$\pi_{1}$} (m-2-1)
(m-1-2) edge [above] node {$\pi_{2}$} (m-2-3)
(m-2-1) edge [above] node {$g_{1}$} (m-3-2)
(m-2-3) edge [above] node {$g_{2}$} (m-3-2)
(m-1-2) edge [right] node {$h$} (m-3-2)
;
\end{tikzpicture}\]
commutes. Under a mild non-factorization condition, one can decompose $\Jac X$
further in terms of Jacobians.
\par
\begin{prop}[\cite{LaRe}]\label{prop:PrymDef} Suppose $g_{1}$ and $g_{2}$ do not
both factorize via the same morphism $Y_{0} \to Y$ of degree $\ge 2$. Then
$\pi_{2}^{*}\Prym(Y_2,g_{2})$ is an abelian subvariety of $\Prym(X,\pi_{1})$. In
particular, $\Jac X$ decomposes, up to isogeny, as
	\[
\Jac X \sim h^{*}(\Jac Y) \times \pi_{1}^{*}\Prym(Y_1,g_{1})
\times \pi_{2}^{*} \Prym(Y_2, g_{2}) \times P\,,
	\]
for some subvariety $P$ of $\Jac X$.
\end{prop}
\par
The subvariety $P$ is called the \emph{Prym variety} $\Prym(X,\pi_{1},\pi_{2})$
of the pair of coverings $(\pi_{1},\pi_{2})$. In the case that $Y=\bP^{1}$
the summand $h^{*}(\Jac Y)$ is of course trivial and $\Prym(Y_j,g_{j})=\Jac Y_{j}$.
\par
We now specialize to the Gothic situation and give explicitly the
various norm endomorphisms for later use. We write $A^\vee$ and $B^\vee$
for the image of $\pi_A^*$ and $\pi_B^*$ respectively.
\par
\begin{prop}\label{prop:NYNP}Let $T$ be a principally polarized abelian variety and $A^\vee,B^\vee
\subset T$  be abelian subvarieties with coprime exponents $e_{\A},e_{\B}$ and such that
$N_{\A}N_{\B}=0$. Then $Y=A^\vee \times \B^\vee$ is a subvariety of~$T$. Moreover, the norm
endomorphisms of~$Y$ and its complementary abelian variety~$P$ satisfy 
	\bes
	N_{Y} \=e_{\B}N_{\A}+e_{\A}N_{\B}\,,\quad\mbox{and}\quad N_{P}\= e_{\A}e_{\B}\Id-N_{Y}\,.
	\ees
\end{prop}
\par
\begin{proof} The injectivity of $Y \to T$ follows from coprimality.
Writing $N=e_{\B}N_{\A}+e_{\A}N_{\B}$, one has $N^{2}= e_{\A} e_{\B}N$ and
$N|_{Y} = e_{\A} e_{\B}\Id_{Y}$. The idempotent $\varepsilon=\frac{1}{e_{\A} e_{\B}}N$
corresponds to the abelian subvariety $Y$ and it is of exponent $e_{\A} e_{\B}$
since $e_{Y}=\min\{n>0\,:\, n\varepsilon_{Y}\in\End(T)\}$ and $(e_{\A},e_{\B})=1$.
The rest of the claims follow.
\end{proof}
\par
\par
\begin{proof}[Proof of Proposition~\ref{prop:Prymis16}]
Thanks to diagram~\ref{eq:Gothic} and since $\gcd(2,3)=1$ the hypothesis of the
Proposition~\ref{prop:PrymDef} is met. Moreover, $\A^\vee \times \B^\vee$ has a polarization of type $(1,6)$ and, by \cite[Corollary~12.1.5]{birkenhakelange} the same holds for the complementary abelian variety.
\end{proof}
\par

%% file: sec_HMVAR.tex
\section{Hilbert modular surfaces and modular embeddings}
\label{sec:HMS}

The Prym-Torelli map~$t$ associates with a flat surfaces
in the Gothic locus, or more generally with any genus four
surface admitting maps $\pi_A$ and $\pi_B$ that fit into
the diagram~\eqref{eq:Gothic}, the dual Prym variety $\Prym^\vee(X,\pi_A,\pi_B)$.
(The reason for dualizing will become apparent in Section~\ref{sec:GothMF}.)
By Proposition~\ref{prop:Prymis16} this gives a map $t: \Omega G \to
\cA_{2,(1,6)}$ to the  the moduli  space of $(1,6)$-polarized abelian
surfaces. 
The goal of this section is to recall some basic properties
of Hilbert modular surfaces that arise from the following observation.
\par
\begin{prop} \label{prop:tGD}
  The Prym-Torelli-image $t(G_D)$ of the Gothic locus is contained
  in the image of a Hilbert modular surface $X_D(\frakb)$ inside the moduli
  space of $(1,6)$-polarized abelian surfaces, where $\frakb$
  is an $\cO_D$-ideal of norm~$6$.
\end{prop}
\par
We compute here the Euler characteristics of these
Hilbert modular surfaces $X_D(\frakb)$ and discuss the modular
embeddings that induce the map $X_D(\frakb) \to \cA_{2,(1,6)}$.

\subsection{Hilbert modular surfaces.}

For any positive discriminant $D\equiv 0,1\mod 4$, write $D=b^2-4ac$ for some
$a,b,c \in\bZ$. The (unique) \emph{quadratic order of discriminant $D$}
is defined as $\cO_D=\bZ[T]/(aT^2+bT+c)$. This order agrees with
$\cO_D=\bZ\oplus \gamma_D\bZ$ inside the quadratic field $K=\bQ(\sqrt{D})$,
where $\gamma\coloneqq\gamma_{D}=\frac{D+\sqrt{D}}{2}$ provided that $D$
is not a square.
\par
For any fractional ideal $\frakc\subset K$, we denote by $\frakc^{\vee}$ the dual with respect to the trace pairing, i.e.\ $\frakc^{\vee}=\{x\in K: \tr_{\bQ}^{K}(x\frakc)\subset\bZ\}$. In particular, $\cO_{D}^{\vee}=\frac{1}{\sqrt{D}}\cO_{D}$.
\par

Let $\frakb$ be an $\cO_{D}$-ideal. The $\cO_{D}$-module $\frakb\oplus\cO_{D}^{\vee}$ is preserved by the \emph{Hilbert modular group}
	\[\SL(\frakb\oplus\cO_{D}^{\vee}) \=
	\begin{pmatrix} \cO_{D} & \sqrt{D}\,\frakb \\ 
	\frac{1}{\sqrt{D}}\,\frakb^{-1} & \cO_{D} \end{pmatrix}\cap \SL_{2}(K)\,.
	\]
Associated with $\frakb$ we can construct the \emph{Hilbert modular surface}
	\[X_D(\frakb) \= \SL(\frakb\oplus\cO_{D}^{\vee})\backslash \bH^{2}\,.\]

\subsection{Abelian surfaces with real multiplication and a $(1,n)$-polarization.}\label{subsec:abeliansurfaces(1,n)}

An abelian surface $T$ admits \emph{real multiplication} by $\cO_{D}$ if there exists an embedding $\cO_{D}\hookrightarrow \End(T)$ by self-adjoint endomorphisms. We will always assume that the action is \emph{proper}, in the sense that it cannot be extended to an action of a larger quadratic discriminant $\cO_{E}\supset \cO_{D}$.
\par
The different components of the moduli space of $(1,n)$-polarized abelian varieties with a choice of real multiplication by $\cO_D$ are parameterized by certain Hilbert modular surfaces (see~\cite[Chapter 7]{HvdG}). 
\par
More precisely, suppose that $(T = \bC^2/\Lambda,\cL)$ is an abelian variety with a $(1,n)$-polarization $\cL$ and a choice of real multiplication by $\cO_{D}$. Then $\Lambda$ is a rank-two $\cO_D$-module with symplectic pairing of signature $(1,n)$.
By \cite{Bass62} such a lattice splits as a direct sum of $\cO_D$-modules. Moreover, 
although $\cO_D$ is not a Dedekind domain for non-fundamental discriminants $D$, 
any rank-two $\cO_D$-module is isomorphic to $\frakb \oplus \cO_D^\vee$ for some $\cO_D$-ideal $\frakb$. The isomorphism can moreover be chosen so that
the symplectic form is mapped to the trace pairing $\langle (a,b)^{T},(\widetilde{a},\widetilde{b})^{T}\rangle = \tr^{K}_{\bQ}(a\widetilde{b}-\widetilde{a}b)$. The type of such a polarization is $(d_{1},d_{2})$, where $d_i \in \bN$
are uniquely determined by $d_1 | d_2$ and
$\cO_{D}/\frakb\cong\bZ/d_{1}\bZ \times \bZ/d_{2}\bZ$. 
\par
In the case of a polarization of type $(1,n)$, it follows (see for
example~\cite[Prop. 5.2.1]{Cohen}) that the ideal $\frakb$ can be generated
as a $\bZ$-module by $(\frac{r+\sqrt{D}}{2}, n )$ for some $0\le r < 2n$. In particular, $N^K_\bQ(\frakb) = n$.
\par
Conversely, for any ideal $\frakb$ of norm $n$ and
$\bm{\tau}=(\tau_{1},\tau_{2}) \in \bH^{2}$, we define the lattice
\bes
\Lambda_{\frakb,\bm{\tau}} \= \{(a + b\tau_1, a^\sigma+b^\sigma \tau_2)^T \,|\,\, 
a \in \frakb, b \in \cO^\vee \}.
\ees
The quotient $T_{\bm\tau}=\bC^{2}/\Lambda_{\frakb,\bm\tau}$ is an abelian surface 
with a $(1,n)$-polarization (given by the trace pairing) and real multiplication
by $\cO_{D}$. The isomorphism class of $T_{\bm\tau}$ depends only on the image
of~$\bm\tau$ in $X_{D}(\frakb)$.
\par
The {\em proof of Proposition~\ref{prop:tGD}} follows from this observation
and the real multiplication built into the definition of the Gothic curves~$G_D$.
\par
It also follows that the locus of $(1,n)$-polarized abelian varieties with a choice of real
multiplication~by $\cO_D$ has as many components as ideals $\frakb$ of
norm~$n$ in $\cO_{D}$, each of these components being parameterized by
the Hilbert modular surface $X_{D}(\frakb)$. Concretely, for the case
we are interested in:
\par
\begin{prop}\label{prop:Xcomponents}
The moduli space of $(1,6)$-polarized abelian surfaces with a choice of real
multiplication by $\cO_D$ is empty for $D \equiv 5\bmod{8}$ or $D \equiv 2\bmod{3}$. \\ It is non-empty and irreducible for $D \equiv  0,12\bmod{24}$, it has has two irreducible components for $D \equiv 4,9,16\bmod{24}$ and four for $D \equiv 1\bmod{24}$.
\end{prop}
\par
\begin{proof}
By the preceding discussion, the locus of $(1,6)$-polarized abelian varieties
with a choice of real multiplication by $\cO_D$ is non-empty if and only if there is
an $\cO_D$-ideal~$\frakb$ with $N^K_\bQ(\frakb) = 6$, 
i.e. if and only if $D \equiv 0,1,4,9,12,16\bmod{24}$.
\par
Each connected component of this locus is parameterized by a Hilbert modular
surface $X_{D}(\frakb)$ for an $\cO_{D}$-ideal $\frakb$ of norm 6. For
$D \equiv 0,12\mod{8}$ there is exactly one prime ideal of norm two
$\frakb_2$ and one prime ideal $\frakb_3$ of norm three, so that the locus
is connected. For $D \equiv 9\bmod{24}$ the prime two splits (but three
is ramified) and for  $D \equiv 4,16\bmod{24}$ the prime three splits
(but two is ramified), resulting in two connected components. For
$D \equiv 1\bmod{24}$ both primes split.
\end{proof}
\par
Note, however, that the locus of real multiplication in $\cA_{2,(1,6)}$
has in general fewer components than the moduli space of abelian surfaces
with a chosen real multiplication by~$\cO_D$. In fact, the
abelian varieties parameterized by $X_D(\frakb)$ and by $X_D(\frakb^\sigma)$
map to the same subsurface in $\cA_{2,(1,6)}$.

\subsection{Euler characteristics.}
The notion Euler characteristic (of curves and of Hilbert modular surfaces)
refers throughout to orbifold Euler characteristics. Let  $D=f^2 D_0$ be
the factorization of the discriminant into a fundamental discriminant $D_0$
and a square of $f \in \bN$. The Euler characteristic of Hilbert modular
surfaces has been computed by Siegel (\cite{siegel36}), for the more usual
Hilbert modular surface $X_D = X_D(\cO_D)$. A reference including also
the case of non-fundamental
discriminants is \cite[Theorem~2.12]{bainbridge07}. Altogether, 
\bes
\chi(X_D) \= 2f^3 \zeta_{\bQ(\sqrt{D})}(-1) \left(\sum_{r|f} \legendre{D_0}{r}
\frac{\mu(r)}{r^2} \right),
\ees
where $\mu$ is the M\"obius function and $\legendre{a}{b}$ is the Jacobi
symbol. The case we are interested in can be deduced from this formula.
\par
\begin{prop} \label{prop:vol16}
  The Euler characteristics of $X_{D}(\frakb)$, for $\frakb$ of norm~$6$,
  and of $X_D$ are related as follows.
\bes
\kappa_{D}\coloneqq \frac{\chi(X_{D}(\frakb))}{\chi(X_D)} = \left\{\begin{array}{lcl}
1   & \text{if} & {\rm gcd}(6,f)=1 \\
3/2 & \text{if} & {\rm gcd}(6,f)=2 \\
4/3 & \text{if} & {\rm gcd}(6,f)=3 \\
2   & \text{if} & {\rm gcd}(6,f)=6. \\
\end{array} \right.
\ees
\end{prop}
\par
\begin{proof}
The groups $\SL(\cO_D \oplus \cO_D^\vee)$ and $\SL(\frakb \oplus\cO_D^\vee)$
are commensurable. To determine the indices in their intersection, 
we conjugate both groups by $\sm{\sqrt{D}} 0 0 1$. This takes the first
group into $\SL(\cO_D \oplus \cO_D)$ and the second group into
$\SL_2(\frakb \oplus \cO_D)$. The two images under conjugation
contain
\bes
\Gamma_\frakb = \left\{ \mat abcd \in \SL_2(K) \,:\, a,d \in \cO_D, \,
b \in \frakb,\, c \in \cO_D \right\}.
\ees
with a finite index that we now calculate. We factorize
$\frakb = \frakp_2\frakp_3$ into the primes of norm two and three and
consider the action of $\SL(\cO_D \oplus \cO_D)$ on
$\bP^1(\cO_D/\frakp_2) \times \bP^1(\cO_D/\frakp_3)$. This action is transitive,
in fact elementary matrices in $\SL_2(\cO_D/\frakp_2)$ and
$\SL_2(\cO_D/\frakp_3)$ generate a transitive group and elementary
matrices can obviously be lifted. Since $\Gamma_\frakb$ is precisely
the stabilizer of $((0:1),(0:1))$, we conclude
$$[\SL(\cO_D \oplus \cO_D): \Gamma_\frakb] \=
|\bP^1(\cO_D/\frakp_2) \times \bP^1(\cO_D/\frakp_3)| \= 12\,.$$
\par
If $\frakb$ is an invertible ideal, we use
$\SL(\frakb \oplus\cO_D) = \SL(\cO_D \oplus \frakb^{-1} )$
and consider the projection
\bes
\pr: \SL(\cO_D \oplus \frakb^{-1} ) \to \SL(\cO_D/\frakp_3 \oplus
\frakb^{-1}/\frakp_2^{-1}) \times  \SL( \cO_D/\frakp_2 \oplus
\frakb^{-1}/\frakp_3^{-1})\,,
\ees
where in the range the modules are considered as $\cO_D/\frakp_3$-module
and $\cO_D/\frakp_2$-module (i.e.\ as vector spaces), respectively.
Even the smaller group $\Gamma_\frakb$ contains the kernel of~$\pr$, 
and in fact $\Gamma_\frakb$ is precisely the stabilizer of $((0:1),(0:1))$.
We conclude that its index is~$12$ in $\SL(\frakb \oplus\cO_D)$
and this completes the case ${\rm gcd}(6,f)=1$.
\par
If ${\rm gcd}(6,f)=2$, we use that $\frakp_2^{-1} = \cO_{D/4}$ and
consider $\SL(\cO_D \oplus \frakb^{-1} )$   
as a subgroup of
$\SL(\cO_{D/4} \oplus \tilde{\frakp}_3^{-1})$, where the tilde indicates that
we now extended scalars of $\frakp_3^{-1}$ to form an $\cO_{D/4}$-module.
We consider the projection
\bes
\pr: \SL(\cO_{D/4} \oplus \tilde{\frakp}_3^{-1}) \to \SL(\cO_{D/4}/\tilde{\frakp}_3
\oplus \tilde{\frakp}^{-1}_3/\cO_{D/4}) \times \SL( \cO_{D/4}/\frakp_2
\oplus  \tilde{\frakp}^{-1}_3/\frakp_2\tilde{\frakp}^{-1}_3 )\,.
\ees

Again, even the smaller group $\Gamma_\frakb$ contains the kernel of~$\pr$.
The image of $\SL(\cO_D \oplus \frakb^{-1})$ under $\pr$ is contained in
the full first factor times the lower triangular matrices in the second
factor, as can be checked using a set of generators for these groups
consisting of elementary matrices. The image of $\Gamma_\frakb$
under $\pr$ is the stabilizer of $(0:1)$ in the first factor times the lower triangular
matrices with $\cO_D/\frakp_2 \subset \cO_{D/4}/\frakp_2 \cong
\tilde{\frakp}^{-1}_3/\frakp_2\tilde{\frakp}^{-1}_3$ in the lower left corner
in the second factor. This subgroup is of index $4 \cdot 2 = 8$ and
this concludes the case ${\rm gcd}(6,f)=2$.
\par
The remaining cases are similar, using $\frakp_3^{-1} = \cO_{D/9}$ if
${\rm gcd}(6,f)=3$ and $\frakb^{-1} = \cO_{D/36}$ if  ${\rm gcd}(6,f)=6$.
\end{proof}

\subsection{Siegel modular embeddings.} \label{sec:SME}

Let $X_{D}(\frakb)$ parameterize a component of the moduli space of $(1,n)$-polarized abelian varieties with a choice of real multiplication by $\cO_D$ as above. The forgetful map $X_{D}(\frakb)\to \cA_{2,(d_{1},d_{2})}$ to the moduli space of $(d_{1},d_{2})$-polarized abelian varieties can be lifted to a holomorphic map $\psi: \bH^{2} \to \bH_2$ which is equivariant with respect to a homomorphism $\Psi :  \SL(\frakb\oplus\cO_{D}^{\vee}) \to G_{P}$, where $P\coloneqq P_{d_{1},d_{2}}= \left(\begin{smallmatrix}d_{1} & 0 \\ 0 & d_{2}\end{smallmatrix}\right)$ and $G_{P}$ is the symplectic group for the polarization type $(d_{1},d_{2})$ (see~\cite[Section 8.2]{birkenhakelange})
\bes
G_{P}=\left\{M\in \Sp_{4}(\bQ)\,:\, M^{T}
\left(\begin{smallmatrix}I_2 & 0 \\ 0 & P\end{smallmatrix}\right) \bZ^{4}
\subseteq
 \left(\begin{smallmatrix}I_2 & 0 \\ 0 & P\end{smallmatrix}\right)\bZ^{4}\,\right\}\,.
\ees
Such a lift $(\psi,\Psi)$ is called a {\em Siegel modular embedding}, and will be used
to pull back classical theta functions, given in standard coordinates on the universal
family over $\bH_2/G_{P}$, to $X_{D}(\frakb)$. We note in passing that there are
two useful conventions for symplectic groups in the case of non-principal polarizations.
The other symplectic group
\bes
\Sp_{2g}^P(\bZ) \= \bigl\{ M \in \bZ^{2g \times 2g}\,:\,
M \cdot
\left(\begin{smallmatrix} 0  & P \\ -P & 0\end{smallmatrix}\right)
\cdot M^T \=
\left(\begin{smallmatrix} 0  & P \\ -P & 0\end{smallmatrix}\right)
\bigr\}\,.
\ees
is convenient, since it has integral entries. Conjugation by
$\left(\begin{smallmatrix}I_2 & 0 \\ 0 & P\end{smallmatrix}\right)$
takes $\Sp_{2g}^P(\bZ)$ into~$G_p$. Whereas the action of $G_P$ is the
standard action, the group $\Sp_{2g}^P(\bZ)$ acts on $\bH_2$ by 
\bes
\Sp_{2g}^P(\bZ) \ni M = \left(\begin{smallmatrix} A & B \\
  C & D \end{smallmatrix}\right):
\quad Z \,\, \mapsto \,\, (AZ +BP) (P^{-1}CZ+P^{-1}DP)^{-1}
\ees
\par
\medskip
In order to construct Siegel modular embeddings, one needs
to find an appropriate $\bZ$-basis of $\frakb\oplus\cO_{D}^{\vee}$.
Let $\frakc$ be any fractional $\cO_{D}$-ideal and $\bm\eta=(\eta_1,\eta_2)$
an ordered basis of $\frakc^{\vee}$, and define the matrices
\be\label{eq:matrixB}
B=B_{\bm\eta} = \begin{pmatrix} \eta_1 & \eta_1^\sigma \\ \eta_2 & \eta_2^\sigma
\end{pmatrix}
\quad\mbox{and}\quad
C=(B_{\bm\eta}^{-1}P)^T = \begin{pmatrix} \nu_1 & \nu_1^\sigma \\
  \nu_2 & \nu_2^\sigma \end{pmatrix}\,.
\ee
We say that $(\eta_1,\eta_2)$ is a basis {\em symplectically adapted to~$P$}
(or a \emph{$(d_{1},d_{2})$-symplec\-tically adapted basis}) if $(\nu_{1},\nu_{2})$
is the basis of an $\cO_D$-ideal. In this case we may factor 
the ideal as~$\frakc\frakb$, where $\frakb$ is necessarily an ideal of
norm~$n=d_{1}\cdot d_{2}$. Accordingly, the basis $\bm\eta$ determines the
rank-2 $\cO_{D}$-module
$\frakc\frakb\oplus\frakc^{\vee}$, that, provided with the trace pairing,
becomes a $(d_{1},d_{2})$-polarized module with symplectic basis
	\[(\nu_{1},0)\,,\ (\nu_{2},0)\,,\ (0,\eta_{1})\,,\ (0,\eta_{2})\,.\]
We do not necessarily assume $d_{1}|d_{2}$ here.
\par
To give an example in the particular case of $\frakc=\cO_{D}$ and
$\frakb=\langle \tfrac{r+\sqrt{D}}{2},n \rangle$ an ideal of norm $n$, we can always use the basis $\bm\eta=\tfrac{1}{\sqrt{D}}
\langle 1,\tfrac{-r+\sqrt{D}}{2}\rangle$ of $\cO_{D}^{\vee}$, which
is $(1,n)$-symplectically adapted to~$P_{1,n}$ and such that the first column
of $(B_{\bm\eta}^{-1}P_{1,n})^T$ agrees with the given basis of $\frakb$.
\par
The period matrix for $T_{\bm{\tau}} = \bC^2/\Lambda_{\frakb,\bm{\tau}}$ with respect to eigenforms for the $\cO_D$-action becomes
	\bes
	\Pi_{\bbu}\=\left(\begin{pmatrix}\tau_{1} & 0 \\ 0 & \tau_{2} \end{pmatrix}
	\cdot B_{\bm\eta}^{T} \ \middle\vert \ C^{T} \right).
	\ees
We refer to the corresponding coordinates of $\bC^2$ as 
\emph{eigenform coordinates} $\bbu=(u_{1},u_{2})$.
By multiplying on the left by $B_{\bm\eta}$, one gets the period matrix in
\emph{standard coordinates} $\bbv=B_{\bm\eta}\cdot \bbu$
\[\Pi_{\bbv}\=\left(\Omega_{\bm{\tau}}\ \middle\vert \
\begin{matrix} 1 & 0 \\ 0 & n \end{matrix}\right),\quad
\text{where} \quad \Omega_{\bm{\tau}}=B_{\bm\eta}\cdot
\begin{pmatrix}\tau_{1} & 0 \\ 0 & \tau_{2} 
\end{pmatrix}\cdot B_{\bm\eta}^{T}\in\bH_{2}\,.\]
Let us remark that, with the notation of Section~\ref{subsec:1npol}, one can assume that 
the columns of $\Pi_{\bbv}$ correspond to the lattice vectors
$\lambda_{1},\, \lambda_{2},\, \mu_{1},\, \mu_{2}$, respectively.
\par
We claim that the following is a well-defined homomorphism
\be \label{eq:SMEPsi}
\begin{array}{cccl}
	\Psi : & \SL(\frakb\oplus\cO_{D}^{\vee}) & \to & G_{P} \\[4pt]
	& \delta=\begin{pmatrix}a & b \\ c & d \end{pmatrix} & \mapsto & 
	\begin{pmatrix}B_{\bm\eta} & 0 \\ 0 & B_{\bm\eta}^{-T} \end{pmatrix}
	\begin{pmatrix}\widehat{a} & \widehat{b} \\ \widehat{c} & \widehat{d} \end{pmatrix}
	\begin{pmatrix}B_{\bm\eta}^{-1} & 0 \\ 0 & B_{\bm\eta}^{T} \end{pmatrix}\,.
	\end{array}\ee
Here we denote by $\widehat{k}$ the matrix $\left(\begin{smallmatrix}k & 0\\0& k^{\sigma}\end{smallmatrix}\right)$, for $k\in K$. The claim can be easily checked by studying the action on integral column vectors of the four blocks forming $\Psi(\delta)$.
\par
It is clear that $(\psi,\Psi)$ defined by $\psi : \bH^{2}  \to  \bH_{2}$, $\bm{\tau}=(\tau_{1},\tau_{2})  \mapsto  \Omega_{\bm{\tau}}$ and $\Psi$ as above induces the forgetful
map $X_{D}(\frakb)\to\cA_{2,(1,n)}$ and is therefore a Siegel modular embedding.

\

We finish this section with a criterion for some specific bases to be symplectically adapted. Recall that to given a triple of integers $Q=(a,b,c)$ such that $D=b^{2}-4ac$, one can associate the fractional ideal $\fraka^{\vee}=\langle 1, \lambda_{Q}\rangle$ of $\cO_{D}$, where $\lambda_{Q}=\frac{-b+\sqrt{D}}{2a}$ is the quadratic irrationality of $Q$.

\begin{lemma}\label{lem:symplecticallyadapted}Let $(d_{1},d_{2})$ be the type of a polarization such that $\gcd(d_{1},d_{2})=1$.\\
The basis $(1, \lambda_{Q})$ is a $(d_{1},d_{2})$-symplectically adapted basis of $\fraka^{\vee}$ if and only if $a\equiv 0 \bmod{d_{1}}$ and $c\equiv 0 \bmod{d_{2}}$. Moreover, $\fraka\frakb=\frac{a}{\sqrt{D}}\left\langle d_{2},-d_{1}\lambda^{\sigma}\right\rangle$.
\end{lemma}

Note that the choice of the type $(d_{1},d_{2})$ of the polarization does not follow the usual convention $d_{1}|d_{2}$ except in the case $d_{1}=1$.

\begin{proof}Let
	\[B = \begin{pmatrix} 1 & 1 \\ \frac{-b+\sqrt{D}}{2a} & \frac{-b-\sqrt{D}}{2a} \end{pmatrix} \quad\mbox{and}\quad
	P=\begin{pmatrix} d_{1} & 0 \\ 0 & d_{2} \end{pmatrix}\,.\]

The basis $(1, \lambda_{Q})$ is a $(d_{1},d_{2})$-symplectically adapted basis if and only if the columns of $(B^{-1}P)^{T}$ generate an ideal. This is equivalent to the existence of an integral matrix $R$ satisfying
	\[\begin{pmatrix}\frac{D+\sqrt{D}}{2} & 0 \\ 0 & \frac{D-\sqrt{D}}{2} \end{pmatrix} B^{-1}P = B^{-1}PR\,.\]
Now a simple calculation shows that 
	\[P^{-1}B\begin{pmatrix}\frac{D+\sqrt{D}}{2} & 0 \\ 0 & \frac{D-\sqrt{D}}{2} \end{pmatrix} B^{-1}P = \begin{pmatrix}\frac{D+b}{2} & \frac{a d_{2}}{d_{1}} \\ -\frac{c d_{1}}{d_{2}} & \frac{D-b}{2} \end{pmatrix}\,.\]
Since $b\equiv D\bmod{2}$ the claim follows. 
The generators of $\fraka\frakb$ correspond to the first column of the matrix $(B^{-1}P)^{T}$.
\end{proof}

\subsection{Cusps of \texorpdfstring{$X_{D}(\frakb)$}{X_D(b)}}
\label{subsec:cusps}

Cusps of $X_{D}(\frakb)$ are orbits of $\bP^{1}(K)$ under the action of $\SL(\frakb\oplus\cO_{D}^{\vee})$. Via the map $(\alpha: \beta)\mapsto \fraka=\alpha\cO_{D}+\beta \sqrt{D}\frakb^{-1}$, they correspond to ideal classes of invertible $\cO_{D}$-ideals $\fraka$ (see~\cite[\S I.4]{vdG}). In order to study the behavior of modular forms around the different cusps and to avoid the problem of changing coordinates in $\SL(\frakb\oplus\cO_{D}^{\vee})\backslash \bH^{2}$, one can instead change the Hilbert modular surface in the following way.

Let $\fraka$ be an invertible $\cO_{D}$-ideal. The trace pairing defined in the previous subsection induces again a symplectic pairing of type $(1,n)$ on the ``shifted'' $\cO_{D}$-module $\fraka\frakb\oplus\fraka^{\vee}$. In particular, one can define a lattice $\Lambda_{\frakb,\bm\tau}^{\fraka}$ for each $\bm{\tau}=(\tau_{1},\tau_{2})\in\bH^{2}$ as above and the Hilbert modular surface
\[X_{D}^{\fraka}\coloneqq X_{D}^{\fraka}(\frakb) \=
   {\SL(\fraka\frakb\oplus\fraka^{\vee})} \backslash \bH^2\,,\]
where 
	\[\SL(\fraka\frakb\oplus\fraka^{\vee})=
	\begin{pmatrix} \cO_{D} & \sqrt{D}\,\fraka^{2}\frakb \\ 
	\frac{1}{\sqrt{D}}\,\fraka^{-2}\frakb^{-1} & \cO_{D} \end{pmatrix}\cap \SL_{2}(K)\,,
	\]
parameterizes $(1,n)$-polarized abelian surfaces with a choice of real multiplication by $\cO_{D}$ too. In fact, for any element
\be\label{eq:cuspchangematrix}
M \= \begin{pmatrix} \alpha & \beta \\ \gamma & \delta
\end{pmatrix} \,\in \,
\begin{pmatrix} \fraka & \sqrt{D}\,\fraka\frakb \\
\frac{1}{\sqrt{D}}\,(\fraka\frakb)^{-1} & \fraka^{-1}
\end{pmatrix}\,\cap\, \SL_{2}(K)
\ee
the map
$$\phi: \bH^2 \to \bH^2, \quad (\tau_1,\tau_2) \mapsto (M\tau_1,M^{\sigma}\tau_2)$$
is equivariant with respect to the action of $U \in
\SL(\frakb\oplus\cO_{D}^{\vee})$ on its domain and
$M U M^{-1} \in \SL(\fraka\frakb\oplus\fraka^{\vee})$ on its range.
Via the map $\phi$ the cusp of $X_{D}(\frakb)$ corresponding to $\fraka$ is sent to the
cusp at infinity of $X_{D}^{\fraka}(\frakb)$.
\par
The matrices defined in the last section for the usual Hilbert modular group
can be changed accordingly. Let $\bm\xi=(\xi_1,\xi_2)$ now be an ordered
basis of~$\fraka_D^\vee$ that is symplectically adapted to~$P$ and such
that the first column of $(B_{\bm\xi}^{-1}P)^T$ forms a basis of the
ideal~$\fraka\frakb$. Then the matrix $B_{\bm\xi}$ determines a
Siegel modular embedding $(\psi_\fraka,\Psi_\fraka)$ by setting
$\psi_\fraka(\tau_1,\tau_2) = B_{\bm\xi} \left(\begin{smallmatrix} \tau_1 & 0 \\
0 & \tau_2 \end{smallmatrix}\right)B_{\bm\xi}^T$ and by defining
$\Psi_\fraka$ as in~\eqref{eq:SMEPsi}.
\par
As expected, by changing the cusp at infinity we are changing the Hilbert modular surface, but the Siegel modular embedding $(\psi_{\fraka},\Psi_{\fraka})$ and
the general one $(\psi,\Psi)$ constructed in the last section are compatible. 
\par
\begin{prop}Let $\bm\eta=(\eta_1,\eta_2)$ and $\bm\xi=(\xi_1,\xi_2)$ be symplectically adapted bases of $\cO_{D}^{\vee}$ and $\fraka^{\vee}$ determining the $\cO_{D}$-modules $\frakb\oplus\cO_{D}^{\vee}$ and $\fraka\frakb\oplus\fraka^{\vee}$ respectively. Moreover, let $M$ be the matrix in~\eqref{eq:cuspchangematrix} and define the matrix $\widetilde M = \left(\begin{smallmatrix}
a & b \\ c & d \\ \end{smallmatrix}\right)$ by
\bes
a = B_{\bm\xi} \widehat{\alpha} B_{\bm\eta}^{-1}, \quad b=B_{\bm\xi}\widehat{\beta}B_{\bm\eta}^T,
\quad c = B_{\bm\xi}^{-T}\widehat{\gamma}B_{\bm\eta}^{-1}, \quad d = B_{\bm\xi}^{-T} \widehat{\delta}B_{\bm\eta}^T\,.
\ees
Then $\widetilde{M}$ belongs to the symplectic group $G_P$ and the left action map
\[\widetilde{\psi} (\Omega) = \widetilde M \cdot \Omega := (a\Omega+b)(c\Omega+d)^{-1}\]
lifts the map~$\phi$ to the Siegel upper half space, i.e. $\widetilde{\psi} \circ \psi \= \psi_\fraka \circ \phi$.
\end{prop}
\par
\begin{proof}Proceeding as in the last section, one can easily check that $\widetilde{M}\in G_P$.
Now, by definition and using the abbreviation $\bm\tau =(\tau_1,\tau_2)$ we have
\bas
 \psi_\fraka \circ \phi (\tau_1,\tau_2)  & \=  B_{\bm\xi}
\Bigl(
\left(\begin{smallmatrix}\alpha & 0\\0& \alpha^{\sigma}\end{smallmatrix}\right)
\left(\begin{smallmatrix}\tau_1 & 0\\0& \tau_2\end{smallmatrix}\right)
+
\left(\begin{smallmatrix}\alpha & 0\\0& \alpha^{\sigma}\end{smallmatrix}\right)
\Bigr) \cdot
\Bigl(
\bigl(\begin{smallmatrix}\gamma & 0\\0& \gamma^{\sigma}\end{smallmatrix}
\bigr)
\left(\begin{smallmatrix}\tau_1 & 0\\0& \tau_2\end{smallmatrix}\right)
+
\left(\begin{smallmatrix}\delta & 0\\0& \delta^{\sigma}\end{smallmatrix}\right)
\Bigr)^{-1} B_{\bm\xi}^T \\
&\= 
\left(B_{\bm\xi}\widehat{\alpha} B_{\bm\eta}^{-1} \psi(\bftau) + B_{\bm\xi}\widehat{\beta}B_{\bm\eta}^T\right)
\cdot \left(B_{\bm\xi}^{-T}\widehat{\gamma}B_{\bm\eta}^{-1}\psi(\bftau)
+ B_{\bm\xi}^{-T} \widehat{\delta}B_{\bm\eta}^T\right)^{-1}
\eas
and thus the map $\widetilde\psi$ has the required commutation property.
\end{proof}
\par

%% file: sec_linebundles.tex
\section{Line bundles on \texorpdfstring{$(1,n)$}{(1,n)}-polarized
  abelian surfaces} \label{sec:linbd}

Classical theta functions are sections of line bundles on the abelian surface
$T = \bC^2/\Lambda$ where $\Lambda = \Pi \bZ^4$ is the period lattice
generated by the period matrix $\Pi = (\Omega,P_{1n})$. They are given by
the Fourier expansion 
\bes
  \th{c_{1}}{c_{2}}{2}:\bH_2 \times\bC^2 \to \bC, \quad 
	\th{c_{1}}{c_{2}}{2}(\Omega,\bbv) = 
	\sum_{\bbx\in\bZ^{2}+c_{1}} \bbe\left( \bbx^{T}\Omega\bbx \right) \bbe\left( 2\bbx^{T}(\bbv+c_{2}) \right)\,,
\ees
where $\bbe(t)=e^{\pi i t}$. (We consider all vectors inside the formula 
as column vectors).
The argument $c$ is called the {\em characteristic} of the theta
functions. Theta functions that differ only in their characteristics
correspond to sections of line bundles that are translates of
each other. For the moment we think of~$\Omega$ fixed and consider
the dependence on~$\Omega$ in the image of a Siegel modular embedding
starting from Section~\ref{sec:PD2T}
\par
The purpose of this section is to give a basis of sections of
a line bundle on a ${(1,n)}$-polarized abelian surface for
a characteristic chosen with the application in Lemma~\ref{lem:characteristic2}
in mind. Moreover we compute the Fourier expansions of these
line bundles with respect to a symplectically adapted basis. The main goal
are consequently the Fourier expansions in Proposition~\ref{prop:FEDi}
and the relation among the values of these theta-functions at two-torsion
points in~\eqref{eq:derivativesD2}. The miraculous reduction
of the number of constraints appearing in the next section
relies on this. 
\par
Most statements in this section are essentially in Sections 3.1, 4.6 and 4.7
of~\cite{birkenhakelange} that we rewrite for our purposes.
Since this reference use the (equivalent) language of canonical (as opposed
to classical) theta functions, we provide a short introduction
and conversion between the languages.

\subsection{Canonical theta functions}\label{subsec:cantheta}

Let~$V$ be a complex vector space and let~$\Lambda$ be a lattice in~$V$.
To a line bundle $\cL$ on the complex torus $T=V/\Lambda$ one associates its 
first Chern class $H=c_{1}(\cL)$, that we view as a  
Hermitian form on~$V$ whose imaginary part takes integral values on $\Lambda$.
To a line bundle $\cL$ one can associate a semicharacter 
$\chi: \Lambda \to S^1$ such that conversely~$\cL$ is the line bundle
associated with (cf.\ \cite[Appendix B]{birkenhakelange})
the {\em canonical factor of automorphy}
\be \label{eq:defaL}
a_\cL(\lambda,\bbu) \= \chi(\lambda)\exp(\pi H(\bbu,\lambda)) 
\in Z^1(\Lambda, H^0(\cO^*_T)), \quad \bbu \in V, \lambda \in \Lambda.
\ee
This correspondence can be made more concrete in the case that~$H$
is positive definite, i.e.\ the line bundle $\cL$ is ample on~$T$
and hence $T$ an abelian variety. A {\em decomposition} of~$V$ for~$H$
is a direct sum $V=V_{1}\oplus V_{2}$ so that 
$\Lambda_{i}\coloneqq V_{i}\cap \Lambda$ are isotropic with respect to
$E = \Im H$. For such a decomposition there is a standard semicharacter
	\be\label{eq:semicharacter}
	\chi_0(\bbu)\= \exp(\pi i E(u_1,u_2)), \quad \text{where} \quad
\bbu= u_1 + u_2, u_i \in V_i\,.
	\ee
with associated line bundle $\cL_0 = \cL(H,\chi_0)$. For every other
line bundle $\cL$ with $c_1(\cL) = H$ there is a point $c \in V$, 
such that $\cL = t_c^* \cL_0$.
The point is called the {\em characteristic} of $\cL$ for the chosen
decomposition. It is uniquely determined up to translation by an 
element in 
$$\Lambda(H) \=\{\bbu \in V\,|\, E
(\bbu,\lambda)\in\bZ\}.$$
(Here and in the sequel we often write e.g.\ $\Lambda(\cL)$ and $\Lambda(H)$
interchangeably for notions depending only on the first Chern class
of the line bundle.)
Consequently, characteristics for a given decomposition are in bijection
with $V/\Lambda(H)$.
\par
For a given line bundle~$\cL$ the global sections $H^0(T,\cL)$ can
be identified with functions
$ \vartheta: V \to \bC, \quad \vartheta(\bbu+\lambda)
= f(\lambda,u) \vartheta(u) $
where~ $f$ is a factor of automorphy for~$\cL$. More concretely, in the 
case $f = a_\cL$ as in~\eqref{eq:defaL} the functions
$$ \vartheta: V \to \bC, \quad \vartheta(\bbu+\lambda) 
\= a_\cL(\lambda,\bbu) \vartheta(\bbu) $$
are called {\em canonical theta functions} for~$\cL$, that we now construct.
We define for every $c \in V$
\ba \label{eq:thetafourier}
\vartheta^c(\bbu) &\= \exp \bigl(-\pi H(\bbu,c) -\tfrac{\pi}{2}H(c,c) 
+ \tfrac{\pi}{2}
\MB(\bbu+c,\bbu+c)\bigr) \\ 
& \phantom{\=} \cdot \sum_{\lambda \in \Lambda_1} \exp\bigl(\pi(H-\MB)(\bbu+c,\lambda)
- \tfrac{\pi}2 (H-\MB)(\lambda,\lambda)\bigr)\,,
\ea
where $\MB$ is the symmetric bilinear extension of $H|_{V_2}$.
For every $w \in K(\cL)$ we use the bilinear extension 
\be \label{eq:aL}
a_\cL(\bbu,\bbv) \= \chi_{0}(\bbu)\exp\left(2\pi i E(c,\bbu) + \pi H(\bbv,\bbu) + \frac{\pi}{2}H(\bbu,\bbu)\right)
\ee
of $a_\cL$ to a function $V \times V \to \bC$ and we set
\be \label{eq:defthetawc}
\vartheta_w^c(\bbu) \= a_\cL(w,\bbu)^{-1}\vartheta^c(\bbu+w)\,.
\ee
\par
Let us denote by $K(H)$ the kernel ${\ker}(\phi_{\cL})=\Lambda(H)/\Lambda$ of the canonical isogeny $\phi_\cL:T \to T^\vee$ defined by $\cL$. For the following theorem we note that the choice of a decomposition $V=V_{1}\oplus V_{2}$ induces direct sum decompositions of the lattice of integral points 
$\Lambda(H)=\Lambda(H)_{1}\oplus\Lambda(H)_{2}$ and of 
$K(H) = K(H)_1 \oplus K(H)_2$, where $K(H)_i = \Lambda(H)_i/(\Lambda \cap 
\Lambda(H)_i)$. In this notation \cite[Theorem~3.2.7]{birkenhakelange} gives:
\par
\begin{theorem} \label{thm:thetabasis}
The function $\vartheta^c_w$ is a canonical theta
function for $\cL = t_c^* \cL_0$. More precisely, if $c$ is a
characteristic with respect to a decomposition of~$V$ then the set
$\{ \vartheta^c_w\,:\, w \in K(\cL)_1\}$ is a basis of $H^0(\cL)$.
\end{theorem}
\par
Next we prove that actually the theta function $\vartheta^{c}_{w}$ only depends on the $K(\cL)_{1}$ component of $w$. This fact will be crucial to get extra relations between the values of theta 
functions at torsion points.
\par
\begin{lemma}\label{lem:thetaonw}
Let $w=w_{1}+w_{2}\in \Lambda(H)/\Lambda$. Then 
$\vartheta^{c}_{w}=\vartheta^{c}_{w_{1}}$.
\end{lemma}
\par
\begin{proof} The definition of the canonical theta function implies
\bes
\vartheta^{c}_{w}(\bbu) \= \exp\bigl(-\pi H(\bbu,c) - \tfrac{\pi}2 H(c,c) \bigr)
\vartheta^{0}_{w}(\bbu)
\ees
and hence it is enough to prove the claim for the characteristic~0.
By the definition~\eqref{eq:defthetawc} of~$\vartheta^{0}_{w}$ and the properties
of the factor $a_{\cL}$ (see~\cite[Lemma~3.1.3]{birkenhakelange})
\bas
\vartheta^{0}_{w}(\bbu) &\=  a_{\cL_{0}}(w_{1}+w_{2},\bbu)^{-1}
\vartheta^{0}_{0}(\bbu+w_{1}+w_{2}) \\
&\= a_{\cL_{0}}(w_{1},\bbu)^{-1}a_{\cL_{0}}(w_{2},w_{1}+\bbu)^{-1}
\vartheta^{0}_{0}(\bbu+w_{1}+w_{2}) \,.
\eas
Applying the Fourier expansion~\eqref{eq:thetafourier} of $\vartheta^{0}_{0}$, 
using~\eqref{eq:aL} and $\chi_0(w_2)=1$, we obtain
\bas
\vartheta^{0}_{w}(\bbu)
&= a_{\cL_{0}}(w_{1},\bbu)^{-1}\exp\left(-\pi(H-\MB)\Bigl(\bbu+w_{1}
+\tfrac{1}{2}w_{2},w_{2}\Bigr)+\tfrac{\pi}{2}\MB(\bbu+w_{1},\bbu+w_{1}) \right) \\
&\quad\ \cdot \sum_{\lambda\in\Lambda_{1}}\exp\left(\pi(H-\MB)(\bbu+w_{1},\lambda)
+\pi(H-\MB)(w_{2},\lambda)-\tfrac{\pi}{2}(H-\MB)(\lambda,\lambda)\right)\,.
\eas
Now~\cite[Lemma 3.2.2]{birkenhakelange} implies
$\pi(H-\MB)\Bigl(\bbu+w_{1}+\tfrac{1}{2}w_{2},w_{2}\Bigr)=0$,
since $w_{2}\in V_{2}$, and 
$\pi(H-\MB)(w_{2},\lambda)=2\pi i\, E(w_{2},\lambda)\,\in 2\pi i\,\bZ$,
since $w_{2}\in\Lambda(H)$.
Applying~\eqref{eq:thetafourier} and~\eqref{eq:defthetawc} again we obtain
	\begin{align*}
	\vartheta^{0}_{w}(\bbu) \=&
	a_{\cL_{0}}(w_{1},\bbu)^{-1}\exp\left(\tfrac{\pi}{2}\MB(\bbu+w_{1},\bbu+w_{1}) \right) \\
	&\cdot \sum_{\lambda\in\Lambda_{1}}\exp\left(\pi(H-\MB)(\bbu+w_{1},\lambda) -\tfrac{\pi}{2}(H-\MB)(\lambda,\lambda)\right)\\
	\=& a_{\cL_{0}}(w_{1},\bbu)^{-1}\vartheta^{0}_{0}(\bbu+w_{1})
	\=\vartheta^{0}_{w_{1}}(\bbu)\,
	\end{align*}
as claimed.
\end{proof}
\par

\subsection{Specialization to \texorpdfstring{$(1,6)$}{(1,6)}-polarization} \label{subsec:1dpol}
From now on we suppose $\dim(T) = 2$ and that $\cL$ is a line
bundle of {\em type $(1,6)$}, i.e.\  there exists a decomposition 
$V=V_{1}\oplus V_{2}$ for $H = c_1(\cL)$ and bases 
$\Lambda_{1}=\langle\lambda_{1},\lambda_{2}\rangle$ and 
$\Lambda_{2}=\langle\mu_{1},\mu_{2}\rangle$ 
in which $\Im H$ has a representation
\be \label{eq:ImHmatrix}
\Im H\=\begin{pmatrix} 0 & P \\ -P & 0 \end{pmatrix}\,,\quad
\text{where $P=\diag(1,6)$\,.}
\ee
Under these assumptions, $\Lambda(H)=\langle\lambda_{1},\tfrac16\lambda_{2},
\mu_{1},\tfrac16\mu_{2}\rangle$ and 
$K(H)=\Lambda(H)/\Lambda\cong \left(\bZ/6\bZ\right)^{2}$.
\par
Recall that a divisor $D$ on~$T$ is {\em symmetric} if $(-1)^*D = D$.
A line bundle $\cL$ is defined to be {\em symmetric} if
the corresponding semi-character~$\chi$ takes values in~$\pm 1$.
This notion is designed so that the line bundle $\cL = \cO(D)$ of
a symmetric divisor is symmetric (cf.\ \cite[Section~4.7]{birkenhakelange}).
For such a line bundle $(-1)^*$ induces an involution on $H^0(\cL)$, 
hence on the vector space generated by canonical theta functions.
\par
With the application to Prym varieties in mind, we focus on the line 
bundle $\cL =t_{c}^{*}{\cL}_{0}$ of characteristic 
$c=\tfrac12\lambda_{1}+\tfrac12\mu_{1}$. The space $H^{0}(\cL)$ is generated by $\{ \vartheta^{c}_{j\lambda_{2}/6}\,:\, j=0,\ldots,5 \}$ and, in this situation, the inverse formula~\cite[Formula~4.6.4]
{birkenhakelange}) gives $ (-1)^* \vartheta_w^c \= (-1) \cdot \vartheta_{-w}^c$
for all $w \in K(H)_1$. Consequently, the spaces of even and odd theta functions are given respectively by 
\ba \label{eq:thetaidef}
H^{0}(\cL)_{+} \= \Bigl\langle
&\vartheta^{c}_{\tfrac{1}{6}\lambda_{2}}-\vartheta^{c}_{\tfrac{5}{6}\lambda_{2}}\,,\ \vartheta^{c}_{\tfrac{2}{6}\lambda_{2}}-\vartheta^{c}_{\tfrac{4}{6}\lambda_{2}}
\Bigr\rangle\quad\mbox{and} \\
H^{0}(\cL)_{-} \= \Bigl\langle & \theta_{0}= \vartheta^{c}_{0}\,,\,
\theta_{1}= \vartheta^{c}_{\tfrac{1}{6}\lambda_{2}}+\vartheta^{c}_{\tfrac{5}{6}\lambda_{2}}\,, 
\theta_{2}= \vartheta^{c}_{\tfrac{2}{6}\lambda_{2}}+\vartheta^{c}_{\tfrac{4}{6}\lambda_{2}}\,,
\theta_{3}= \vartheta^{c}_{\tfrac{1}{2}\lambda_{2}}
\Bigr\rangle\,.
\ea
%
\par
We will need the following result relating the values of odd theta functions 
at certain 2-torsion points, more precisely the set of $2$-torsion
points in the kernel~$K(H)$ of the map $\phi_\cL$ to the dual torus.
\par
\begin{lemma}\label{lem:thetarelations} Let $\theta_{0}(\bbu),\ldots,
\theta_{3}(\bbu)$ be the generators of $H^{0}(\cL)_{-}$. Then
\bas
	\theta_{0}(\bbu) &= a_{\cL}\left(\tfrac{1}{2}\lambda_{2},\bbu\right)^{-1}\theta_{3}\left(\bbu+\tfrac{1}{2}\lambda_{2}\right) &&= & & a_{\cL}\left(\tfrac{1}{2}\mu_{2},\bbu\right)^{-1}\theta_{0}\left(\bbu+\tfrac{1}{2}\mu_{2}\right) \,, \\
	\theta_{1}(\bbu) &= a_{\cL}\left(\tfrac{1}{2}\lambda_{2},\bbu\right)^{-1}\theta_{2}\left(\bbu+\tfrac{1}{2}\lambda_{2}\right) &&= - & & a_{\cL}\left(\tfrac{1}{2}\mu_{2},\bbu\right)^{-1}\theta_{1}\left(\bbu+\tfrac{1}{2}\mu_{2}\right) \,, \\
	\theta_{2}(\bbu) &= a_{\cL}\left(\tfrac{1}{2}\lambda_{2},\bbu\right)^{-1}\theta_{1}\left(\bbu+\tfrac{1}{2}\lambda_{2}\right) &&=  & & a_{\cL}\left(\tfrac{1}{2}\mu_{2},\bbu\right)^{-1}\theta_{2}\left(\bbu+\tfrac{1}{2}\mu_{2}\right) \,, \\
	\theta_{3}(\bbu) &= a_{\cL}\left(\tfrac{1}{2}\lambda_{2},\bbu\right)^{-1}\theta_{0}\left(\bbu+\tfrac{1}{2}\lambda_{2}\right) &&= - & & a_{\cL}\left(\tfrac{1}{2}\mu_{2},\bbu\right)^{-1}\theta_{3}\left(\bbu+\tfrac{1}{2}\mu_{2}\right) \,.
\eas
\end{lemma}
\par
\begin{proof} For any $w=w_{1}+w_{2}$ and $\widetilde{w}=\widetilde{w}_{1}+\widetilde{w}_{2}\in\Lambda(H)/\Lambda$ we find, using~\eqref{eq:defthetawc}
and the transformation law of the canonical factor of automorphy (cf.\ Exercise 3.7(2)
in~\cite{birkenhakelange}), that 
\bes
\vartheta^{c}_{w}(\bbu) \=
\exp\left(2\pi\i \Im H(\widetilde{w}_{1},\widetilde{w}_{2}-w_{2})\right)a_{\cL_{X}}(w-\widetilde{w},\bbu)^{-1}\vartheta^{c}_{\widetilde{w}}(\bbu+w-\widetilde{w})\,.
\ees
\par
The first equalities claimed in the lemma are a direct application 
of this formula to $\widetilde{w}=\tfrac{j}{6}\lambda_{2}$ 
and $w=\tfrac{j+3}{6}\lambda_{2}$, where indices should be taken $\mod 6$.
\par
The second ones follow from the same formula applied to $\widetilde{w}=\tfrac{j}{6}\lambda_{2}$ and $w=\tfrac{j}{6}\lambda_{2} + \tfrac{1}{2}\mu_{2}$ together with the fact that, by Lemma~\ref{lem:thetaonw}, $\theta^{c}_{\widetilde{w}}=\theta^{c}_{w}$.
\end{proof}
\par

\subsection{Partial derivatives at two-torsion points} \label{sec:PD2T}

So far the computations were for a general abelian surface
and we now restrict to real multiplication loci, i.e.\ to a period
matrix $\Omega_{\bm{\tau}} = \psi(\bm{\tau})$ in the image of a
Siegel modular embedding determined by a $(d_1,d_2)$-symplectically
adapted basis $(\omega_1,\omega_2)$ as in Section~\ref{sec:SME}.
Since on a surface with real multiplication there are two eigendirections,
that we have given the coordinates $u_i$, for a general theta
function~$\vartheta$ the partial derivatives
\bes
D_{i}\vartheta(\bm\tau,\bbu_{0}) \coloneqq \frac{\partial}{\partial u_{i}}\vartheta(\bm\tau,\bbu)|_{\bbu=\bbu_{0}}\,,
\ees
will be of particular interest in the sequel.
As a direct consequence of Lemma~\ref{lem:thetarelations} together with the
vanishing of the $\theta_{j}$ at the given 2-torsion points we obtain
the analogous results for the derivatives $D_{i}\theta_{j}$, for $i=1,2$
\ba \label{eq:derivativesD2}
	D_{i}\theta_{0}(0) &= a_{\cL}\left(\tfrac{1}{2}\lambda_{2},0\right)^{-1}D_{i}\theta_{3}\left(\tfrac{1}{2}\lambda_{2}\right) &&= & & a_{\cL}\left(\tfrac{1}{2}\mu_{2},0\right)^{-1}D_{i}\theta_{0}\left(\tfrac{1}{2}\mu_{2}\right) \,, \\
	D_{i}\theta_{1}(0) &= a_{\cL}\left(\tfrac{1}{2}\lambda_{2},0\right)^{-1}D_{i}\theta_{2}\left(\tfrac{1}{2}\lambda_{2}\right) &&= - & & a_{\cL}\left(\tfrac{1}{2}\mu_{2},0\right)^{-1}D_{i}\theta_{1}\left(\tfrac{1}{2}\mu_{2}\right) \,, \\
	D_{i}\theta_{2}(0) &= a_{\cL}\left(\tfrac{1}{2}\lambda_{2},0\right)^{-1}D_{i}\theta_{1}\left(\tfrac{1}{2}\lambda_{2}\right) &&=  & & a_{\cL}\left(\tfrac{1}{2}\mu_{2},0\right)^{-1}D_{i}\theta_{2}\left(\tfrac{1}{2}\mu_{2}\right) \,, \\
	D_{i}\theta_{3}(0) &= a_{\cL}\left(\tfrac{1}{2}\lambda_{2},0\right)^{-1}D_{i}\theta_{0}\left(\tfrac{1}{2}\lambda_{2}\right) &&= - & & a_{\cL}\left(\tfrac{1}{2}\mu_{2},0\right)^{-1}D_{i}\theta_{3}\left(\tfrac{1}{2}\mu_{2}\right) \,.
\ea
Another consequence is that an odd theta function behaves near those
non-trivial two-torsion points like an odd function in the following sense.
\par
\begin{cor} \label{cor:multiplicitycor}
 Let $f\in H^{0}(\cL)_{-}$ be an odd theta function, let $Q$ be one
of the two-torsion points $\{0,\tfrac{1}{2}\lambda_{2},\tfrac{1}{2}\mu_{2},
\tfrac{1}{2}(\lambda_{2}+\mu_{2})\}$ and fix $i=1$ or $i=2$. If
$D_{i}^{2k-1}f(Q)=0$ for all $k=1,\ldots,n$, then $D_{i}^{2n}f(Q)=0$.
\end{cor}
\par
\begin{proof}The proof is trivial for $Q=0$ since $f$ is an odd function
of $\bC^{2}$. To discuss the other two-torsion points, write $f=f_{1}+f_{2}$,
where $f_{1}\in\langle\theta_{0},\theta_{2}\rangle$
and $f_{2}\in\langle \theta_{1},\theta_{3}\rangle$.
For $Q=\tfrac{1}{2}\mu_{2}$ we can write for each $N>0$
\bes
D_{i}^{N}f(\tfrac{1}{2}\mu_{2}) \=
\sum_{j=0}^{N} {{N}\choose{j}} D_{i}^{N-j}a_{\cL}(\tfrac{1}{2}\lambda_{2},0)
\left( D_{i}^{j}f_{1}(0) - D_{i}^{j}f_{2}(0)\right)\,
\ees
by Lemma~\ref{lem:thetarelations}. Since we work with a space of
odd theta functions, $D_{i}^{2k}\theta_{j}(0)=0$ for every $k$ and
$j\in\{0,1,2,3\}$. Consequently, we can use this formula inductively
to show that the hypothesis $D_{i}^{2k-1}f(\tfrac{1}{2}\mu_{2})=0$ for all
$k=1,\ldots,n$ holds if and only if $D_{i}^{2k-1}f_{1}(0) - D_{i}^{2k-1}f_{2}(0)=0$
for all $k=1,\ldots,n$. As a consequence, 
\[D_{i}^{2n}f(\tfrac{1}{2}\mu_{2}) = \sum_{k=1}^{n} {{2n}\choose{2k-1}} D_{i}^{2n-2k+1}a_{\cL}(\tfrac{1}{2}\lambda_{2},0)  \left( D_{i}^{2k-1}f_{1}(0) - D_{i}^{2k-1}f_{2}(0)\right)=0\,.\]
\par
For $Q=\tfrac{1}{2}\lambda_{2}$, we write $\widetilde{f}_{j}$ for $f_{j}$
with $\theta_{0}$ and $\theta_{1}$ exchanged with $\theta_{3}$ and $\theta_{2}$
respectively. With this notation
\[
D_{i}^{N}f(\tfrac{1}{2}\lambda_{2}) \= \sum_{j=0}^{N} {{N}\choose{j}} D_{i}^{N-j}a_{\cL}(\tfrac{1}{2}\lambda_{2},0)  \left( D_{i}^{j}\widetilde{f}_{1}(0) - D_{i}^{j}\widetilde{f}_{2}(0)\right)\,,
	\]
Again, the hypothesis $D_{i}^{2k-1}f(\tfrac{1}{2}\lambda_{2})=0$ for all
$k=1,\ldots,n$ holds if and only if $D_{i}^{2k-1}\widetilde{f}_{1}(0) - D_{i}^{2k-1}\widetilde{f}_{2}(0)=0$ for all $k=1,\ldots,n$. Consequently, 
	\[D_{i}^{2n}f(\tfrac{1}{2}\lambda_{2}) = \sum_{k=1}^{n} {{2n}\choose{2k-1}} D_{i}^{2n-2k+1}a_{\cL}(\tfrac{1}{2}\lambda_{2},0)  \left( D_{i}^{2k-1}\widetilde{f}_{1}(0) - D_{i}^{2k-1}\widetilde{f}_{2}(0)\right)=0\,.\]
The proof for $Q=\tfrac{1}{2}(\lambda_{2}+\mu_{2})$ follows the same lines.
\end{proof}

\subsection{Fourier expansions}

For a $(1,6)$-symplectically adapted basis $\bm\eta = (\eta_1,\eta_2)$
we define $\rho_{\bm\eta}(x_1,x_2)=x_{1}\eta_{1}+x_{2}\eta_{2}$, hence
$\bbx^{T}B_{\bm\eta}= \left(\rho_{\bm\eta}(\bbx),\rho_{\bm\eta}^{\sigma}(\bbx)\right)$
for the matrix~$B_{\bm\eta}$ used in~\eqref{eq:matrixB} to define a Siegel
modular embedding of the Hilbert modular surface~$X_D(\frakb)$.
Recall that the choice of such a basis also determines a decomposition
of~$V$ using
\be \label{eq:V1V2}
V_1 \= \langle (\nu_{1},0)\,, (\nu_{2},0)\rangle_\bR \,,\quad
V_2 \= \langle (0,\eta_{1})\,, (0,\eta_{2})\rangle_\bR\,.
\ee
We moreover define the shifted lattice $\Lambda_{\epsilon,\delta}=\bZ^{2}+(\epsilon,\delta)^{T}$ and abbreviate
$\rho = \rho_{\bm\eta}$ if~${\bm\eta}$ has been fixed.
\par
\begin{prop} \label{prop:FEDi}
The Nullwerte of the derivatives of the theta functions
$\theta_j$ for $j \in \{0,1,2,3\}$, as defined in~\eqref{eq:thetaidef}, have the
Fourier expansion
\ba \label{eq:thetaseries}
	\frac{\partial}{\partial u_{1}}\theta_j(\bm\tau,0)
	&= 2\pi\i\sum_{\bbx\in\Lambda_{\frac{1}{2},\frac{j}{6}}}  \bbe\left(x_{1}\right) \rho(\bbx) q_{1}^{\rho(\bbx)^{2}} q_{2}^{\rho^{\sigma}(\bbx)^{2}}, \\
	\frac{\partial}{\partial u_{2}}\theta_j(\bm\tau,0)
	&= 2\pi\i\sum_{\bbx\in\Lambda_{\frac{1}{2},\frac{j}{6}}}  \bbe\left(x_{1}\right) \rho^{\sigma}(\bbx) q_{1}^{\rho(\bbx)^{2}} q_{2}^{\rho^{\sigma}(\bbx)^{2}}, \\
\ea
where $q_i = \bbe(\tau_i)$ and $\bbe(\cdot) = \exp(\pi i\, \cdot)$.
\end{prop}
\par
\begin{proof} By  \cite[Lemma~8.5.2]{birkenhakelange} the canonical theta
function with characteristic~$c$ is given by 
\bes
	\vartheta^{c}(\bm\tau,\bbv) = e^{\frac{\pi}{2}B(\bbv,\bbv) -\pi i c_1^T c_2}\, \th{c_1}{c_2}{3}(\bm\tau,\bbv)\,
\ees
in terms of classical theta functions. We differentiate this, use that
the $\theta_j$ are odd, hence vanish at zero and use the Fourier expansions
\bas
\frac{\partial}{\partial u_{1}}\th{(\frac{1}{2},\frac{j}{6})}{(\frac{1}{2},0)}{3}(\bm\tau,0)
&= 2\pi\i\sum_{\bbx\in\Lambda_{3j}}  \bbe\left(x_{1}\right) \rho(\bbx) q_{1}^{\rho(\bbx)^{2}} q_{2}^{\rho^{\sigma}(\bbx)^{2}}\,, \\
\frac{\partial}{\partial u_{2}}\th{(\frac{1}{2},\frac{j}{6})}{(\frac{1}{2},0)}{3}(\bm\tau,0)
&= 2\pi\i\sum_{\bbx\in\Lambda_{3j}} \bbe\left(x_{1}\right) \rho^{\sigma}(\bbx) q_{1}^{\rho(\bbx)^{2}} q_{2}^{\rho^{\sigma}(\bbx)^{2}}\,.
\eas
This immediately gives the expansion for $\theta_0$ and $\theta_3$.
For the two remaining generators we moreover use that
	\[
	\frac{\partial}{\partial u_{i}}\th{(\frac{1}{2},\frac{-j}{6})}{(\frac{1}{2},0)}{3}(\bm\tau,0) 
	\= \frac{\partial}{\partial u_{i}}\th{(\frac{1}{2},\frac{j}{6})}{(\frac{1}{2},0)}{3}(\bm\tau,0)\,,\quad\mbox{for $j=1,2$.}
\]
as we see by changing the order of summation in~\eqref{eq:thetaseries}
using the observation that~$\rho$ is odd.
\end{proof}
\par
\medskip

\subsection{Derivatives of theta functions as Hilbert modular forms}

The set of all Siegel theta functions for characteristics
in $\tfrac1N \bZ$ (with $N$ fixed) satisfies a modular transformation law,
(see~\cite[Section 8.4]{birkenhakelange} for the complete formula). This
implies that the restriction via a Siegel modular embedding satisfies
a modular transformation law for the Hilbert modular group. In general,
this action still permutes characteristics, but here we make use of
the following fact.
\par
\begin{lemma}\label{lem:actiononcharacteristic}The space $H^{0}(\cL)$ of theta functions of characteristic $c=\tfrac12\lambda_{1}+\tfrac12\mu_{1}$ is preserved by the whole modular group $\SL(\frakb\oplus\cO_{D}^{\vee})$.
\end{lemma}
\par
\begin{proof}The action of the modular group on characteristics preserves the
set of characteristics corresponding to symmetric line bundles, and the action
on theta functions preserves the even and odd subspaces. Let $\cL$ be a
symmetric line bundle of characteristic $c$ that provides the $(1,6)$-polarization. 
Since $h^{0}(L)=6$ the space of odd theta functions of $L$ has dimension
	\[h^{0}_{-}=\frac{1}{2}(6-\# S)+\# S^{-}\,,\]
by~\cite[Proposition~4.6.5]{birkenhakelange}, where
$S=\{\overline{w}\in K(L)_{1}\,:\, 2\overline{w}=2\overline{c}_{1}\}$ and
$S^{-}=\{\overline{w}\in S\,:\, e(4\pi i \Im H(w+c_{1},c_{2}))=-1\}$.
One now computes that line bundle of characteristic $\tfrac{1}{2}\lambda_{1}+
\tfrac{1}{2}\mu_{1}$ is the only one with a 4-dimensional space of odd theta
functions. Thus every element of the modular group fixes this characteristic.
\end{proof}
\par
\medskip
Recall that a Hilbert modular form~$f$ of bi-weight~$(k,\ell)$ with
character~$\chi$ for the subgroup $\Gamma$ of a Hilbert modular group
is a holomorphic function~$f: \bH^2 \to \bC$ with the transformation law
\bes
f(\gamma \tau_1, \gamma^\sigma \tau_2) \= \chi(\gamma)(c\tau_1+d)^k
(c^\sigma \tau_2 +d^\sigma)^\ell f(\tau_1,\tau_2)
\ees
for all $\sm{a}bcd \in \Gamma$. 
The specialization of the theta transformation law implies that for
an even theta function $\vartheta$ of characteristic~$c$ the Nullwert
$\vartheta(\bm\tau,0)$ is a Hilbert modular form of
bi-weight~$(\tfrac12,\tfrac12)$ with some finite character
for some finite index subgroup of the Hilbert modular group. The partial
derivatives $D_{1}\vartheta(\bm\tau,0)$ and $D_{2}\vartheta(\bm\tau,0)$ of odd theta functions are
modular forms of bi-weight $(3/2,1/2)$ and $(1/2,3/2)$ respectively,
see~\cite[Section~9]{MZ} for the details.
\par

\subsection{Line bundles of type \texorpdfstring{$(2,3)$}{(2,3)}}\label{subsec:1npol}

The usual convention for the type $(d_{1},d_{2})$ of a polarization is that $d_{1}|d_{2}$. However, 
it will be convenient in our particular case to consider also the polarization type $(2,3)$ 
rather than just of type $(1,6)$. In this subsection we translate the results between the two 
different conventions.
\par
Let $\cL$ be a line bundle of type $(1,6)$ and let $V=V_{1}\oplus V_{2}$ be a decomposition for~$\cL$, 
so that $\Lambda=\Lambda_{1}\oplus \Lambda_{2} = \langle \lambda_{1},\lambda_{2}\rangle \oplus \langle 
\mu_{1},\mu_{2}\rangle$ gives a symplectic basis of the lattice with canonical $(1,6)$-symplectic matrix, 
i.\ e.\ the non-trivial intersection are $E(\lambda_{1},\mu_{1})=1$ and $E(\lambda_{2},\mu_{2})=6$. 
\par
The matrices $R_{\lambda}=\left(\begin{smallmatrix}2&1\\3&1\end{smallmatrix}\right)$ and 
$R_{\mu}=\left(\begin{smallmatrix}-2&1\\3&-1\end{smallmatrix}\right)$ give a change of basis to a symplectic basis $\langle \widetilde{\lambda}_{1},\widetilde{\lambda}_{2}\rangle \oplus \langle \widetilde{\mu}_{1},\widetilde{\mu}_{2}\rangle$ with canonical $(2,3)$-symplectic matrix while preserving the
chosen decomposition of~$V$. In particular we may identify the characteristics in the two situations
and we may identify the basis elements of $H^0(\cL)$ named in~\eqref{eq:thetaidef} in the two 
conventions. The  distinguished characteristic $c=\frac{1}{2}\lambda_{1} + \frac{1}{2}\mu_{1} \in 
\tfrac12\Lambda(H)/\Lambda(H)$ is expressed in the new basis as 
$c = \frac{1}{6}\widetilde{\lambda}_{2} + \frac{1}{6}\widetilde{\mu}_{2}$ 
since $\Lambda(H) = \langle \tfrac12 \widetilde{\lambda}_{1},  \tfrac13  \widetilde{\lambda}_{2},
\tfrac12 \widetilde{\mu}_{1}, \tfrac13\widetilde{\mu}_{2}\rangle$.
\par
Let now $\bm\eta=(\eta_1,\eta_2)$ be a $(2,3)$-symplectically adapted basis
of~$\fraka^\vee$, determining the $\cO_{D}$-module
$\fraka\frakb\oplus\fraka^{\vee}$, and consider the Siegel modular embedding
given by the matrix
$B\coloneqq B_{\bm\eta}=\left(\begin{smallmatrix} \eta_1 & \eta_1^\sigma \\
\eta_2 & \eta_2^\sigma \end{smallmatrix}\right)$ as in Section~\ref{subsec:cusps},
so that the cusp $\fraka$ of $X_{D}(\frakb)$ corresponds to the cusp at 
infinity of $X_{D}^{\fraka}(\frakb)$. Then $\bm\eta R_\mu^{-1}$ is $(1,6)$-symplectically 
adapted, and this base change together with the action of $R_\lambda^{-1}$ on $\bm\nu$ preserves
the decomposition~\eqref{eq:V1V2}, so that we are in indeed in the situation considered above.
\par
\begin{lemma} With $\rho(\bbx) = \rho_{\bm\eta}(\bbx)\coloneqq x_{1}\eta_{1}
+ x_{2}\eta_{2}$ stemming from a $(2,3)$-symplec\-tically adapted basis
the global sections $\theta_j \in H^0(\cL)$ of the line bundle
with characteristic $c = \frac{1}{6}\widetilde{\lambda}_{2} + \frac{1}{6}\widetilde{\mu}_{2}$ 
have Fourier expansions as in~\eqref{eq:thetaseries} 
with the lattice coset $\Lambda_{\frac{1}{2},\frac{j}{6}}$ for the series~$\theta_j$ replaced by
$\Lambda_{\frac{j}{2},\frac{2j+3}{6}}$ and the character $\bbe(x_1)$ replaced by $\bbe(x_2)$.
\end{lemma}
\par
\begin{proof} By definition, the lattice coset is $(1/2, j/6)$ in the basis $\lambda_1,\lambda_2$, 
which is equal to $(j/2, (2j+3)/6)$ in the basis $\widetilde{\lambda}_1,\widetilde{\lambda}_2$
and the character is determined by the $\mu$-component of the characteristic.
\end{proof}

%% file: sec_thetagothic.tex
\section{The Gothic modular form and the Gothic theta function}
\label{sec:GothMF}

We now specialize again to curves $(X,\pi_A,\pi_B)$ in the Gothic locus.
Abel-Prym maps denote, in analogy to the classical Abel-Jacobi map from
a curve to its Jacobian, the map from~$X$ to its Prym variety. Since
the Prym variety is not principally polarized, there are two natural choices
that we analyze here: to the Prym variety and to its dual. The main player
is the pre-Abel-Prym map $\varphi: X \to \Prym^\vee(X,\pi_A,\pi_B)$
to the dual Prym variety defined in Section~\ref{sec:complementary}. Since the Prym variety $\Prym^\vee(X)$ of a point in $G_{D}$ admits real multiplication by $\cO_{D}$, we can see the \Teichmuller curve $G_{D}$ inside some Hilbert modular surface $X_{D}(\frakb)$. Let us denote by $G_{D}(\frakb)$ the union of those components of the Torelli-image of $G_{D}$ in $X_{D}(\frakb)$ for which $du_{1}$ induces the eigenform $\omega$ at each point $(X,\omega)$.  
\par
Our goal is to describe $\varphi(X)$ in terms of theta functions
and nearly determine the Torelli-image of~$G_D$. 
\par
\begin{theorem} \label{thm:GDmodform}
The Torelli-image $G_D(\frakb)$ is contained in the vanishing locus of
the Hilbert modular form
\bes \label{eq:modularform}
\cG_{D}(\bm{\tau})\coloneqq D_2\theta_{0}(\bm\tau,0) \cdot D_2 \theta_{1}
(\bm{\tau},0) - D_2 \theta_{2}(\bm\tau,0) \cdot D_2 \theta_{3}(\bm{\tau},0)\,.
\ees
of bi-weight $(1,3)$. Consider the locus
\bes
\tRED(\frakb) =  \{\cG_{D}(\bm{\tau}) = 0 \} \cap \{\cF_a(\bm{\tau}) = 0 \}
\cap \{ \cF_b(\bm{\tau}) = 0\} \,,\ees
where we define the modular forms
\bas
\cF_a(\bm{\tau}) &= D_1\theta_{0}(\bm\tau,0)\cdot D_2\theta_{2}(\bm{\tau},0)
- D_1\theta_{2}(\bm\tau,0)\cdot D_2 \theta_{0}(\bm{\tau},0)\quad\mbox{and} \\
\cF_b(\bm{\tau}) &= D_1\theta_{1}(\bm\tau,0)\cdot D_2\theta_{3}(\bm{\tau},0)
- D_1\theta_{3}(\bm\tau,0)\cdot D_2 \theta_{1}(\bm{\tau},0)\,.
\eas
Then for all points in $\{\cG_{D}(\bm{\tau}) = 0 \} \setminus \tRED(\frakb)$
the theta function 
\bes
\theta_{X}(\bbu) \=
\left| \begin{array}{cccc}
\theta_{0}(\bbu) & \theta_{1}(\bbu) & \theta_{2}(\bbu) & \theta_{3}(\bbu) \\
D_{1}\theta_{0}(0) & D_{1}\theta_{1}(0) & D_{1}\theta_{2}(0) & D_{1}\theta_{3}(0) \\
D_{2}\theta_{0}(0) & 0 & D_{2}\theta_{2}(0) & 0 \\
0 & D_{2}\theta_{1}(0) & 0 & D_{2}\theta_{3}(0) \\
\end{array}\right|
\ees
is non-zero and the vanishing locus of this theta function is equal to
the pre-Abel-Prym
image $\varphi(X)$ of a Gothic Veech surface~$X$.
\end{theorem}
\par
We will discuss the exceptional set where the modular forms
$\cG_D$, $\cF_a$ and $\cF_b$ jointly vanish in Section~\ref{sec:modularcurves}. It is part of the reducible locus, as suggested by the notation
and as we will see in Proposition~\ref{prop:RedisRed}.

\subsection{The Abel-Prym map and the pre-Abel-Prym map}

Let $(X,\omega, \pi_A, \pi_B)$ be a flat surface in the Gothic
locus $\Omega G$. For each choice of a `base point' $p\in X$
there is the usual Abel-Jacobi map $\alpha_{p}: X \to \Jac X$
centered at $p$. We will fix once and for all the center of the
Abel-Jacobi map to be $p=p_{1}$ one of the fixed points of~$J$
where $\omega$ does not vanish.
\par
We have defined $\Prym(X) = \Prym(X,\pi_A,\pi_B)$ as the subvariety complementary to $A^\vee \times B^\vee$, hence there is a natural inclusion
$\iota: \Prym(X) \to \Jac(X)$. Its dual is thus a quotient map
$\iota^\vee: \Jac(X) \to \Prym^\vee(X)$ and the norm endomorphism $N_P$
defined in Proposition~\ref{prop:NYNP} is also such a quotient map. Using~\eqref{eq:NYdef}
we conclude that they fit into the following commutative diagram
\bes
\begin{tikzpicture}
\matrix (m) [matrix of math nodes, row sep=1.75em, column sep=3.5em, text height=1.5ex, text depth=0.25ex]
{ 
X & \Jac X & \Prym X \\
& & \Prym^{\vee}\!X \\
};
\path[->,font=\scriptsize]
(m-1-1) edge [above] node {$\alpha_{p}$} (m-1-2)
(m-1-1) edge [out=30,in=150] node [above] {$\abel$} (m-1-3)  
(m-1-2) edge [above] node {$N_{P}$} (m-1-3)
(m-1-2) edge [above] node {$\check{\iota}$} (m-2-3)
(m-2-3) edge [right] node {$\psi_{\cL}$} (m-1-3)
(m-1-1) edge [out=315,in=180] node [above] {$\pabel$} (m-2-3)  
;
\end{tikzpicture}
\ees
\par
The composition $\abel\coloneqq N_{P}\circ\alpha_{p_{1}}$ is called the
\emph{Abel-Prym map} and the composition
$\pabel\coloneqq\check{\iota}\,\alpha_{p_{1}}$ is called
the \emph{pre-Abel-Prym map} centered at $p_{1}$, respectively.
By Proposition~\ref{prop:NYNP} we can write the Abel-Prym map in terms of divisors as
	\bes
	\abel(x)\=\left[x^{(1)} - 3\,J(x^{(1)}) - 2\,x^{(2)} - 2\,x^{(3)} + 2\,p_{1} + 2\,q_{1} + 2\,J(q_{1}) \right]\,.
	\ees
Moreover $\abel(x)=\abel(y)$ if and only if
	\begin{equation}\label{eq:abelprymnotinjective}
	x^{(1)} - 3\,J(x^{(1)}) - 2\,x^{(2)} - 2\,x^{(3)} -y^{(1)} + 3\,J(y^{(1)}) + 2\,y^{(2)} + 2\,y^{(3)} \sim 0\,.
	\end{equation}
As a consequence of this formula we obtain:
\par
\begin{lemma}\label{lem:AbelPrym}The Abel-Prym map $\abel$ maps
$Z\cup P$ to a single point, i.e.\
	\[\abel(z_{i})= \abel(p_{i})= 0\,\quad\mbox{for $i=1,2,3$.}\]
\end{lemma}
\par
\begin{proof} Using~\eqref{eq:abelprymnotinjective} and the fact
that points in $Z\cup P$ are fixed under~$J$ the claim is equivalent to
$$ 2p_i^{(1)} + 2p_i^{(2)}+ 2p_i^{(3)} \sim 
2p_j^{(1)} + 2p_j^{(2)}+ 2p_j^{(3)} $$
and 
$$ 2p_i^{(1)} + 2p_i^{(2)}+ 2p_i^{(3)} \sim 
2z_j^{(1)} + 2z_j^{(2)}+ 2z_j^{(3)} $$
for any $p_i,p_j \in P$ and $z_j \in Z$.
This follows from the preimage diagram~\eqref{eq:hpreim} and the fact that
each of the points involved appears with coefficient~$2$ in $h^*(e_i)$.
\end{proof}
\par
\subsection{The natural line bundles on $\Prym^\vee(X)$}

There are several natural line bundles on the Prym varieties.
The restriction of the principal polarization on $\Jac X$ to $\Prym X$
via~$\iota$ yields a polarization of type $(1,6)$ that we denote given
by a line bundle~$\cL$. But we are rather interested in $\Prym^\vee(X)$.
There, we first have the bundle $\cL_{X}\coloneqq\cO_{\Prym^{\vee}\!X}(\pabel(X))$
generated by the image of the Gothic Veech surface that we are mainly
interested in. Second, there is the following general construction.
\par
Let $H = c_1(\cL$) and let $\phi_{\cL}:\Prym X \to  \Prym^{\vee}\!X$ be the 
isogeny associated with~$\cL$. Since $\cL$ is of type $(1,6)$ there is an 
isogeny $\psi:\Prym^{\vee}\!X \to \Prym X$ such that 
$\psi \circ\phi_{\cL}=[6]$ (cf.~\cite[Section~14.4]{birkenhakelange}). 
More precisely, $\psi = \psi_{\check{\cL}}$ for a line bundle $\check{\cL}$
on~$\Prym^{\vee}\!X$, well defined only up to translations, with the
same polarization $H = c_1(\check{\cL})$. To fix a precise point
of reference, we fix a decomposition for the universal covering $V$
of $\Prym^{\vee}\!X$
in which~$\Im H$ has the form~\eqref{eq:ImHmatrix}. Such a decomposition
distinguishes a line bundle in the algebraic class of~$\check{\cL}$,
namely the symmetric line bundle $\check{\cL}_{0}=L(H,\chi_{0})$ of
characteristic 0 (see Section~\ref{subsec:cantheta}) associated to the 
semicharacter $\chi_{0}(v_{1}+v_{2})=e(\pi\i\Im H(v_{1},v_{2}))$.
\par
\begin{lemma}\label{lem:characteristic}
The line bundles $\cL_{X}$ and $\check{\cL}_{0}$ are algebraically equivalent.
\end{lemma}
\par
\begin{proof} We use the endomorphism $\delta(C,D)$ associated
with a curve $C$ and a divisor $D$ of an abelian variety $T$. It is
defined by mapping $a \in T$ to the sum of the intersection points of the curve
$C$ translated by $a$ and the divisor $D$, see \cite[Section~5.4 and~11.6 ]
{birkenhakelange}. By \cite[Theorem~11.6.4]{birkenhakelange} we need to
show that $\delta(\varphi(X),\check{\cL}) =  \delta(\check{\cL},\check{\cL})$.
By \cite[Proposition~5.4.7]{birkenhakelange} and Riemann-Roch
$\delta(\check{\cL},\check{\cL}) = - 6{\rm id}_{\Prym^{\vee}\!X}$. On the
other hand
$$ \delta(\varphi(X),\check{\cL}) \= - \check{\iota} \circ \iota
\circ \psi_{\check{\cL}}
\= - \phi_{\cL} \circ \psi_{\check{\cL}}\= -6 \,{\rm id}_{\Prym^{\vee}\!X}\,$$
by \cite[Proposition~11.6.1]{birkenhakelange}.
\end{proof}

\subsection{The pre-Abel-Prym map}

Next, we study the pre-Abel-Prym map. We write $\Prym^{\vee}\!X = V / \Lambda$.
\par
\begin{lemma}\label{lem:preAbelPrym} The pre-Abel-Prym map $\pabel$ with base
point $p_1$ sends the $p_i$ to zero, i.e.\ 
	\bes
	\pabel(p_{1})\= \pabel(p_{2})\= \pabel(p_{3})\= 0\,\quad\text{for $i=1,2,3$.}
	\ees
The points in $Z$ are sent to three different non-trivial two-torsion points
in a Lagrangian subspace of $\Lambda$, i.e.\ 
$\pabel(Z)=\{\tfrac12\lambda_{2},\tfrac12\mu_{2},\tfrac12(\lambda_{2}+\mu_{2})\}$
for some decomposition of~$V$. 
\par
Moreover, the endomorphism $(-1)$ of $\Prym^{\vee}\!X$ induces the 
involution~$J$ on $\pabel(X)$ and $\pabel$ is injective on $X \setminus P$. 
\end{lemma}
\par
\begin{proof}The inclusion $\A^{\vee}\times \B^{\vee}\subset \Jac X$ is given
in terms of degree-zero divisors~$D$ and~$E$ by $(D,E)\mapsto D+J(D)+E^{(1)}+E^{(2)}+E^{(3)}$ .
In particular, on the images of $q_i$ in~$A$ and~$B$ (as in~\eqref{eq:hpreim}) 
this inclusion map is given by $(\overline{q}_{1}-\overline{q}_{i},e'_{i}-e'_{1})
\mapsto [p_{i}-p_{1}]=\pabel(p_{i})$. This proves that the points $p_{i}$ are sent to zero.
\par
Next, for each $x\in X$ the divisor $x+J(x)-2p_{1}$ belongs to 
$\pi_A^* {\rm Div}^0(A)$, hence maps to zero in~$\Prym^{\vee}\!X$ 
and therefore
	\[\pabel(x)\=[x-p_{1}]=[-J(x)+p_{1}]=-\pabel(J(x))\,.\]
In particular the points $\pabel(z_{i})=[z_{i}-p_{1}]$ have order two and
	\[\sum_{i=1}^3 \pabel(z_{i})\=[z_{1}+z_{2}+z_{3}-3p_{1}]=0\,.\]
As a consequence all three of the~$z_i$ are 2-torsion points and 
by Lemma~\ref{lem:AbelPrym} they moreover lie in $\Lambda(H)$. It remains
to exclude that $\pabel(Z) = 0$.
\par
By the preceding Lemma~\ref{lem:characteristic} and Riemann-Roch, the
curve $\pabel(X)$ is of arithmetic genus~$7$. If $\pabel(X)$ is generically injective
then $\pabel(Z)=0$ would imply that there are $6$ branches passing through
zero and the arithmetic genus had to be larger than~$7$, contradiction. 
On the other hand, the geometric genus of $\pabel(X)$ is at least two, since 
this curve generates $\Prym^{\vee}\!X$, hence the degree of $\pabel$ is 
at most three. In this case, the differential $\omega$ has to be a
pullback of a differential on (the normalization of) the genus two 
curve $\Prym^{\vee}\!X$. This is impossible, as discussed in \cite[Lemma~6.2]{MMW}.
\end{proof}
\par
We can now complete the identification of the line bundles begun
in Lemma~\ref{lem:characteristic}.
\par
\begin{lemma}\label{lem:characteristic2} 
Let $(X,\omega)\in \Omega G$. With the above choice of 
a decomposition of~$V$, the line bundles~$\cL_{X}$ and~$\check{\cL}_{0}$ 
differ by the characteristic $c=\tfrac12\lambda_{1}+\tfrac12\mu_{1}$, i.e.\ 
$\cL_{X}=t_{c}^{*}\check{\cL}_{0}$.
\end{lemma}
\par
\begin{proof} 
To compute the characteristic, note that by Lemma~\ref{lem:preAbelPrym} the image $\pabel(X)$ is a symmetric divisor, that is $(-1)^{*}\pabel(X)=\pabel(X)$. The requirement on a line bundle in the algebraic class of $\check{\cL}_{0}$ to be symmetric, narrows the number of choices down to $2^{4}$ possibilities, which agree with the translates of $\check{\cL}_{0}$ by half-integral points. As a consequence $\cL_{X}=t_{c}^{*}\check{\cL}_{0}$ for some half-integral character $c\in\tfrac12\Lambda(H)/\Lambda(H)$.
\par
In order to compute explicitly the characteristic of $\cL_{X}$, let us first note that by Lemma~\ref{lem:preAbelPrym} the only 2-torsion points in $\pabel(X)$ are $\pabel(z_{i})$ for $i=1,2,3$, all of them with multiplicity 1. By~\cite[Proposition~4.7.2]{birkenhakelange} the semicharacter $\chi$ associated to the line bundle $\cL_{X}$ takes the value
	\[\chi(\lambda)=(-1)^{\mult_{\frac{1}{2}\lambda}(\pabel(X))-\mult_{0}(\pabel(X))}\]
for each lattice element $\lambda\in\Lambda$. Since $\mult_{0}(\pabel(X))=3$, we deduce that $\chi$ takes values $+1$ at $\lambda_{2}$, $\mu_{2}$, $\lambda_{2}+\mu_{2}$ and $-1$ at $\lambda_{1}$, $\mu_{1}$, $\lambda_{1}+\mu_{1}$.
\par
Recall that $\Lambda(H)=\langle\lambda_{1},\tfrac{1}{6}\lambda_{2},\mu_{1},\tfrac{1}{6}\mu_{2}\rangle$, and let $c=a_{1}\lambda_{1}+\tfrac{a_{2}}{6}\lambda_{2}+b_{1}\mu_{1}+\tfrac{b_{2}}{6}\mu_{2}$, where $a_{1},a_{2},b_{1},b_{2}\in\{0,\tfrac{1}{2}\}$. Using the fact that $\chi=\chi_{0}\cdot \exp(2\pi\i\Im H(c,\cdot))$ and the expression~\eqref{eq:semicharacter} for $\chi_{0}$, one gets $a_{1}=b_{1}=\frac{1}{2}$ and $a_{2},b_{2}=0$.
\end{proof}

\subsection{Identifying the theta function} \label{sec:IdThety}

Our main objective now is to describe $\pabel(X)$ as the vanishing locus of
some theta function $\theta_{X}$ in $H^{0}(\cL_{X})$. For this purpose, we 
restrict furthermore to the case that $(X,\omega)$ is a Gothic eigenform
for real multiplication by $\cO_D$. This implies that on the Prym
variety we have the distinguished eigenform coordinates introduced
in Section~\ref{sec:SME}. 
\par
\begin{lemma}\label{lem:thetaconditions}Let $(X,\omega)\in \Omega G_{D}$ for some $D$. Then $\pabel(X)$ is the vanishing locus of a global section $\theta_{X}\in H^{0}(\cL_{X})_{-}$ satisfying
	\begin{itemize}
	\item[(C1)] $D_{1}\theta_{X}(0) = 0$,
	\item[(C2)] $D_{2}\theta_{X}(0) = 0$,
	\item[(C3)] $D_{2}\theta_{X}(\tfrac{1}{2}\mu_{2}) = 0$,
	\item[(C4)] $D_{2}\theta_{X}(\tfrac{1}{2}\lambda_{2}) = 0$,
	\item[(C5)] $D_{2}\theta_{X}(\tfrac{1}{2}(\lambda_{2}+\mu_{2})) = 0$.
	\end{itemize}
\end{lemma}
\par
\begin{proof}By definition and Lemma~\ref{lem:characteristic2}, $\pabel(X)$ is 
the vanishing locus of some theta function $\theta_{X}\in H^{0}(\cL_{X})$. Since $\mult_{0}(\pabel(X))=3$, this theta function is necessarily odd by~\cite[Lemma 4.7.1]{birkenhakelange} and the comments after that lemma.
\par
Being $\theta_{X}$ an odd function, both $\theta_{X}$ and its second derivatives vanish at 0. 
Since $\mult_{0}(\pabel(X))=3$, also its first derivatives must vanish, that 
is $D_{1}\theta_{X}(0)=D_{2}\theta_{X}(0)=0$.
\par
Let us assume that $du_{1}$ is the eigenform in $\Omega\cM_{4}(2^{3},0^{3})$. Note that the condition of this eigenform having a zero of order $k$ at a point $p$ translates into $\partial^{j}\theta_{X}/\partial u_{2}^{j}$ vanishing at $\pabel(p)$ for $j=0,\ldots,k$.
\end{proof}
\par
Recall the definition of the generators $\theta_{0},\theta_{1},
\theta_{2},\theta_{3}$ of $H^{0}(\cL_{X})_{-}$ from~\eqref{eq:thetaidef},
and let $\theta_{X}(\bbu)=\sum_{i}a_{i}\,\theta_{i}(\bbu)$ be a theta function
cutting out $\varphi(X)$. By~\eqref{eq:derivativesD2}, the conditions
in Lemma~\ref{lem:thetaconditions} correspond to the following system of
equations
	\bes
	\begin{array}{ll}
	\mbox{\emph{(C1)}} & a_{0}\,D_{1}\theta_{0}(0) + a_{1}\,D_{1}\theta_{1}(0) + a_{2}\,D_{1}\theta_{2}(0) + a_{3}\,D_{1}\theta_{3}(0) = 0 \,, \\
	\mbox{\emph{(C2)}} & a_{0}\,D_{2}\theta_{0}(0) + a_{1}\,D_{2}\theta_{1}(0) + a_{2}\,D_{2}\theta_{2}(0) + a_{3}\,D_{2}\theta_{3}(0) = 0 \,, \\
	\mbox{\emph{(C3)}} & a_{0}\,D_{2}\theta_{0}(0) - a_{1}\,D_{2}\theta_{1}(0) + a_{2}\,D_{2}\theta_{2}(0) - a_{3}\,D_{2}\theta_{3}(0) = 0 \,, \\
	\mbox{\emph{(C4)}} & a_{0}\,D_{2}\theta_{3}(0) + a_{1}\,D_{2}\theta_{2}(0) + a_{2}\,D_{2}\theta_{1}(0) + a_{3}\,D_{2}\theta_{0}(0) = 0 \,, \\
	\mbox{\emph{(C5)}} & a_{0}\,D_{2}\theta_{3}(0) - a_{1}\,D_{2}\theta_{2}(0) + a_{2}\,D_{2}\theta_{1}(0) - a_{3}\,D_{2}\theta_{0}(0) = 0 \,. \\
	\end{array}
	\ees
Note that conditions \emph{(C2)}--\emph{(C3)} and conditions
\emph{(C4)}--\emph{(C5)} can be rephrased as
	\be\label{eq:C2toC5}
	\left\{\begin{array}{l}
	a_{0}\,D_{2}\theta_{0}(0) + a_{2}\,D_{2}\theta_{2}(0) = 0 \\
	a_{1}\,D_{2}\theta_{1}(0) + a_{3}\,D_{2}\theta_{3}(0) = 0
	\end{array}\right.\quad\mbox{and}\quad
	\left\{\begin{array}{l}
	a_{0}\,D_{2}\theta_{3}(0) + a_{2}\,D_{2}\theta_{1}(0) = 0 \\
	a_{1}\,D_{2}\theta_{2}(0) + a_{3}\,D_{2}\theta_{0}(0) = 0
	\end{array}\right.\,,
	\ee
respectively. This already allows us to get some necessary conditions on
the derivatives of theta functions for a point to belong to the Gothic locus.
\par
\begin{prop}\label{prop:GDinmodular}
  If the point $\bm\tau \in \bH^2$ has the property that there is a
  non-zero odd theta function
$\theta_{X}(\bbu)=\sum_{i}a_{i}(\bm\tau) \,\theta_{i}(\bbu)$ on $T_{\bm\tau}$
satisfying \emph{(C2)}-\emph{(C5)}, then $\cG_{D}(\bm\tau)=0$.
In particular, for any $(X,\omega)\in \Omega G_{D}$, the Prym variety
$\Prym^{\vee}\!X$ belongs to the vanishing locus of
the Gothic modular form $\cG_{D}(\bm\tau)$.
\end{prop}
\par
\begin{proof}By~\eqref{eq:C2toC5}, the coefficients must satisfy
	\bes
	M \begin{pmatrix}
	a_{0} \\ a_{1} \\ a_{2} \\ a_{3}
	\end{pmatrix}
	\coloneqq
	\begin{pmatrix}
	D_{2}\theta_{0}(0) & 0 &  D_{2}\theta_{2}(0)  & 0 \\
	0 & D_{2}\theta_{1}(0) & 0  &   D_{2}\theta_{3}(0) \\
	D_{2}\theta_{3}(0) & 0 &  D_{2}\theta_{1}(0) & 0  \\
	0 & D_{2}\theta_{2}(0) & 0 & D_{2}\theta_{0}(0)
	\end{pmatrix}
	\begin{pmatrix}
	a_{0} \\ a_{1} \\ a_{2} \\ a_{3}
	\end{pmatrix} = 
	\begin{pmatrix}
	0 \\ 0 \\ 0 \\ 0
	\end{pmatrix}\,.
	\ees
This system of equations must have a non-trivial solution, and therefore
	\[\det(M)=\left(D_{2}\theta_{0}(0)\cdot D_{2}\theta_{1}(0)-D_{2}\theta_{2}(0)\cdot D_{2}\theta_{3}(0)\right)^{2}=\cG_{D}(\bm\tau)^{2}=0\,.\]
The second claim follows from Lemma~\ref{lem:thetaconditions}.
\end{proof}
\par

\subsection{The vanishing locus of the Gothic modular form}
\label{sec:vanGothic}

We now start in the converse direction and analyze the  vanishing locus
of the Gothic modular form $\cG_{D}$. For this purpose we note
that the theta function $\theta_{X}$ defined in Theorem~\ref{thm:GDmodform} equals
\bes
\theta_{X}(\bbu)  \= \Theta_a \cF_b - \Theta_b \cF_a\,,
\ees
where
\bas
\Theta_{a}(\bm{\tau},\bbu) &\=
\theta_{0}(\bm\tau,\bbu)\cdot D_2\theta_{2}(\bm{\tau},0)
- \theta_{2}(\bm\tau,\bbu)\cdot D_2\theta_{0}(\bm{\tau},0)\,, \\
\Theta_b (\bm{\tau},\bbu) &\=
\theta_{1}(\bm\tau,\bbu)\cdot D_2\theta_{3}(\bm{\tau},0)
- \theta_{3}(\bm\tau,\bbu)\cdot D_2 \theta_{1}(\bm{\tau},0)\,
\eas
and where $\cF_a(\bm{\tau}) = D_1 \Theta_a|_{\bbu=0}$ and
$\cF_b(\bm{\tau}) = D_1 \Theta_b|_{\bbu=0}$ as defined in
Theorem~\ref{thm:GDmodform}, too. 
\par
\begin{proof}[Proof of Theorem~\ref{thm:GDmodform}]
The non-vanishing of $\theta_X$ on the complement of $\tRED$
follows from the factorization given above and the linear
independence of the $\theta_i$. Given Proposition~\ref{prop:GDinmodular}
it remains to show that on the complement of $\tRED$ the divisor
$Y = Y_{\bm{\tau}} :=\{\theta_X = 0\}$ is indeed the $\varphi$-image of a Gothic
Veech surface.
\par
We first check the conditions~$(C1)$--$(C5)$ for~$Y$.
Differentiating $\theta_X$ implies that~$Y$ satisfies $(C1)$ using the second
row of the defining matrix, and $Y$ satisfies $(C2)$ and $(C3)$ in the
reformulation~\eqref{eq:C2toC5}, as can be seen from the last two rows.
From~\eqref{eq:derivativesD2} we deduce
	\bas
	D_{2}\theta_{X}\left(\tfrac{1}{2}\lambda_{2}\right) &\= 
	\left|\begin{array}{cccc}
	D_{2}\theta_{3}(0) & D_{2}\theta_{2}(0) & D_{2}\theta_{1}(0) & D_{2}\theta_{0}(0) \\
	D_{1}\theta_{0}(0) & D_{1}\theta_{1}(0) & D_{1}\theta_{2}(0) & D_{1}\theta_{3}(0) \\
	D_{2}\theta_{0}(0) & 0 &  D_{2}\theta_{2}(0) & 0 \\
	0 & D_{2}\theta_{1}(0) & 0 & D_{2}\theta_{3}(0) \\
	\end{array}\right| = \left(\mathcal{F}_{b}+\mathcal{F}_{a}\right)\mathcal{G}_{D}\,,\\
	D_{2}\theta_{X}\left(\tfrac{\lambda_{2}+\mu_{2}}2\right) &\= \left|\begin{array}{cccc}
	- D_{2}\theta_{3}(0) & D_{2}\theta_{2}(0) & -D_{1}\theta_{3}(0) & D_{2}\theta_{0}(0) \\
	D_{1}\theta_{0}(0) & D_{1}\theta_{1}(0) & D_{1}\theta_{2}(0) & D_{1}\theta_{3}(0) \\
	D_{2}\theta_{0}(0) & 0 & D_{2}\theta_{2}(0)  & 0 \\
	0 & D_{2}\theta_{1}(0) & 0  & D_{2}\theta_{3}(0) \\
	\end{array}\right| =\left(\mathcal{F}_{b}-\mathcal{F}_{a}\right)\mathcal{G}_{D}\,.
	\eas
we deduce that for $Y$ the conditions \emph{(C4)}-\emph{(C5)} hold as well.
\par
Since $\theta_X$ is a section of the line bundle $\cL_X$ of characteristic
$c = \tfrac12\lambda_1 + \tfrac12\mu_1$, the argument in Lemma~\eqref{lem:characteristic2} implies that the multiplicity of~$Y_{\bm\tau}$ at
each point in the set $S^+ = \{\tfrac12 \lambda_2,\tfrac12  \mu_2,
\tfrac12(\lambda_2+\mu_2)\}$ is odd, in particular $Y$ contains these points.
By \cite[Proposition~4.7.5~a)]{birkenhakelange} each of the other $2$-torsion
points is not contained in~$Y$ or~$Y$ has even multiplicity there.
\par
{\em The case~$Y$ reduced with zero as its only singular point.} 
Then $X_{\bm\tau} = Y_{\bm\tau}^{\rm norm}$ is the desingularization at zero.
We check the properties of a Gothic eigenform.
Since~$Y_{\bm\tau}$ is non-singular at~$S^+$ in the case under consideration,
the one-form $du_1$ is an eigenform for real multiplication and has double
zero at each of the three the points in~$S^+$.
\par
The involution $(-1)$ on $T_{\bm{\tau}}$ induces an involution~$J$ on~$X_{\bm\tau}$
that has $6$ fixed points. The quotient $A = A_{\bm{\tau}} = X_{\bm\tau}/J$
is therefore a smooth curve of genus one. The complement~$T^\perp_{\bm\tau}$
of $T_{\bm\tau}$ in $\Jac(X)$ is $(1,6)$-polarized (see
\cite[Corollary~12.1.5]{birkenhakelange} or the proof of
Proposition~\ref{prop:Prymis16}). The pullback of the theta divisor
on $\Jac(X)$ to~$A^\vee$ has degree two since $X_{\bm\tau} \to A_{\bm{\tau}}$
is a double covering. We claim that the restriction of the theta divisor
on~$\Jac(X)$ to the complement $B^\vee$ of~$A^\vee$ in~$T^\perp_{\bm\tau}$
has degree three. In fact, we may view $B^\vee$ as the complement
to the image of the addition map $+: A^\vee \times T^\perp_{\bm\tau} \to \Jac(X)$.
If $+$ factors through an isogeny (necessarily of exponent two), then
the image is $(1,1,3)$-polarized, hence $B^\vee$ has a polarization
of type~$(3)$, again by \cite[Corollary~12.1.5]{birkenhakelange}.
The case that $+$ is injective, hence the image is $(1,2,6)$-polarized,
contradicts loc.\ cit. Consequently, the map $\piB: X \to \Jac(X) \to B$ is
a covering of degree three.
\par
We claim that the map $\piB:X_{\bm\tau} \to B_{\bm\tau}$ is odd. In fact, 
writing $j=(-1)$ on the elliptic curve~$B_{\bm\tau}$ we compute that
\[j\circ\piB(x)=[p_{1}-x]=[p_{1}-x]+[J(x)+x-2p_{1}]=[J(x)-p_{1}]=\piB(J(x))\,,\]
since $x+J(x)-2p_{1}\in A^{\vee}$. This argument also shows that the images of
the points in~$P$ and~$Z$ are $2$-torsion points in any quotient of 
$\Jac(X_{\bm\tau})/(A_{\bm\tau}^{\vee})$, in particular in~$B_{\bm\tau}$. Since
$|\piB(Z)|=1$ on points in the Gothic locus, we deduce that $|\piB(Z)|=1$ over
all of~$X_D$. We have indeed checked that $(X,du_1,\pi_A, \pi_B)$
has all the Gothic properties, under our assumptions on~$Y$.
\par
{\em The case~$Y$ reduced with other singularities besides zero}
does not occur. In fact, if $Y = \sum Y_i$ then $Y^2 = 12$ for a
line bundle of type $(1,6)$ by Riemann-Roch. A triple point such as zero
contributes~$6$ to $Y^2$. Each of the points in~$S^+$ is either a triple
point or $du_1$ has a double zero there, contributing~$2$ to~$Y^2$
by increasing the genus of the component passing through this point.
The total count implies that~$Y$ is non-singular at~$S^+$ and also
non-singular elsewhere besides zero, since the three double zeros at~$S^+$
and the contribution at $0 \in T_{\bf{\tau}}$ already add up to~$12$.
\par
{\em The case~$Y$ non-reduced.} 
The above counting argument has to be refined for~$Y$ non-reduced,
since e.g.\ a triple point might consist of~$2Y_1$ and~$Y_2$ intersecting
transversally, hence contributing only~$4$ to~$Y^2$. We first note that
there are at most two branches through zero, since if~$Y$ contained
non-reduced $a_1Y_1 + a_2Y_2 + a_3Y_3$ all meeting at zero, the odd multiplicity at the origin implies that $a_{1}+a_{2}+a_{3}$ is at least 5, and therefore $Y^{2}>12$. 
\par
We now write $Y = a_1Y_1 + a_2Y_2 + Y_R$
with $a_1 \geq a_2$, with $Y_1$ and $Y_2$ irreducible and passing through zero
while $Y_R$ is potentially reducible with no component passing through zero.
In particular $a_1 + a_2$ is odd. 
\par
{\em Case $(a_1,a_2)=(3,2)$.} In this case~$0$ is the only intersection
point of $Y_1$ and $Y_2$ and $Y_{i}^{2}=0$, so both components are elliptic curves. Consider the product $Y_{1}\times Y_{2}$ with the polarization $2\,p_{1}^{*}\cO_{Y_{1}}(0)\otimes 3\,p_{2}^{*}\cO_{Y_{2}}(0)$. The addition map $Y_{1}\times Y_{2} \to T_{\bm\tau}$ is an isomorphism at the level of complex tori since $Y_{1}\cdot Y_{2}=1$, and the pullback of $\cL_{X}=\cO_{T_{\bm,\tau}}(Y)$ agrees with the $(2,3)$-product polarization. In particular the map is an isomorphism of abelian surfaces and hence we are in~$\tRED$ (see Proposition~\ref{prop:RedisRed}).
\par
{\em Case $(a_1,a_2)=(6,1)$.} Again~$0$ is the only intersection point, and $Y_{i}^{2}=0$, so both components are again elliptic curves. Odd parity of the theta function implies that $S^{+}\subset Y_{2}$, but then $du_{1}$ induces an abelian differential on $Y_{2}$ with 3 zeroes of order $\ge 2$, which is a contradiction.
\par
{\em Case $(a_1,a_2)=(4,1)$.} Again~$0$ is the only intersection point, and the case $Y_{i}^{2}=0$ for $i=1,2$ yields the same contradiction as in the case before. Hence we have $Y_{2}^{2}=4$ and $S^{+}\subset Y_{2}$. This implies that on the one hand $Y_{2}$ has genus 3, and on the other $du_{1}$ induces an abelian differential on $Y_{2}$ with 3 zeroes of order $\ge 2$, which is again a contradiction.
\par
{\em Case $(a_1,a_2)=(2,1)$.} We have the following possibilities:
	\begin{enumerate}[nosep]
	\item $Y_{1}\cdot Y_{2} = 1\,,\ Y_{2}^{2}=4$: the curve $Y_{2}$ has genus 3 and $du_{1}$ induces an abelian differential on it with 3 zeroes of order $\ge 2$.
	\item $Y_{1}\cdot Y_{2} = 2\,,\ Y_{2}^{2}=0\mbox{ or }2$: the curve $Y_{2}$ has genus 1 (or 2) and $du_{1}$ induces an abelian differential on it with 2 zeroes of order $\ge 2$.
	\item $Y_{1}\cdot Y_{2} = 3\,,\ Y_{2}^{2}=0$: the curve $Y_{2}$ has genus 1 and $du_{1}$ induces an abelian differential on it with a zero of order $\ge 2$.
	\end{enumerate}
All these cases yield contradictions with the genus of the curve $Y_{2}$ and this completes the claim.
\end{proof}

%% file: sec_modularcurves.tex
\section{Modular curves and the reducible locus} \label{sec:modularcurves}

The main result in this section is an explicit parametrization of the
{\em reducible locus}, the locus where the $(1,6)$-polarized abelian varieties
with real multiplication split as a product of two elliptic curves $E_1$
and $E_2$, which are necessarily isogenous.
This locus is a union of modular curves (also known as Hirzebruch-Zagier
cycles or Shimura curves), in fact exclusively non-compact modular curves.
\par
There are interesting similarities and differences to the reducible locus in
the principally polarized case and the well-studied case of genus two
\Teichmuller\ curves. The main similarity is that the \Teichmuller\ curves are
disjoint from the reducible locus in both situations, Gothic and genus two.
The two cases also agree in the fact that the reducible locus has many components, several but not
all of which can be distinguished by the precise endomorphism ring.
\par
The main difference starts with the fact that the reducible locus
decomposes into two sub-loci that can already be distinguished by degree
of restriction of the polarization line bundle to $E_1$ and $E_2$. Since the
product of these degrees is~$6$, the reducible locus decomposes into
$\RED$ and $\REDIRR$, where the indices give the degree of the restricted line
bundles. These loci are indeed disjoint, as we show in
Section~\ref{sec:Red16}. The main result of this section is a description
of the components of~$\RED$ and a computation of their volumes.
\par

\subsection{Modular curves on Hilbert modular surfaces}
\label{subsec:modularcurves}

The reducible locus consists of modular curves (also known as
Hirzebruch-Zagier cycles or Shimura curves). Modular curves are
the images of graphs of M\"obius transformations in $\bH^2$ that descend
to algebraic curves in the Hilbert modular surface. We recall the
precise definition, adapted to our Hilbert modular surfaces
$X_{D}^{\fraka}(\frakb) = \SL(\fraka\frakb\oplus\fraka^{\vee})\backslash \bH^2$.
\par
Let us define the ideal $\avdg = \sqrt{D}\fraka^{2}\frakb$. We say that $U\in\SL_{2}(K)$ is a \emph{generator matrix} for the Hilbert
modular group $\SL(\fraka\frakb\oplus\fraka^{\vee})$,  if it is of the form
	\bes
	U=\begin{pmatrix} a\sqrt{D} & \mu \\ -\mu^{\sigma} & Ab\sqrt{D} \end{pmatrix}\,,\quad
	\mbox{where $a,b\in\bZ$, $\mu\in \avdg$ and $A=N(\avdg)$.}
	\ees
and we define the modular curves  $F_{U}$ to be the image
in $X_{D}^{\fraka}(\frakb)$ of the set
\[\left\{(\tau_{1},\tau_{2})\in\bH^{2}\,:\, 
\begin{pmatrix} \tau_{2} & 1 \end{pmatrix} U
\begin{pmatrix} \tau_{1} \\ 1 \end{pmatrix} 
\= a\sqrt{D}\tau_{1}\tau_{2} - \mu^{\sigma}\tau_{1} + \mu\tau_{2} + Ab\sqrt{D}
\= 0 \,\right\}\,.\]
The generator matrix $U$ is \emph{primitive} if it is not divisible by any
natural number $m>1$. For any integer $N>0$, the \emph{modular curve}~$F_{N}$
is defined as the union
\[F_{N} \= \bigcup_{\substack{U\ \mbox{\scriptsize primitive} \\ \det(U)=AN}} F_{U} \,.\]
The components of $F_N$ and their geometry (cusps, fixed points) were
intensely studied by Hirzebruch and his students (see the survey in
\cite[Chapter~V]{vdG}). Most notably the volumes of the union
\bes
T_{N} \= \bigcup_{\det(U)=AN} F_{U}  = \bigcup_{\ell^{2}|N} F_{N/\ell^{2}}\,
\ees
are the coefficients of a modular form, in fact an Eisenstein series
of weight~$2$ for some character.
\par
This however does not yet yield formulas for the volume of $\RED$, since
the latter turns out to be a union of modular curves, but not of the entire
curves $F_N$. In fact, $F_N$ can be decomposed as the union of the
curves $F_{N}(\nu)$ for $\nu\in M / \sqrt{D}M$,
where
\[F_{N}(\nu) \= \bigcup\{F_{U} :\, \mbox{$U$ is primitive with $\det(U)=AN$
and $\nu(U)=\nu$} \}\,.\]
In the case of abelian surfaces with principal polarization the reducible
locus was written in terms of $F_N(\nu)$ by~\cite{mcmHMS}. However,
the $F_N(\nu)$ are sometimes still reducible and this decomposition
does not directly yield a volume formula, so we proceed differently
for our $(1,6)$-polarization.
\par 

\subsection{The $(2,3)$-reducible locus}

Let us define the \emph{$(2,3)$-reducible locus~$\RED$} as the locus inside the moduli space $\cA_{2,(2,3)}$ of $(2,3)$-polarized abelian surfaces consisting of products $E_{1}\times E_{2}$ of elliptic curves with the natural
$(2,3)$-polarization $2\,p_{1}^{*}\cO_{E_{1}}(0)\otimes 3\,p_{2}^{*}\cO_{E_{2}}(0)$. For each $\cO_{D}$-ideal $\frakb$ of norm 6, we will write $\RED(\frakb)$ for the pullback of $\RED$ to $X_{D}(\frakb)$.
\par
\begin{theorem}\label{thm:ProtToHMS} Let $D=f^{2}D_{0}$ be a positive
quadratic discriminant with conductor~$f$. There is a bijective
correspondence between irreducible components of the $(2,3)$-reducible locus
with a chosen proper real multiplication by $\cO_{D}$ and the set of
prototypes
\bes
\cP_{D} \= \left\{[\ell,e,m]\in\bZ^{3}\, :\,
\mbox{$\ell,m>0$, $D=e^{2}+24\ell^{2}m$ and $\gcd\left(e,\ell,f \right)=1$}\right\}\,.
\ees
More precisely, the component parametrized by the prototype $\cP=[\ell,e,m]$ is
the image of a Shimura curve in the Hilbert modular surface
$X^{\fraka}_{D}(\frakb)$ corresponding to the ideals
$\fraka=\tfrac{1}{\sqrt{D}}(2\ell,\frac{e+\sqrt{D}}{2})$ and
$\frakb=(6,\frac{r+\sqrt{D}}{2})$, where
	\[r=\left\{\begin{array}{ll}
	e\,, & \mbox{if $D\equiv 0 \bmod{2}$,} \\
	e+6\,, & \mbox{if $D\equiv 1 \bmod{2}$.}
	\end{array}\right.\]
The image in $\cA_{2,(2,3)}$ of the Shimura curve given by $\cP=[\ell,e,m]$
is isomorphic to~$\Gamma_0(m)\backslash\bH$.
\end{theorem}
\par
We split the proof into a series of lemmas.
\par
\begin{lemma} \label{le:23stdper} 
The period matrix of an abelian surface parametrizing
a point in $\RED$ with real multiplication can be assumed to be 
	\[\Pi_{m}(\tau) = 
	\left(\begin{matrix} 2\tau & 0 \\ 0 & 3m\tau \end{matrix}
	\ \middle\vert \ 
	\begin{matrix} 2 & 0 \\ 0 & 3 \end{matrix} \right)
	,\quad\mbox{for some $\tau\in\bH$ and $0<m\in\bZ$}
	\]
with the polarization given by standard form $\sm{0}{P_{23}}{-P_{23}}0$.
\end{lemma}
\par
\begin{proof}
Since we will be interested in the components of this locus that lie in some
Hilbert modular surface, let us assume furthermore that $E_{1}$ and $E_{2}$
are isogenous elliptic curves, so the left block of of the period matrix
$\Pi_{m,n}(\tau)$ is a diagonal matrix with entries $(2\tau, 3(m\tau +n))$
with $m,n \in \bQ$. Positive definiteness of the period matrix implies
$m>0$. We define the matrices
\bes
M_1 = \left(\begin{matrix}
a& 0 & b  &0 \\
0 & aU+cV & 0 & \tfrac{dV +b}{L} \\
b & 0 & dU  & 0 \\
0 & cL & 0  & d
\end{matrix}\right)
\quad \text{and} \quad
M_2 = \left(\begin{matrix}
x & 0 & y  &0 \\
0 & xq & 0 & yp \\
p & 0 & q  & 0 \\
0 & 1 & 0  & 1
\end{matrix}\right)\,.
\ees
\par
We first argue that we can take $n=0$.
Write $m=U/L$ and $n=V/L$ with $\gcd(U,V,L)=1$.
Take~$d$ such that $\gcd(dU, L-dV)=1$. (To show the existence,
consider $d_i$ with $\gcd(d_i,L)=1$. Among a collection of $d_i$
with $\gcd(d_i-d_j,U)=1$ with more elements than~$B$ has divisors,
one will work.) Let $b = L-dV$ and take $a,c$ such that
$adU -c(L-dV)=1$. Then the matrix~$M_1$ 
has integral coefficients, belongs to the symplectic group~$\Sp_{2g}^P(\bZ)$
and takes $\Pi_{m,n}(\tau)$ to $\Pi_{m,0}(\tau')$ for some~$\tau'$.
\par
To show that we may assume $m \in \bZ$ we write $m = p/q$
and take $x,y \in \bZ$ such that $xq - yp =1$. Then the matrix~$M_2$
belongs to~$\Sp_{2g}^P(\bZ)$ and takes $\Pi_{p/q,n}(\tau)$ to $\Pi_{pq,0}(\tau')$
for some~$\tau'$.
\end{proof}
\par
\begin{lemma} \label{le:23unique}
An abelian surface in the $(2,3)$-reducible locus
contains a unique elliptic curve with a polarization of type~$(2)$
and a unique elliptic curve with a polarization of type~$(3)$.
\par
In particular a matrix $M \in \Sp_{2g}^P(\bZ)$ taking the locus
$\{ \Pi_m(\tau), \tau \in \bH \}$ into
some locus  $\{ \Pi_{m_2}(\tau), \tau \in \bH \}$ consists of matrices
diagonal in each of its four blocks (like the matrices $M_{1}$ and $M_{2}$ above).
\end{lemma}
\par
\begin{proof} The type of a polarization is translation invariant.
So we may assume that the elliptic curve in question passes through the
origin. Such an elliptic curves~$E$ in a product of elliptic curves
is determined by a rational slope in the universal cover. We
may assume this slope is $(2x,3y)$ with $x,y \in \bZ$ coprime and
both different from zero, since we already know the polarizations
of the curves with slope $(1,0)$ and $(0,1)$. If we denote
by $a_1,a_2,b_1,b_2$ the symplectic basis corresponding to the column
vectors of $\Pi_m(\tau)$, lattice points in~$E$ are given by the
multiples of $f_1 = xa_1 + (y/m)a_2$ and $f_2 = xb_1 + yb_2$
that have integral coefficients. This implies that $m|y$ and
that the type of the polarization on~$E$ is
$\langle f_1, f_2 \rangle = 2x^2 + 3my^2$, therefore proving the claim.
\end{proof}
\par
\begin{lemma} \label{le:RedHasRM}
The analytic representation of real multiplication
by $\gamma=\tfrac{D+\sqrt{D}}{2}$ on an abelian surface with
period matrix $\Pi_{m}(\tau)$ with $m \in \bZ$ is given by
\bes
A_{\gamma}\=\left(\begin{matrix} \frac{D+e}2 &2\ell \\ 3\ell m &
\frac{D-e}2
\end{matrix}\right)\,,
\ees
with $e,\ell \in \bZ$ and $D = e^2 + 24 \ell^{2} m$.
\par
The real multiplication defined by $[\ell,e,m]$ and $[-\ell,e,m]$ are equivalent, whereas the real multiplication defined by $[\ell,e,m]$ and $[-\ell,-e,m]$ are Galois conjugate.
\end{lemma}
\par
\begin{proof}
The abelian surface $T_{\tau,m}$ given by the period matrix $\Pi_{m}(\tau)$
admits real multiplication by $\cO_{D}$, if and only if there are matrices
$A_{\gamma} \in \GL_2(\bQ)$ and $R_{\gamma} \in \Sp(4,\bZ)$ that are
the analytic and rational representations of $\gamma=\tfrac{D+\sqrt{D}}{2}$,
i.e.\ such that $A \Pi_{m}(\tau)  = \Pi_{m}(\tau) R_{\gamma}$,
$\tr(A_{\gamma}) = D$ and $\det(A_{\gamma}) = (D^2-D)/4$. Together with
the self-adjointness of~$R_{\gamma}$ this implies that
	\[
	A_{\gamma} \=\left(\begin{matrix}a&b \\c&d\end{matrix}\right)\,, \qquad 
	R_{\gamma} \=
	\left(\begin{matrix}
	a& 3b m/2 &0 &0 \\
	2c/3m & d & 0&0 \\
	0 & 0& a & 3b/2 \\
	0 & 0& 2c/3 & d
	\end{matrix}\right)\,.	
	\]
such that $d = D-a$ and $ad-bc= (D^2-D)/4$, and moreover that
$c=\tfrac32bm \in 3\bZ$.
Integrality of $R_{\gamma}$ implies that $a,d, \ell = b/2 \in \bZ$
and we set $e = 2a-D$.
\par
Finally, the real multiplications defined by $[\ell,e,m]$ and $[-\ell,e,m]$ are conjugate under the isomorphism $-\mathrm{Id}|_{E_{2}}$. The claim about Galois conjugation is obvious.
\end{proof}
\par
\begin{proof}[Proof of Theorem~\ref{thm:ProtToHMS}]
Suppose we are given a tuple $[\ell,e,m]$ as in the theorem. We check
that the real multiplication on the locus of matrices $\Pi_m(\tau)$
given by Lemma~\ref{le:RedHasRM} is indeed proper. The action is not
proper if $\gamma/k$ also acts for some $1<k \in \bZ$ , i.e.\ if
all the entries of $R_\gamma$ are divisible by~$k$. This implies
$k|\gcd(e,\ell,f)$ and conversely this divisibility is also sufficient
for the action to be non-proper.
\par
Next we show that  the images in $\cA_{2,(2,3)}$ of the loci given
by $\Pi_m(\tau)$ for $m \in \bZ$ are pairwise disjoint. Otherwise there
exists a symplectic matrix taking the locus $\Pi_m(\tau)$ into $\Pi_{m_2}(\tau)$.
By Lemma~\ref{le:23unique} this matrix is diagonal in each block. It suffices
thus consider only matrices of the form
\bes
M = \left(\begin{matrix}
a& 0 & b  &0 \\
0 & ka/m & 0 & kb \\
b & 0 & d  & 0 \\
0 & c/k & 0  & md/k
\end{matrix}\right)
\ees
with integral entries and $ad-bc = 1$, where a priori $k \in \bQ$
and $m_2 = k^2/m$, which implies $k \in \bZ$. Since~$c$ and~$d$ have
no common divisor, this implies $k|m$, hence $k=m=m_2$.
\par
This argument also gives the stabilizer of the locus $\{\Pi_m(\tau), \tau
\in \bH\}$ in the symplectic group. Such a matrix is of the
form of $M$ with $k=m$ and integrality of the entries
implies that the quotient curve is isomorphic to~$\Gamma_0(m)\backslash\bH$.
\par
Finally, we determine for each component in $\RED$ with chosen real
multiplication given by prototype $\cP=[\ell,e,m]$ a Hilbert modular surface
and a Siegel modular embedding that maps to this component. To exhibit a
Siegel modular embedding, we need to find eigenform coordinates, i.e.\
a matrix that diagonalizes the analytic representation of real multiplication given by $A_{\gamma}$ in Lemma~\ref{le:RedHasRM}. Such a matrix is given by
\bes
V_{\cP}\=
\begin{pmatrix} 1 & 1 \\ \lambda
& \lambda^{\sigma} \end{pmatrix},\quad\mbox{where $\lambda\coloneqq\lambda_{\cP}=\frac{-e+\sqrt{D}}{4\ell}$.}
\ees
Indeed, associated to the prototype $\cP$ one can produce the quadratic form $Q_{\cP}=[2\ell,e,-3\ell m]$ of discriminant $D$, so that $\lambda$ is precisely the quadratic irrationality of $Q_{\cP}$, and Lemma~\ref{lem:symplecticallyadapted} ensures that the first column $(1,\lambda)$ of the matrix
$V_{\cP}=B_{\bm\eta}$ is a $(2,3)$-symplectically adapted basis for a
fractional ideal $\fraka^{\vee}$ of $\cO_{D}$, and the first column of the
matrix $(V_{\cP}^{-1}P_{23})^{T}$ is a basis
$\frac{1}{\sqrt{D}}(-4\ell\lambda^{\sigma},6\ell)$ of the ideal $\fraka\frakb$. A
simple calculation shows that
$\fraka=\frac{1}{\sqrt{D}}\left\langle 2\ell,-2\ell\lambda^{\sigma}\right\rangle$ and
therefore, writing $\frakb=\langle 6,\frac{r+\sqrt{D}}{2}\rangle$ for $r\in\bZ$, the following
equality of ideals determines~$r$
	\begin{align*}
	\sqrt{D}\fraka\frakb & = \left\langle 4\ell\lambda^{\sigma},6\ell\right\rangle = \\
		&= \left\langle 2\ell,-2\ell\lambda^{\sigma}\right\rangle \left\langle6,\tfrac{r+\sqrt{D}}{2}\right\rangle = \left\langle 12\ell , -12\ell\lambda^{\sigma}, \ell(r+\sqrt{D}) , \tfrac{(D+er)+\sqrt{D}(e+r)}{4} \right\rangle.
	\end{align*}
To verify this, it is enough to prove that the second ideal
lies in the first one, and one checks that this holds for~$r$ as stated in
the theorem.
\end{proof}
\par
In order to translate the theorem into Euler characteristics,
we define another set of prototypes, closely related to standard
quadratic irrationalities. For a quadratic discriminant $D=f^{2}D_{0}$ with conductor $f$, we let
\begin{align}
\cP_k(D) &\=  \bigl\{ [a,b,c] \in \bZ^{3}\, :\,
\mbox{$a>0>c$ , $D=b^{2} - 4\cdot k \cdot ac$}  \label{eq:defPkD} \\
& \qquad\  \mbox{and $\gcd\left(f,b,c/c_{0}\right)=1$, where $c_{0}$ is the square-free part of $c$ }\bigr\}\,.\nonumber
\end{align}
The following result gives an explicit formula for the Euler characteristics of the reducible loci $\RED(\frakb)$ in terms of prototypes.
\par
\begin{lemma} \label{lem:volRed}The Euler characteristic of the reducible locus
$\RED(\frakb)$ in the Hilbert modular surface $X_{D}(\frakb)$ is given by
\bes \label{eq:chired1}
	\chi(\RED(\frakb)) \=  -\frac{1}{6k} \, \sum_{[a,b,c]\in \cP_6(D) } a \,,
\ees
for each of the ideals $\frakb$ of norm~$6$ in $\cO_{D}$, where $k$ is the number of $\cO_{D}$-ideals of norm 6.
\end{lemma}
\par
\begin{proof} By Theorem~\ref{thm:ProtToHMS}, the different components of $\RED$ in $\cA_{2,(2,3)}$ are isomorphic to certain $\Gamma_0(m)\backslash\bH$. Note that
\bes
\chi(\Gamma_0(m)) \= -\frac{1}{6}\, m\!\! \prod_{\substack{p|m \\ p\,\, \text{prime}}}\!
\bigl(1 + \frac{1}p\bigr)\,.
\ees
Moreover, it is easy to show that for each $D$
\bes
\sum_{[\ell,e,m]\in\cP_{D}} \chi(\Gamma_0(m))
\= -\frac{1}{6} \sum_{[a,b,c]\in\cP_{6}(D)} a\,.
\ees

Let us now suppose that $D \equiv 4,9,16 \bmod 24$, so that there exist two ideals $\frakb\neq\frakb^{\sigma}$ of norm 6. This implies that precomposition
of a chosen real multiplication $\cO_D \to \End_{T_\tau}$ with Galois
conjugation gives a point on a different Hilbert modular surface,
the one with the conjugate $\frakb^\sigma$. Each component of $\cP_D$ is in the image of some
Hilbert modular surface $X_{D}^\fraka(\frakb)$ with $\frakb$ determined
in Theorem~\ref{thm:ProtToHMS} and thus, by the change of cusp
explained in Section~\ref{subsec:cusps}, also on the standard Hilbert modular
surface $X_{D}(\frakb)$. Precomposition with Galois conjugation
corresponds to $e \mapsto -e$. Consequently, on
$X_{D}^\fraka(\frakb)$
\bes
\chi(\RED(\frakb)) \= \frac{1}{2}\sum_{[\ell,e,m]\in \cP_D} \chi(\Gamma_0(m)) = -\frac{1}{2}\frac{1}{6} \sum_{[a,b,c]\in\cP_{6}(D)} a\,.
\ees
\par

In the case $D \equiv 0,12 \bmod 24$ there exists only one ideal $\frakb = \frakb^\sigma$ of norm 6, and thus the map $X_{D}(\frakb) \to \cA_{2,(2,3)}$ is generically
$2:1$ onto its image. In the particular case of $\RED(\frakb)$, components corresponding to prototypes $[\ell,e,m]$ and $[\ell,-e,m]$ are sent to the same component of $\RED$, whereas components corresponding to prototypes $[\ell,0,m]$ lie in the ramification locus of $X_{D}(\frakb) \to \cA_{2,(2,3)}$. As a consequence, 
\bes
\chi(\RED(\frakb)) \= \sum_{[\ell,e,m]\in \cP_D} \chi(\Gamma_0(m)) = -\frac{1}{6} \sum_{[a,b,c]\in\cP_{6}(D)} a\,.
\ees
\par 
Finally, if $D \equiv 1 \bmod 24$ there exist four ideals $\frakb_{1},\frakb_{1}^{\sigma},\frakb_{2},\frakb_{2}^{\sigma}$ of norm 6. For the same reason as above, the forgetful map from~$\RED(\frakb_{i})$ to $\cA_{2,(2,3)}$ is an isomorphism onto its image. Precomposition with Galois conjugation corresponds again 
to $e \mapsto -e$. We conclude
\bes
\chi(\RED(\frakb_{1})) + \chi(\RED(\frakb_{2})) \= \frac{1}{2}\sum_{[\ell,e,m]\in \cP_D} \chi(\Gamma_0(m)) = -\frac{1}{2}\frac{1}{6} \sum_{[a,b,c]\in\cP_{6}(D)} a\,.
\ees
Using Lemma~\ref{le:samevol} we deduce that $\chi(\RED(\frakb_{1})) = \chi(\RED(\frakb_{2}))$ and the result follows.
\end{proof}
\par
\begin{lemma} \label{le:samevol}
For $D\equiv 1 \mod 24$ not a square
\bes
\sum_{b \equiv 1,11 \mod 12 \atop 0 < b < \sqrt{D}}
\sigma_1\Bigl(\frac{D-b^2}{24}\Bigr) \=
\sum_{b \equiv 5,7 \mod 12 \atop 0 < b < \sqrt{D}}
\sigma_1\Bigl(\frac{D-b^2}{24}\Bigr)\,.
\ees
\end{lemma}
\par
\begin{proof} Recall the definition 
\bes
\eta(q) \= q^{1/24}\prod_{n=1}^{\infty}(1-q^{n}) = \sum_{b \geq 1} \Bigl(\frac{12}{b}\Bigr) q^{b^2/24}
\ees
of the Dedekind $\eta$-function and recall that
\bes
E_2(q) \= -\frac1{24} + \sum_{n \geq 1} \sigma_1(n) q^n \=  
\frac{\eta'(q)}{\eta(q)}\,.
\ees
where $' = q \tfrac{\partial}{\partial q}$. The statement of the lemma 
is now equivalent to
\bes
0 \= [q^{D/24}] (E_2(q) \eta(q)) =  [q^{D/24}] \eta'(q)\,,
\ees
which obviously holds for $D$ non-square by definition of~$\eta$.
\end{proof}
\par
Finally, we relate the components given by prototypes at least
coarsely to the usual classification of modular curves. 
\par
\begin{prop}\label{prop:RedModCurve} Let $\cP=[\ell,e,m]\in \cP_{D}$ be a prototype for real multiplication
by $\cO_{D}$ belonging to $X^{\fraka}_{D}(\frakb)$. The corresponding component $F_{\cP}$ of $\RED(\frakb)$ is an irreducible
component of the modular curve $F_{g^{2}m}(\mu)$, where $\mu=g(e+\sqrt{D})/\sqrt{D}$ and $g=\gcd(e,\ell)$.
\end{prop}
\par
\begin{proof} By the proof of Theorem~\ref{thm:ProtToHMS}, we know that the Siegel modular embedding determined by this prototype is given by $V_{\cP}$, and therefore one has
\bes
\begin{pmatrix} 2\tau & 0 \\ 0 & 3m\tau \end{pmatrix}
\=
V_{\cP}\begin{pmatrix} \tau_{1} & 0 \\ 0 & \tau_{2} \end{pmatrix} V_{\cP}^{T} 
\ = \begin{pmatrix} \tau_{1}+\tau_{2} & \lambda\tau_{1}+\lambda^{\sigma}\tau_{2} \\ \lambda\tau_{1} +\lambda^{\sigma}\tau_{2} & \lambda^2\tau_{1} + (\lambda^{\sigma})^{2}\tau_{2}\end{pmatrix}
\ees
for a curve $(\tau_{1},\tau_{2})=(\tau_{1}(\tau),\tau_{2}(\tau))$ in
$\bH^{2}$, where $\lambda\coloneqq \lambda_{\cP}=\frac{-e+\sqrt{D}}{4\ell}$.
\par
In particular, this curve necessarily lies in the curve $\lambda\tau_{1}+\lambda^{\sigma}\tau_{2}=0$, which obviously agrees with $F_{g^{2}m}(\mu)$. The only thing left to prove is that $\mu$ is primitive in $M$ and $N(\mu)=N(M)g^{2}m$.

From the calculations in the proof of Theorem~\ref{thm:ProtToHMS} one gets
	\[M=\sqrt{D}\fraka^{2}\frakb = \frac{4\ell}{\sqrt{D}}\left\langle 3\ell,\ell\lambda^{\sigma},e\lambda^{\sigma}\right\rangle\,,\]
and $N(M)=24\ell^{2}/D$. Since $\mu = \tfrac{4\ell}{\sqrt{D}}g\lambda^{\sigma}$, it is clear that $\mu$ is primitive in $M$, and $N(\mu)= 24\ell^{2}g^{2}m/D$.
\end{proof}
\par

\subsection{The $(1,6)$-reducible locus} \label{sec:Red16}

To put the results of the previous section in perspective we compare
here loci of reducible abelian surfaces according to their polarization.
The moduli space $\cA_{2,(1,6)}$ of $(1,6)$-polarized abelian surfaces
is of course isomorphic to~$\cA_{2,(2,3)}$ used in the previous section,
an isomorphism being induced by multiplication of period matrices by
$\diag(1/2,2)$ from the left. 
\par
In $\cA_{2,(1,6)}$ (and by the above isomorphism thus also in~$\cA_{2,(2,3)}$)
one can similarly define the \emph{$(1,6)$-reducible
locus~$\REDIRR$} of products $E_{1}\times E_{2}$ of isogenous elliptic curves
with the natural $(1,6)$-polarization $p_{1}^{*}\cO_{E_{1}}(0)\otimes 6\,
p_{2}^{*}\cO_{E_{2}}(0)$. With the arguments of Lemma~\ref{le:23stdper} 
we can put period matrices in~$\REDIRR$ in the form
	\[\Pi_{\tau,m} = 
	\left(\begin{matrix} \tau & 0 \\ 0 & m\tau \end{matrix}
	\ \middle\vert \ 
	\begin{matrix} 1 & 0 \\ 0 & 6 \end{matrix} \right)
	,\quad\mbox{for some $\tau\in\bH$ and $0<m\in\bQ$,}
	\]
with the polarization given by the standard form 
$\sm{0}{P_{16}}{-P_{16}}0$. The remaining arguments in the previous
section work verbatim in this case as well and yield:
\par
\begin{theorem}
Let $D=f^{2}D_{0}$ be a positive quadratic discriminant with conductor $f$.
There is a bijective correspondence between irreducible components of
the $(1,6)$-reducible locus admitting proper real multiplication by $\cO_{D}$
and the set of prototypes~$\cP_D$ as defined in Theorem~\ref{thm:ProtToHMS}.
\end{theorem}
\par
In particular, $\RED$ and $\REDIRR$ have the same Euler characteristics.
However:
\par
\begin{prop} The loci $\RED$ and $\REDIRR$ are disjoint in~$\cA_{2,(2,3)}$.
\end{prop}
\par
\begin{proof} The degrees of elliptic curves on an abelian surface
in~$\RED$ are the values of the quadratic form $2x^2+3my^2$
for $x,y, \in \bZ$, as computed in the proof of Lemma~\ref{le:23unique}.
This form never takes the value~$1$.
\end{proof}
\par

%% file: sec_cusps.tex
\section{The divisor of the Gothic modular form}\label{sec:divGD}

In this section we calculate the vanishing locus of the Gothic modular form.
\par
\begin{theorem}\label{thm:divGD} Let $G_{D}(\frakb)$ denote the union of components 
of the Torelli-image of $G_{D}$ lifted to $X_{D}(\frakb)$ such that $du_{1}$ 
induces the eigenform $\omega$ at each point $(X,\omega)$. Then
	\[\div(\cG_{D})=G_{D}(\frakb)+2\RED(\frakb)\,.\]
\end{theorem}

The theorem will be a direct consequence of Propositions~\ref{prop:GothicOnRed} and~\ref{prop:GothicOnGD} below, together with Theorem~\ref{thm:GDmodform}.

\subsection{The Fourier expansion of the Gothic modular form}

For each cusp $\fraka\in X_{D}(\frakb)$ let $\bm\eta=(\eta_1,\eta_2)$ be a basis of $\fraka^\vee$ which is $(2,3)$-symplectically adapted, determining the $\cO_{D}$-module $\fraka\frakb\oplus\fraka^{\vee}$. 
\par
We want to write down the Fourier expansion of $\cG_{D}$ around this cusp using the Siegel modular embedding given by the matrix $B\coloneqq B_{\bm\eta}=\left(\begin{smallmatrix} \eta_1 & \eta_1^\sigma \\\eta_2 & \eta_2^\sigma \end{smallmatrix}\right)$ as in Section~\ref{subsec:cusps}, so that the cusp $\fraka$ of $X_{D}(\frakb)$ corresponds to the cusp at infinity of $X_{D}^{\fraka}(\frakb)$. The stabilizer of $\infty$ agrees with the subgroup 
	\[\SL(\fraka\frakb\oplus\fraka^{\vee})_{\infty}=\left\{\begin{pmatrix}\varepsilon & \mu \\ 0 & \varepsilon^{-1} \end{pmatrix}\,:\, \varepsilon\in\cO_{D}^{*}\,,\ \mu\in M\coloneqq \sqrt{D}\fraka^{2}\frakb\, \right\}\,.\]
For any Hilbert modular form $f$ one has $f(\bm\tau+\mu)=f(\bm\tau)$ for $\mu \in M$, and therefore one can write the Fourier expansion 
	\bes
	f(\bm\tau)=\sum_{\nu\in M^{\vee}}a_{\nu}\bbe\left(\tr(\nu \bm\tau)\right) \,,
	\ees
where $\tr(\nu \bm\tau)=\nu \tau_{1}+\nu^{\sigma}\tau_{2}$ and $M^{\vee}=(\sqrt{D}\fraka^{2}\frakb)^{\vee}=\frac{1}{\sqrt{D}}\fraka^{\vee}(\fraka\frakb)^{-1}$.

Denote by $\rho_{\bm\eta}(\bbx)\coloneqq x_{1}\eta_{1} + x_{2}\eta_{2}$, for $\bbx=(x_{1},x_{2})\in\bQ^{2}$. We will drop the subindex from $\rho_{\bm\eta}$ whenever the choice of basis is clear.
\par

\begin{prop}\label{prop:GothicFourier23} 
The Fourier expansion of $\cG_{D}$ around the cusp~$\fraka$
is given by
	\bas
	\cG_{D}(\bm\tau)  = 8\pi^{2}\i\cdot\,\bigg( &\sum_{\substack{\bba\in\Lambda_{0,\frac{1}{2}}  \\ \bbb\in\Lambda_{\frac{1}{2},\frac{1}{6}}}}\!\!\!  	
	k_{\bba,\bbb} \ 
	q_{1}^{\rho_{\bm\eta}\left(\bba\right)^{2}+\rho_{\bm\eta}\left(\bbb\right)^{2}} 
	q_{2}^{\rho_{\bm\eta}^{\sigma}\left(\bba\right)^{2}+\rho_{\bm\eta}^{\sigma}\left(\bbb\right)^{2}} \\
	- &\sum_{\substack{\bba\in\Lambda_{\frac{1}{2},\frac{1}{2}}  \\ \bbb\in\Lambda_{0,\frac{1}{6}}}}\!\!\!  	
	k_{\bba,\bbb} \ 
	q_{1}^{\rho_{\bm\eta}\left(\bba\right)^{2}+\rho_{\bm\eta}\left(\bbb\right)^{2}} 
	q_{2}^{\rho_{\bm\eta}^{\sigma}\left(\bba\right)^{2}+\rho_{\bm\eta}^{\sigma}\left(\bbb\right)^{2}}\bigg)\,,
	\eas
where $k_{\bba,\bbb}=(-1)^{a_{2}+b_{2}}
	\rho_{\bm\eta}^{\sigma}(\bba)\rho_{\bm\eta}^{\sigma}(\bbb)$ and $\Lambda_{\epsilon,\delta}=\bZ^{2}+(\epsilon,\delta)^{T}$.
\end{prop}

\subsection{Vanishing order along \texorpdfstring{$\RED$}{Red_23}}

The reducible loci $\RED(\frakb)$ turn out to lie in the vanishing locus of the Gothic modular form $\cG_{D}$. We next calculate the corresponding vanishing order.

Recall that, by the results of Section~\ref{subsec:modularcurves}, the reducible loci $\RED(\frakb)$ decompose into different components $F_{\cP}$ indexed by prototypes in $\cP_{D}$.

\begin{prop} \label{prop:GothicOnRed}
The Gothic modular form $\cG_{D}$ vanishes to order 2 along the reducible locus $\RED(\frakb)$.
\end{prop}

\begin{proof}Let $\cP=[\ell,e,m]\in \cP_{D}$ be the prototype corresponding to a component $F_{\cP}\subset F_{g^{2}m}(\mu)$ of $\RED(\frakb)$ as in Proposition~\ref{prop:RedModCurve}. Recall from Theorem~\ref{thm:ProtToHMS} that the curve $F_{\cP}$ lives in the Hilbert modular surface $X_{D}(\frakb)$ determined by the $(2,3)$-symplectically adapted basis $\fraka^{\vee}=\langle 1,\lambda\rangle$, where $\lambda=\frac{-e+\sqrt{D}}{4\ell}$. Note that $\lambda$ is precisely the irrationality associated to the quadratic form $Q_{\cP}=[2\ell,e,-3\ell m]$ of discriminant $D$. Moreover, by Proposition~\ref{prop:RedModCurve} the curve $F_{\cP}$ can be parametrized by $\tau \mapsto (\alpha\tau, \alpha^{\sigma}\tau)$, where $\alpha=-\tfrac{1}{4\ell g}\mu =\frac{\lambda^{\sigma}}{\sqrt{D}}$.

In the chosen basis for $\fraka^{\vee}$, one has $\rho(x_{1},x_{2}) = x_{1}+x_{2}\lambda$ and therefore
	\[\tr(\alpha\rho(x_{1},x_{2})^{2})= \frac{1}{2\ell} \left(x_{1}^{2} + \frac{3}{2}m x_{2}^{2}\right)\,.\]

Now, restricted to $F_{\cP}$ the coordinates $q_{1}$ and $q_{2}$ become $q^{\alpha}$ and $q^{\alpha^{\sigma}}$ respectively, where $q=\bbe(\tau)$. In particular, up to a $8\pi^{2}\i$ factor, the expression for $\cG_{D}$ from Proposition~\ref{prop:GothicFourier23} along $F_{\cP}$ reads
	\begin{align*}
	\cG_{D}(\tau) & =
	\sum_{\substack{\bba\in\Lambda_{0,\frac{1}{2}}  \\ \bbb\in\Lambda_{\frac{1}{2},\frac{1}{6}}}}\!\!\!  (-1)^{a_{2}+b_{2}}(a_{1}+a_{2}\lambda^{\sigma})(b_{1}+b_{2}\lambda^{\sigma}) q^{\frac{1}{g}(a(a_{1}^{2}+b_{1}^{2})-c(a_{2}^{2}+b_{2}^{2}))} \\
	& - \sum_{\substack{\bba\in\Lambda_{\frac{1}{2},\frac{1}{2}}  \\ \bbb\in\Lambda_{0,\frac{1}{6}}}}\!\!\!  (-1)^{a_{2}+b_{2}}(a_{1}+a_{2}\lambda^{\sigma})(b_{1}+b_{2}\lambda^{\sigma}) q^{\frac{1}{g}(a(a_{1}^{2}+b_{1}^{2})-c(a_{2}^{2}+b_{2}^{2}))}\,.
	\end{align*}

Due to the symmetries of the lattices considered, the $q$-exponents of the terms corresponding to different choices of the signs $\pm a_{1}$ and $\pm b_{1}$ are the same. Moreover, the flip $(a_{1},a_{2};b_{1},b_{2})\mapsto(b_{1},a_{2};a_{1},b_{2})$ gives a bijection between the lattice $\Lambda_{0,\frac{1}{2}}\times\Lambda_{\frac{1}{2},\frac{1}{6}}$ appearing in the first summand and the lattice $\Lambda_{\frac{1}{2},\frac{1}{2}} \times \Lambda_{0,\frac{1}{6}}$ appearing in the second one.

As a consequence, the coefficients of the terms corresponding to $(a_{1},a_{2};b_{1},b_{2})$ and $(-a_{1},a_{2};-b_{1},b_{2})$ in the first lattice and their flipped images $(b_{1},a_{2};a_{1},b_{2})$ and $(-b_{1},a_{2};-a_{1},b_{2})$ in the second one give (up to a $(-1)^{a_{2}+b_{2}}$ factor)
	\begin{align*}
	 (a_{1}b_{1}+a_{2}b_{2}(\lambda^{\sigma})^{2}) + \lambda^{\sigma}(a_{1}b_{2}+a_{2}b_{1})
	+&  (a_{1}b_{1}+a_{2}b_{2}(\lambda^{\sigma})^{2}) - \lambda^{\sigma}(a_{1}b_{2}+a_{2}b_{1}) \\
    - (a_{1}b_{1}+a_{2}b_{2}(\lambda^{\sigma})^{2}) - \lambda^{\sigma}(b_{1}b_{2}+a_{1}a_{2})
	-&  (a_{1}b_{1}+a_{2}b_{2}(\lambda^{\sigma})^{2}) + \lambda^{\sigma}(b_{1}b_{2}+a_{1}a_{2})\,,
	\end{align*}
which sums up to zero, and therefore $\cG_{D}$ vanishes along $F_{\cP}$.

In order to determine the vanishing order, we will study the highest order $k$ such that all the $k$-derivatives of $\cG_{D}$ vanish along $F_{\cP}$. The Fourier expansions of the restriction of the derivatives $\partial^{k}\cG_{D}/\partial\tau_{1}^{k}$ and $\partial^{k}\cG_{D}/\partial\tau_{2}^{k}$ to the Shimura curve $F_{\cP}$ are given by the same series as above, with the coefficients replaced by
	\[(-1)^{a_{2}+b_{2}}(a_{1}+a_{2}\lambda^{\sigma})(b_{1}+b_{2}\lambda^{\sigma})(\rho(\bba)^{2}+\rho(\bbb)^{2})^{k}\]
in the case of $\partial^{k}\cG_{D}/\partial\tau_{1}^{k}$ and the equivalent expression with $(\rho^{\sigma}(\bba)^{2}+\rho^{\sigma}(\bbb)^{2})^{k}$ for $\partial^{k}\cG_{D}/\partial\tau_{2}^{k}$.

The coefficients of $\partial\cG_{D}/\partial\tau_{1}$ corresponding to $(a_{1},a_{2};b_{1},b_{2})$ and $(-a_{1},a_{2};-b_{1},b_{2})$ in the first lattice and their flipped images $(b_{1},a_{2};a_{1},b_{2})$ and $(-b_{1},a_{2};-a_{1},b_{2})$ in the second one are given this time by
	\begin{align*}
	  (a_{2}\lambda^{\sigma}+a_{1}) \cdot (b_{2}\lambda^{\sigma}+b_{1}) \cdot &
	  \left[ (a_{1}^{2}+b_{1}^{2}+a_{2}^{2}\lambda^{2}+b_{2}^{2}\lambda^{2}) + 2\lambda (a_{1}a_{2} + b_{1}b_{2}) \right] \\
	+ (a_{2}\lambda^{\sigma}-a_{1}) \cdot (b_{2}\lambda^{\sigma}-b_{1})\cdot &
	\left[ (a_{1}^{2}+b_{1}^{2}+a_{2}^{2}\lambda^{2}+b_{2}^{2}\lambda^{2}) - 2\lambda (a_{1}a_{2} + b_{1}b_{2}) \right] \\
    - (a_{2}\lambda^{\sigma}+a_{1}) \cdot (b_{2}\lambda^{\sigma}+b_{1}) \cdot &
    \left[ (a_{1}^{2}+b_{1}^{2}+a_{2}^{2}\lambda^{2}+b_{2}^{2}\lambda^{2}) + 2\lambda (b_{1}a_{2} + a_{1}b_{2}) \right] \\
	- (a_{2}\lambda^{\sigma}-a_{1}) \cdot (b_{2}\lambda^{\sigma}-b_{1}) \cdot &
	\left[ (a_{1}^{2}+b_{1}^{2}+a_{2}^{2}\lambda^{2}+b_{2}^{2}\lambda^{2}) - 2\lambda (b_{1}a_{2} + a_{1}b_{2}) \right] \,,
	\end{align*}
which again sums up to zero. The same calculation for the derivative $\partial\cG_{D}/\partial\tau_{2}$ shows that it is zero too, and the vanishing order of $\cG_{D}$ along $F_{\cP}$ is therefore at least 2.

Finally, a simple but long calculation shows that the minimum coefficient of $\partial^{2}\cG_{D}/\partial\tau_{1}^{2}$, given by the terms corresponding to $\bba=(0,\pm\tfrac{1}{2})$ and $\bbb=(\pm\tfrac{1}{2}, \tfrac{1}{6})$ in the first lattice and $\bba=(\pm\tfrac{1}{2},\pm\tfrac{1}{2})$ and $\bbb=(0, \tfrac{1}{6})$ in the second one, is $-\tfrac{4}{27}(-1)^{2/3}\lambda^{2}(\lambda^{\sigma})^2$.
\end{proof}

\subsection{Vanishing order along \texorpdfstring{$G_{D}$}{G_D}}

The modular form $\cG_{D}$ vanishes along $G_{D}(\frakb)$ by construction. We next prove that it vanishes only with multiplicity~$1$.
\par
\begin{prop} \label{prop:GothicOnGD} If the Gothic modular form $\cG_{D}$ vanishes to order $>1$ at $\tau_{0}$, then $\tau_{0}\in\tRED(\frakb)$. In particular, $\cG_{D}$ vanishes to order 1 along the Gothic \Teichmuller curve $G_{D}(\frakb)$.
\end{prop}

\begin{proof}Assume that $\cG_{D}$ vanishes to order strictly larger than 1 at a point $\bm\tau_{0}\in G_{D}(\frakb)$, so that in particular $\bm\tau_{0}\in X_{D}\setminus\tRED(\frakb)$. This is equivalent to both derivatives $\partial\cG_{D}/\partial\tau_{1}(\bm\tau_{0})$ and $\partial\cG_{D}/\partial\tau_{2}(\bm\tau_{0})$ being zero.

Recall the theta functions $\theta_{X}(\bm\tau,\bbu)$, $\Theta_{a}(\bm\tau,\bbu)$ and $\Theta_{b}(\bm\tau,\bbu)$ defined at the beginning of Section~\ref{sec:vanGothic}, and the fact that $\div\theta_{X}(\bm\tau_{0})=\varphi(X)$ is the pre-Abel-Prym image of a curve $X$ in the Gothic locus.

By~\eqref{eq:derivativesD2}, $\cG_{D}(\bm\tau)$ is proportional to $D_{2}\Theta_{a}(\bm\tau,\tfrac{1}{2}\lambda_{2})$ and $D_{2}\Theta_{b}(\bm\tau,\tfrac{1}{2}\lambda_{2})$ and therefore, by the heat equation (see~\cite[Proposition 8.5.5]{birkenhakelange}) one has
	\[\frac{\partial}{\partial\tau_{2}}\cG_{D}(\bm\tau_{0}) = \frac{\partial^{3}}{\partial u_{2}^{3}}\Theta_{a}(\bm\tau_{0},\bbu)|_{\bbu=\tfrac{1}{2}\lambda_{2}} = \frac{\partial^{3}}{\partial u_{2}^{3}}\Theta_{b}(\bm\tau_{0},\bbu)|_{\bbu=\tfrac{1}{2}\lambda_{2}} \,.\]

In particular, since all the lower order $u_{2}$-derivatives of $\theta_{X}$ vanish, one has
	\[\frac{\partial^{3}}{\partial u_{2}^{3}}\theta_{X}\left(\tfrac{1}{2}\lambda_{2}\right)=\frac{\partial^{3}}{\partial u_{2}^{3}}\Theta_{a}\left(\tfrac{1}{2}\lambda_{2}\right)\cdot\cF_{b} - \frac{\partial^{3}}{\partial u_{2}^{3}}\Theta_{b}\left(\tfrac{1}{2}\lambda_{2}\right)\cdot\cF_{a}=0\,.
	\]
Therefore, the differential $du_{1}$ induces an abelian differential on $X$ with two double zeroes at $\tfrac{1}{2}\mu_{2}$ and $\tfrac{1}{2}(\lambda_{2}+\mu_{2})$ and a zero of order $\ge 3$ at $\tfrac{1}{2}\lambda_{2}$, which is a contradiction to~$X$ having genus 4.
\end{proof}

Note that we have proved that not even the $\tau_{2}$-derivative of $\cG_{D}$ vanishes anywhere along $G_{D}(\frakb)$. This gives actually a direct proof of the following fact, without knowing that the curves originate as \Teichmuller curves.

\begin{cor}\label{cor:Kobayashi} The vanishing locus of $\cG_D$ is a
  union of Kobayashi geodesics.
\end{cor}
\par
\begin{proof}Being a Kobayashi geodesic is equivalent to always being
transversal to one of the two natural foliations of $X_{D}(\frakb)$
(see~\cite[Proposition~1.3]{moellerprym}), hence modular curves are obviously Kobayashi geodesics and the non-vanishing of
the derivative $\partial/\partial \tau_{2}\cG_{D}(\bm\tau)$ anywhere
in $G_{D}(\frakb)$ proves the statement.
\end{proof}
\par
Finally, the following result shows that the reducible locus agrees indeed with the locus $\tRED(\frakb)$ defined in Theorem~\ref{thm:GDmodform}
\par
\begin{prop}\label{prop:RedisRed}The two definitions of the reducible locus
in $X_D(\frakb)$ agree, that is
	\[\RED(\frakb)=\{\cG_{D}(\bm{\tau}) = 0 \} \cap \{\cF_a(\bm{\tau}) = 0 \} \cap \{ \cF_b(\bm{\tau}) = 0\}\,.\]
\end{prop}

\begin{proof}By the previous proposition, the only thing left to prove is that the intersection on the right hand side is included in the reducible locus.

Let $\bm\tau_{0}\in \{\cG_{D}(\bm{\tau}) = 0 \} \cap \{\cF_a(\bm{\tau}) = 0 \} \cap \{ \cF_b(\bm{\tau}) = 0\}$. Assume without loss of generality that $\Theta_{a}$ is non-zero, otherwise take $\Theta_{b}$. This theta function satisfies:
	\begin{itemize}
	\item $D_{2}\Theta_{a}(0)=D_{1}\Theta_{a}(0)=0$ by definition and by $\cF_{a}(\bm\tau_{0})=0$ respectively;
	\item $D_{2}\Theta_{a}(\tfrac{1}{2}\mu_{2})=D_{1}\Theta_{a}(\tfrac{1}{2}\mu_{2})=0$ by translation to zero via~\eqref{eq:derivativesD2};
	\item $D_{2}\Theta_{a}(\tfrac{1}{2}\lambda_{2})=D_{2}\Theta_{a}(\tfrac{1}{2}(\lambda_{2}+\mu_{2}))=0$, both by translation to zero via~\eqref{eq:derivativesD2} and $\cG_{D}(\bm\tau_{0})=0$.
	\end{itemize}
Hence the theta function $\Theta_{a}$ satisfies all the conditions of $\theta_{X}$ in the proof of Theorem~\ref{thm:GDmodform}, and additionally $D_{1}\Theta_{a}(\tfrac{1}{2}\mu_{2})=0$. As a consequence $Y=\div\Theta_{a}$ is a divisor with self-intersection $Y^{2}=12$ by Riemann-Roch, and multiplicity 3 at the origin and $\tfrac{1}{2}\mu_{2}$. Moreover, since at least the first and second (by odd parity) $u_{2}$-derivatives vanish at $\tfrac{1}{2}\lambda_{2}$ and $\tfrac{1}{2}(\lambda_{2}+\mu_{2})$, either $du_{1}$ induces an abelian differential with zeroes of order $\ge 2$ at those points, or the multiplicity of $Y$ at them is $\ge 3$. The same analysis as in the proof of Theorem~\ref{thm:GDmodform} concludes that the only option is $T_{\bm\tau_{0}}\in\RED(\frakb)$ (case $Y=3Y_{1}+2Y_{2}$). 
\par
Note that the case $T_{\bm\tau_{0}}\in G_{D}$ (case $Y=\varphi(X)$ reduced with zero as its only singular point) is not possible due to the extra vanishing of $D_{1}\Theta_{a}(\tfrac{1}{2}\mu_{2})$, which implies multiplicity $\ge 3$ at that point.
\end{proof}

%% file: sec_examples.tex
\section{Modular embedding of $G_{12}$}

This section is independent of the rest of the paper and illustrates
the parametrization of the Gothic locus in the language of the modular
embeddings. We illustrate this for $D=12$, the unique case where
$G_D$ is a triangle curve and, therefore, the methods of hypergeometric
differential equations are available.
\par
A modular embedding for the Fuchsian group~$\Gamma$ with quadratic
invariant trace field~$K$ is a map $\tau \mapsto (\tau, \varphi(\tau))$ from~$\bH$
to $\bH^2$ such that $\varphi(\gamma \tau) = \gamma^\sigma \varphi(\tau)$.
The universal covering of a map $C \to X_D(\frakb)$ from a Teichm\"uller
curve~$C$ with quadratic trace field to the corresponding Hilbert modular
surface gives rise to a modular embedding, see e.g.\ \cite{MZ} for more
details.
\par
The \emph{hypergeometric differential equation} with parameters $(a,b,c)\in \mathbb{R}$ is given by
\begin{equation}\label{eq:HGDE}
L(a,b,c)(y)\=t(1-t)y'' + \left(c-\left(a+b+1\right)t\right)y'-aby=0\,.
\end{equation}
Whenever $\tfrac{1}{l}=|1-c|$, $\tfrac{1}{m}=|c-a-b|$ and $\tfrac{1}{n}=|a-b|$
for some $l,m,n\in \mathbb{Z}\cup\{\infty\}$ satisfying $1/l + 1/m + 1/n < 1$,
the monodromy group of this equation is the Fuchsian triangle group
$\Delta(l,m,n)$. If $l=\infty$, i.e.\ if $c=1$, the space of solutions
of~\eqref{eq:HGDE} near $t=0$ is generated by $y_{1}(t)$ and $\log(t)y_{1}(t)+y_{2}(t)$ where
\begin{align*}
	y_{1}(t) &= F(a,b,c;t)\coloneqq \sum_{n=0}^{\infty}\frac{(a)_{n}(b)_{n}}{(c)_{n}n!}\,t^{n} \ \mbox{and}\\
	y_{2}(t) &= \sum_{r=0}^{\infty}\frac{(a)_{r}(b)_{r}}{(c)_{r}r!}\left(\sum_{k=0}^{r-1}\frac{1}{a+k}+\frac{1}{b+k}-\frac{2}{c+k}\right)t^{r}\,,
\end{align*}
Here $(x)_{n}$ denotes the Pochhammer symbol and $F$ is the \emph{hypergeometric function} with coefficients $(a,b,c)$ for $c=1$.
\par
By Proposition~\ref{prop:G12summary}, the Veech group of $G_{12}$ is the
triangle group $\Delta=\Delta(\infty,3,6)$, generated by the matrices
	\[M_{\infty} = \left(\begin{array}{cc}
		1 & \alpha \\
		0 & 1
	\end{array}\right),\ 
	M_{3} = \left(\begin{array}{cc}
		\frac{1}{2} & \frac{\sqrt{3}}{2} \\
		-\frac{\sqrt{3}}{2} & \frac{1}{2}
	\end{array}\right),\ 	
	M_{6} = \left(\begin{array}{cc}
		\sqrt{3} + \frac{1}{2} & \frac{5\sqrt{3}}{6} + 1 \\
		-\frac{\sqrt{3}}{2} & -\frac{1}{2}
\end{array}\right)\,,\]
where $\alpha = \frac{2\sqrt{3}}{3} + 2$. It is therefore the monodromy
group of the hypergeometric differential equation 
$L\coloneqq L(\tfrac{5}{12},\tfrac{1}{4},1)=0$
corresponding to $(\tfrac{1}{l},\tfrac{1}{m},\tfrac{1}{n})=(0,\tfrac{1}{3},\tfrac{1}{6})$ and we let $y_1(t)$ and $y_2(t)$ be the functions 
defined above. We will also identify the quotient $\Delta\backslash\bH$ with $\bP^{1}$ via the function $t:\bH\to\bP^{1}$ sending the  elliptic generators $M_{\infty}$, $M_{3}$, and $M_{6}$ of $\Delta(\infty,3,6)$ to 0, 1 and $\infty$ respectively.
\par
Given that the invariant trace field of $\Delta$ is $\bQ(\sqrt{3})$, we will
also be interested in the ``conjugate'' differential equation corresponding to the triangle
group $\Delta^{\sigma}$, for the non-trivial element $\sigma\in\mathrm{Gal}(\bQ(\sqrt{3}))$. One can see that the rotation numbers of the generators of
$\Delta^{\sigma}$ of order 3 and 6 are $e^{4\pi i/3}$ and $e^{2\pi i/6}$,
respectively. As a consequence, the differential equation associated to the
group $\Delta^{\sigma}$ is $
L_{\sigma} \coloneqq L(\tfrac{1}{4},\tfrac{1}{12},1)=0\,,$
which corresponds to $(\tfrac{1}{l_{\sigma}},\tfrac{1}{m_{\sigma}},\tfrac{1}{n_{\sigma}})=(0,\tfrac{2}{3},\tfrac{1}{6})$. We denote the corresponding functions defining the solutions of $L_{\sigma}$ 
by $\widetilde{y}_1(t)$ and $\widetilde{y}_2(t)$.
\par
By~\cite[Formula~(52)]{MZ} the modular embedding~$\varphi$ is given
in terms of these solutions by
\be \label{eq:modemb}
\varphi(\tau) \= \frac{\alpha^\sigma}{\alpha}\tau \,+\, \frac{\alpha^{\sigma}}{2\pi i}\,
\Bigl(\log\frac{A}{\widetilde{A}} \,+\, \frac{\widetilde{f}_2(\tau)}
{\widetilde{f}_1(\tau)} \,+\, \frac{f_2(\tau)}{f_1(\tau)}\Bigr)\,.
\ee
for the constants~$A$ and $\widetilde{A}$, where $f_i(\tau) = y_i(t(\tau))$ and $\widetilde{f}_i(\tau)= \widetilde{y}_i(t(\tau))$. Since $t(\tau)$ is $M_{\infty}$-invariant, 
we can express these functions in terms of the parameter $q=e^{2\pi i\tau/\alpha}$. 
The constants $A$ and $\widetilde{A}$ are determined by $Q(t) \coloneqq te^{y_2(t)/y_1(t)} = Aq$ and $\widetilde{Q}(t) \coloneqq te^{\widetilde{y}_2(t)/\widetilde{y}_1(t)} = \widetilde{A}\widetilde{q}$, where $\widetilde{q}=e^{2\pi i\varphi(\tau)/\alpha^{\sigma}}$. The main remaining task is thus to determine~$A$ and $\widetilde{A}$.

\par

%

Due to the chosen normalization, the function $t(\tau)$ takes the value 1 at the point $i$ with multiplicity 3 and $t(\tau)\neq 1$ whenever $\im(\tau)>1$. It follows that the function $1/(t(\tau)-1)$ has a triple pole at $\tau=i$ and that, as a power series in $q$ (resp. $Q$), the closest singularity to the origin is given by $q_{0}=e^{-2\pi/\alpha}$ (resp. $Q_{0}=Aq_{0}$). This implies that, if one writes $1/(t(Q)-1)^{1/3}=\sum b_{n}Q^{n}$ as a power series in $Q$, the quotients $b_{n}/b_{n+1}$ will tend exponentially fast to $Q_{0}$. This yields a high precision approximation
	\[A \approx 33.9797081543461844465412173813877\ldots\] 
The same calculations for $\widetilde{Q}=\widetilde{A}\widetilde{q}$ yield
	\[\widetilde{A} \approx 3254.6483182744669365311774168770392\ldots\]
These constants can be recognized as the ``conjugate-in-exponent'' pair
	\bas
	A &\= (2+\sqrt{3})^{-6-\sqrt{3}}(1+\sqrt{3})^9(3+\sqrt{3})^3\quad\mbox{and} \\
	\widetilde{A} &\= (2+\sqrt{3})^{-6+\sqrt{3}}(1+\sqrt{3})^9(3+\sqrt{3})^3\,.
	\eas
The resulting modular embedding from formula~\eqref{eq:modemb} is approximately
	\bas
	\varphi(\tau) \=& -\frac{2(1-\sqrt{3})}{\pi i}\log (2+\sqrt{3}) + (2-\sqrt{3})\tau + \frac{3+\sqrt{3}}{3\pi i} \left(-\frac{1}{6} \, Aq -\frac{5}{1152} \, A^{2}q^{2} \right. \\ 
	& \left. -\frac{61}{497664} \, A^{3}q^{3} -\frac{713}{382205952} \, A^{4}q^{4} -\frac{4943}{183458856960} \, A^{5}q^{5}+\ldots \right) \\ 
	\= & 0.963i + 0.268 \,\tau + 9.243\cdot 10^{-7}i \,\tau^2 -1.159\cdot 10^{-6}\,\tau^3 + 8.389 \cdot 10^{-7}i\,\tau^4 \\
	& -2.136 \cdot 10^{-7}\,\tau^5 + 4.611\cdot 10^{-7}i \,\tau^6 + 9.035 \cdot 10^{-8} \,\tau^7 + 1.630 \cdot 10^{-7}i \,\tau^8 \\
	& + 1.053 \cdot 10^{-7} \,\tau^9 + 2.502 \cdot 10^{-8}i \,\tau^{10} 5.408 \cdot 10^{-8} \,\tau^{11} + \ldots
	\eas 
and the modular transformation can be checked numerically.

\par

Finally, note that the group $\Delta$ does not belong to the Hilbert modular group $\SL(\frakb\oplus\cO_{12}^{\vee})$, but the conjugate $\Delta_{C}=C\Delta C^{-1}$ by the matrix 
	\[C=\left(\begin{matrix} 3\sqrt{3}+9 & -3\sqrt{3}-15 \\ 0 & 1 \end{matrix}\right)\] does. Consequently, the map $\tau \mapsto (\tau, C^{\sigma}\circ\varphi\circ C^{-1}(\tau))$, where matrices act on $\bH$ by M\"{o}bius transformations, parametrizes the \Teichmuller curve $G_{12}=\Delta_{C}\backslash\bH \to X_{12}(\frakb)$. Indeed, it can be numerically checked that the image of this map lies in the vanishing locus of the modular form $\cG_{12}$ and, since $\RED(\frakb)$ is empty in this case, it actually equals $\{\cG_{12}(\bm\tau)=0\}$.

%% file: sec_asymptotics.tex
\section{Asymptotics of divisor sums} \label{sec:AsyDiv}

As preparation for computing the asymptotics of volumes and Lyapunov exponents
in the next section, we study here for a fundamental discriminant $D$ the
asymptotics as $D \to \infty$ of 
\bes
e(D,k) \= \sum_{b^2 \equiv D\mod 4k, \atop |b| \leq \sqrt{D}}
\sigma_1\Bigl(\frac{D-b^2}{4k}\Bigr) \= \sum_{[a,b,c] \in \cP_k(D)} a\,,
\ees
where $\sigma_1(\cdot)$ is the divisor sum function and
where $\cP_k(D)$ has been introduced in~\eqref{eq:defPkD}.
Our focus is on the cases $k=1$ and $k=6$, but the method works
for general~$k$.
\par
\begin{theorem} \label{thm:eDk}
The following asymptotic statements hold:
\bas
e(D,1) &\= \phantom{\frac{1}{50}}\frac{\zeta_{\bQ(\sqrt{D})}(-1)}{2 \zeta(-3)}
+ O(D^{5/4}) && \phantom{for}\\
e(D,6) &\=  \frac{1}{50} \frac{\zeta_{\bQ(\sqrt{D})}(-1)}{2 \zeta(-3)}
+ O(D^{5/4}) && \text{for}\quad D \equiv 0,12 \mod 24 \\
e(D,6) &\= \frac{2}{50} \frac{\zeta_{\bQ(\sqrt{D})}(-1)}{2 \zeta(-3)}
+ O(D^{5/4}) && \text{for}\quad D \equiv 4,9,16 \mod 24 \\
e(D,6) &\= \frac{4}{50} \frac{\zeta_{\bQ(\sqrt{D})}(-1)}{2 \zeta(-3)}
+ O(D^{5/4})  && \text{for}\quad D \equiv 1 \mod 24 \\
\eas
as $D \to \infty$ among fundamental discriminants.
\end{theorem}
\par
Note that $\zeta_{\bQ(\sqrt{D})}(-1) > C D^{3/2}$, so the theorem captures
indeed the asymptotics for large~$D$. Our proof here follows closely an
application of the circle method used by Zagier in \cite[Section~4]{DonNZeta}.
To set the stage, we define the one-variable theta series and the Eisenstein
series to be the modular forms
\bes
\theta(\tau) \= \sum_{\ell = -\infty}^\infty e^{\pi i \ell^2 \tau},
\quad
G_{2}(\tau) \= -\frac{1}{24} 
 \,+\,
\sum_{a=1}^\infty \sigma_{1}(a) e^{2\pi i a\tau} \,.
\ees
Then the modular form 
\bes
F(\tau,k) \,:=\, G_2(2k\tau) \theta(\tau)
\= \sum_{n=0}^\infty e(n,k) e^{\pi i n\tau} \,.
\ees
has a Fourier expansion with coefficients that generalize
the coefficients we are interested in. The basic idea idea is to compute
the Fourier coefficients of $F(\tau,k)$ by integration at small height~$\ve$.
The dominating term of the asymptotics then comes from the expansions near
each rational point. Consequently, we use the modular transformation law
to obtain the expansions
\ba \label{eq:THGtransf}
\theta\Bigl(\frac{a}c +iy\Bigr) &\= \lambda(a,c)(cy)^{-1/2} + O(y^{-1/2}
e^{-\pi/4 c^2y})\,, \\
G_2 \Bigl(\frac{a}c +iy\Bigr) & \= -\zeta(2)(cy)^{-2} + O(y^{-2}
e^{-\pi/c^2y})
\ea
as $y \to \infty$ where $a,c \in \bZ$ with $\gcd(a,c)=1$ and where
$\lambda(a,c)$ is a Legendre symbol times a power of~$i$, depending
on the parities of~$a$ and~$c$. Here we mainly need to know that
the Gauss sum
\bes
\gamma_c(n) = c^{-1/2} \sum_{a=1}^{2c} \lambda(a,c) e^{-\pi ina/c}
\ees
is computed in \cite[Theorem~2]{DonNZeta} for $D$ fundamental  to be
a weakly multiplicative function in~$c$ given on prime powers by
\be \label{eq:gacdef}
2^r \mapsto 
\begin{cases}
1 & \text{if}\,\, r \in \{0,1\} \\
2 \chi(2) & \text{if}\,\, r=2 \\
2 & \text{if}\,\, r=3 \,\,\text{and}\,\, 2|D \\
0 & \text{otherwise} \\
\end{cases} 
\quad \text{and} \quad 
p^r \mapsto 
\begin{cases}
1 & \text{if}\,\, r=0 \\
\chi(p) & \text{if}\,\, r=1 \\
-1 & \text{if}\,\, r=2 \,\,\text{and}\,\, p|D \\
0 & \text{otherwise} \\
\end{cases} 
\ee
for odd primes~$p$, where $\chi(m) = (\frac{D}m)$. Define
\be \label{eq:defstarenk}
e^\ast(n,k) \= \sum_{c=1}^\infty \frac{\gcd(c,2k)^2}{c^{2}}\gamma_c(n)\,.
\ee
\par
\begin{lemma} \label{le:east}
For $k$ square-free and $D$ a fundamental discriminant
\bes
\frac{e^\ast(D,k)}{e^\ast(D,1)} \= \prod_{p | k\,\,\text{prime}} \frac{1+\chi(p)}{1+p^{-2}}\,.
\ees
\end{lemma}
\par
\begin{proof}
Since the summands in~\eqref{eq:defstarenk} are weakly multiplicative in~$c$,
the function $e^\ast(D,1)$ admits an Euler product expansion. For $p\neq 2$
equation~\eqref{eq:gacdef} directly implies that the local factor is
\bes
1+ \frac{\chi(p)}{p^2} + \frac{\chi(p)^2-1}{p^4} \= \frac{1-p^{-4}}
{1-{\chi(p)}{p^{-2}}}\,.
\ees
For $p=2$ the same conclusion holds up to global factor~$2$ after taking
the factor $\gcd(c,2k)^2$ into account. In total
\be \label{eq:eveast}
e^\ast(D,1) \= 2 \frac{L(2,\chi)}{\zeta(4)}\,,
\ee
where $L(s,\chi) = \zeta_K(s)/\zeta(s)$ is the $L$-series associated with
the character~$\chi$. The passage from $\gcd(c,2)$ to $\gcd(c,2k)$ only
changes the factors at the primes dividing~$k$. For $p \neq 2$ the local
factor now is $1+\chi(p) + (\chi(p)^2-1)p^{-2}$, and the ratio compared
to the original factor results in the modification claimed in the lemma.
For $p=2$ the same final conclusion holds.
\end{proof}
\par
\begin{proof}[Proof of Theorem~\ref{thm:eDk}] Since~$F$ is periodic
under $\tau \mapsto \tau+2$ we can compute the coefficients
\bes e(n,k)\= \frac12 \int_{i\ve}^{2+i\ve} e^{\pi in\tau} F(\tau,k) dy
\ees
using Cauchy's formula by integration at small height~$\ve$. We replace the
right hand side in a neighborhood of $a/c \in [0,2)$ by the dominating term in
\bes
F\Bigl(\frac{a}c +iy ,k\Bigr) \= \frac{\zeta(2)}{16\pi^2}
\lambda(a,c)\frac{\gcd(c,2k)^2}{k^2\, c^{5/2}}y^{-5/2} +O(y^{-5/2}e^{-\pi/4 c^2y} )\,
\ees
obtained as combination of~\eqref{eq:THGtransf}\,. The sum over all
`major arcs' of the circle method is the summation of these neighborhoods.
It is computed in \cite[Equation~(32)]{DonNZeta} using the integral
representation of the Gamma-function to be 
\be \label{eq:defolenk}
\overline{e}(n,k) \= \frac{\pi^{1/2}\zeta(2)}{16\,\Gamma(5/2)} \,
\frac{n^{3/2}}{k^2}\, \sum_{c=1}^\infty \frac{\gcd(c,2k)^2}{c^{2}}\gamma_c(n)\,.
\ee
To see that the major arcs indeed give the dominating term, we can
argue as in \cite[p.~81]{DonNZeta} (referring to Hardy), for any
fixed~$k$. The equations~\eqref{eq:defolenk} and~\eqref{eq:eveast}
can now be combined as in \cite{DonNZeta} to the case $k=1$ of the
theorem. The cases for $k=6$ differ by the factor $1/k^2$
in~\eqref{eq:defolenk} and the factors in Lemma~\ref{le:east}.
\end{proof}

%% file: sec_lyapunov.tex
\section{Volumes and Lyapunov exponents}\label{sec:volandlyap}

The results of the previous sections can now be assembled to compute
the Euler characteristic of the Gothic \Teichmuller curves and
their Lyapunov exponents. We first state a more precise version of
Theorem~\ref{thm:introVol}. Recall the definition of $\kappa_{D}$ in
Proposition~\ref{prop:vol16}.
\par
\begin{theorem}\label{thm:volumes}
Let $D$ be a non-square discriminant. The Gothic \Teichmuller curve $G_D$ is
non-empty if and only $D \equiv 0,1,4,9,12,16\bmod{24}$. \\ 
For $D \equiv 0,12\bmod{24}$ the Gothic \Teichmuller curve $G_D$ has Euler characteristic
\bes \label{eq:fullvol1comp}
-\chi(G_D) = \frac{1}{20}\,\kappa_{D}\!\!\!  \sum_{[a,b,c]\in \cP_{1}(D)} \!\!\! a\quad  - \frac{1}{3} \sum_{[a,b,c] \in \cP_6(D)} \!\!\! a
\ees
For $D \equiv 4,9,16\bmod{24}$
the Gothic Teichm\"uller curve $G_D = G_D^0 \cup G_D^1$ consists
of two sub-curves $G_D^\epsilon$ of the same volume equal to
\bes \label{eq:fullvol2comp}
-\chi(G_D^\epsilon) = \frac{1}{20}\,\kappa_{D}\!\!\!  \sum_{[a,b,c]\in \cP_{1}(D)} \!\!\! a\quad  - \frac{1}{6} \sum_{[a,b,c] \in \cP_6(D)}  \!\!\! a \,,\qquad \epsilon \in \{0,1\}.
\ees 
For $D \equiv 1\bmod{24}$
there is a decomposition $G_D = G_D^{00} \cup G_D^{01}
\cup G_D^{10} \cup G_D^{01}$ of the  Gothic Teichm\"uller curve into
four sub-curves $G_D^{\epsilon \delta}$ of the same volume
equal to 
\bes \label{eq:fullvol4comp}
-\chi(G_D^{\epsilon \delta}) = \frac{1}{20}\,\kappa_{D}\!\!\! \sum_{[a,b,c]\in \cP_{1}(D)} \!\!\! a\quad - \frac{1}{12} \sum_{[a,b,c] \in \cP_6(D)} \!\!\!a\,,\qquad \epsilon,\delta \in \{0,1\}.
\ees
\end{theorem}
\par
\medskip
To state the other theorems, we provide a brief introduction to Lyapunov
exponents, in particular for flat surfaces $(X,\omega)$ in the Gothic locus.
\par
Lyapunov exponents measure the growth rate of cohomology classes
in $H^1(X,\bR)$ under parallel transport along the geodesic flow in~$\overline{\SL(2,\bR)\cdot(X,\omega)}$, the closure of the $\SL(2,\bR)$-orbit of $(X,\omega)$ 
(see e.g.\ \cite{ZorichFlat} or \cite{MoePCMI} for background). The Lyapunov
spectrum of
a genus four surface consists of Lyapunov exponents $\lambda_1 = 1 \geq
\lambda_2 \geq \lambda_3 \geq \lambda_4$ and their negatives.
\par
In the case of flat surfaces $(X,\omega)$ in the Gothic locus, the existence of the maps $\pi_{A}$ and $\pi_{B}$ in~\eqref{eq:Gothic} decomposes the local system~$\bV$ with fiber $H^1(X,\bR)$
over~$\Omega G$ into local
subsystems $\bV_A$ and $\bV_B$ of rank two, corresponding to the elliptic
curves~$A$ and~$B$,
and the (`Prym') complement~$\bV_P$. Since the generating differential of
the Gothic form belongs to the Prym part, the exponent $\lambda_1=1$ is one
of the two positive exponents $\{1,\lambda_P\}$ of~$\bV_P$. If we denote
by $\lambda_A$ and $\lambda_B$ the Lyapunov exponents from the elliptic curves,
then the sets 
	\[\{1,\lambda_2, \lambda_3, \lambda_4\} \= \{1,\lambda_A, \lambda_B, \lambda_P\}\]
coincide. Since by definition $\omega^2 = \pi_A^* q$ for some quadratic
differential~$q$ on~$A$, the double covering formula of Eskin-Kontsevich-Zorich
(\cite{EKZ}) implies
\bes
\lambda_1 + \lambda_P + \lambda_B - \lambda_A \= \tfrac14 \cdot 3 \cdot
\bigl(\tfrac1{1+2} + \tfrac1{-1+2}\bigr) \= 1
\ees
\par
\begin{theorem} \label{thm:LyapOmegaG}
The Prym Lyapunov exponent~$\lambda_P$ of a generic surface in the
Gothic locus is equal to $3/13$.
\end{theorem}
\par
This is a direct consequence of the asymptotics formulas in Theorem~\ref{thm:eDk},
the following proposition and the convergence of individual Lyapunov exponents
(\cite{bew}), since the curves $G_D$ equidistribute
towards (the Lebesque measure on) the Gothic locus by \cite{esmimo}.
\par
\begin{prop} \label{prop:LyapGD}
The Prym Lyapunov exponent of a
Gothic Veech surface on~$G_D$ is equal to
\bes \lambda_P(G_D(\frakb)) \=  1+\frac{\chi(X_{D}(\frakb))}{\chi(G_{D}(\frakb))} \,.
\ees
\end{prop}
\par
Note that we do not claim that the curves~$G_D(\frakb)$ 
are connected, although we expect this to be true. Therefore, the statement of
the proposition has to be interpreted as volume-weighted average of the~$\lambda_P$
of the connected components.
\par
\begin{proof}[Proof of  Theorem~\ref{thm:volumes}] 
The arguments in the following work for any good compactification $\overline{X_{D}(\frakb)}$ of $X_{D}(\frakb)$ (see~\cite{moellerprym}). Since the specific choice of compactification is not relevant, we will denote simply by $[C]$ the class of the closure $\overline{C}$ in $\overline{X_{D}(\frakb)}$.

Let $[\omega_i]$ be the classes of the two foliations of the Hilbert modular surface
$X_D(\frakb)$. Then the uniformization of $X_D(\frakb)$
implies $\chi(X_D(\frakb)) = [\omega_1]\cdot [\omega_2]$ and the vanishing locus of a modular
for of bi-weight $(k,\ell)$ has class $\tfrac12(k[\omega_1] + \ell[\omega_2])$.
\par
Theorem~\ref{thm:GDmodform} and Proposition~\ref{prop:GothicOnRed} together show that
the vanishing locus of the Gothic modular form~$\cG_D$ is a union of Kobayashi curves
and the first coordinate can be used as parameter for each of these curves. By Theorem~\ref{thm:divGD}, $\div(\cG_{D})=G_{D}(\frakb)+2\RED(\frakb)$, where $G_{D}(\frakb)$ denotes the union of those components of the Torelli-image of $G_{D}$ in $X_{D}(\frakb)$ for which $du_{1}$ induces the eigenform $\omega$ at each point $(X,\omega)$.
\par
This implies that integration of $\omega_1$ along $\div(\cG_D)$ (equivalently, the intersection product $-[\omega_{1}]\cdot[\div(\cG_{D})]$) computes the
sum of the Euler characteristics of these curves with the multiplicity determined
in Propositions~\ref{prop:GothicOnRed} and~\ref{prop:GothicOnGD} (see
\cite[Corollary~10.4]{bainbridge07} or \cite[Proposition~1.3]{moellerprym}). We obtain
\be\label{eq:eulerchar}
- \frac32 \chi(X_D(\frakb)) \= -[\omega_{1}]\cdot (\tfrac12[\omega_1] + \tfrac32[\omega_2]) = \chi(G_D(\frakb)) + 2\, \chi(\RED(\frakb))\,.
\ee
Proposition~\ref{prop:vol16} together with the well-known expression for the Euler characteristic $\chi(X_{D})$ of standard Hilbert modular surfaces in terms of prototypes (see for example~\cite{Hir}) give
	\[\chi(X_{D}(\frakb)) = \frac{1}{30}\,\kappa_{D}\!\!\!  \sum_{[a,b,c]\in \cP_{1}(D)} \!\!\! a\,.\]
Formula~\eqref{eq:eulerchar} together with Lemma~\ref{lem:volRed} proves the result for $G_{D}(\frakb)$. The only thing left to do is to prove the decomposition of $G_{D}$ into sub-curves as claimed.
\par
We claim that, for different ideals $\frakb_{1}$ and $\frakb_{2}$, the images of $G_{D}(\frakb_{1})$ and $G_{D}(\frakb_{2})$ in $\cA_{2,(2,3)}$ are different. In fact, if $\frakb_{2}\neq\frakb_{1}^{\sigma}$, the images of the whole Hilbert modular surfaces $X_{D}(\frakb_{1})$ and $X_{D}(\frakb_{2})$ are disjoint in $\cA_{2,(2,3)}$, since the lattices of the corresponding abelian surfaces are not even isomorphic as $\cO_{D}$-modules. On the other hand, if $\frakb^{\sigma}\neq \frakb$, the sub-curves $G_{D}(\frakb)$ and $G_{D}(\frakb^{\sigma})$ can both be seen in $X_{D}(\frakb)$ as Kobayashi geodesics with $\omega_{1}$ and $\omega_{2}$ as parameters, respectively. In particular, if their images under $X_{D}(\frakb)\to\cA_{2,(2,3)}$ agreed, their associated eigenforms for real multiplication would map to two different eigenforms on each point $X\in G_{D}$.
\par
Finally, by construction $G_{D}$, is covered by the union of the images of $G_{D}(\frakb)$ for the different ideals $\frakb$ of norm 6.
\end{proof}
\par
\begin{proof}[Proof of Proposition~\ref{prop:LyapGD} and Theorem~\ref{thm:LyapOmegaG}] 
The Lyapunov exponent $\lambda_{P}(C)$ of a Kobayashi geodesic $C$ in $X_{D}(\frakb)$ is given by the following quotient (see~\cite{bainbridge07} or~\cite{moellerprym})
	\[\lambda_{P}(C)=\frac{[\omega_{2}]\cdot [C]}{[\omega_{1}]\cdot [C]}\,.\]
The reducible locus $\RED(\frakb)$ is a union of Shimura curves, and therefore one has $\lambda_{P}(\RED(\frakb))=1$ and $-[\omega_{1}]\cdot[\RED(\frakb)]=-[\omega_{2}]\cdot[\RED(\frakb)]=\chi(\RED(\frakb))$. In the case of the Gothic Teichm\"{u}ller curves, since $[G_{D}(\frakb)]=[\div(\cG_{D})]-2[\RED(\frakb)]$, one has
	\[\lambda_{P}(G_{D}(\frakb)) = \frac{[\omega_{2}]\cdot\left(\tfrac12[\omega_1] + \tfrac32[\omega_2] -2[\RED(\frakb)]\right) }{[\omega_{1}]\cdot\left(\tfrac12[\omega_1] + \tfrac32[\omega_2] -2[\RED(\frakb)]\right) }= \frac{\tfrac12\chi(X_{D}(\frakb)) + 2\chi(\RED(\frakb))}{\tfrac32\chi(X_{D}(\frakb)) + 2\chi(\RED(\frakb))}\,.\]
By Theorem~\ref{thm:introVol}, this is exactly $1+\chi(X_{D}(\frakb))/\chi(G_{D}(\frakb))$.


Theorem~\ref{thm:LyapOmegaG} follows by taking the limit and using Theorem~\ref{thm:volumes}
and Theorem~\ref{thm:eDk}
\end{proof}

%% file: sec_tables.tex
\section{Tables of invariants}\label{sec:tables}

\begin{table}[htp!]\centering
\hspace{-2cm}\begin{minipage}[t]{0.4\linewidth}
 \begin{tabular}[t]{| l | c | c | c | c |}  \hline $D$ & $\#$ 
 & $ \chi(X_{D}(\frakb))$ & $ \chi(\RED(\frakb))$ & $ \chi(G_{D}^{\epsilon\delta})$ \\ [0.5ex]  \hline  
12 & 1 & 1/3 & 0 & $-1/2\,^{*}$ \\ 
24 & 1 & 1 & -1/6 & $-7/6$ \\ 
28 & 2 & 4/3 & -1/6 & $-5/3$ \\ 
33 & 2 & 2 & -1/6 & $-8/3\,^{*}$ \\ 
40 & 2 & 7/3 & -1/6 & $-19/6\,^{*}$ \\ 
48 & 1 & 4 & -1/2 & $-5$ \\ 
52 & 2 & 5 & -1/2 & $-13/2$ \\ 
57 & 2 & 14/3 & -1/2 & $-6\,^{\dag}$ \\ 
60 & 1 & 4 & -1/3 & $-16/3$ \\ 
72 & 1 & 20/3 & -2/3 & $-26/3\,^{*}$ \\ 
73 & 4 & 22/3 & -2/3 & $-29/3\,^{\dag}$ \\ 
76 & 2 & 19/3 & -2/3 & $-49/6\,^{*}$ \\ 
84 & 1 & 10 & -1 & $-13\,^{\dag}$ \\ 
88 & 2 & 23/3 & -5/6 & $-59/6$ \\ 
96 & 1 & 12 & -1 & $-16$ \\ 
97 & 4 & 34/3 & -7/6 & $-44/3$ \\ 
105 & 2 & 12 & -4/3 & $-46/3$ \\ 
108 & 1 & 12 & -4/3 & $-46/3$ \\ 
112 & 2 & 16 & -3/2 & $-21\,^{\dag}$ \\ 
120 & 1 & 34/3 & -1 & $-15\,^{*}$ \\ 
124 & 2 & 40/3 & -7/6 & $-53/3\,^{*}$ \\ 
129 & 2 & 50/3 & -3/2 & $-22\,^{*}$ \\ 
132 & 1 & 18 & -2 & $-23$ \\ 
136 & 2 & 46/3 & -5/3 & $-59/3$ \\ 
145 & 4 & 64/3 & -13/6 & $-83/3$ \\ 
148 & 2 & 25 & -5/2 & $-65/2\,^{\dag}$ \\ 
153 & 2 & 80/3 & -8/3 & $-104/3\,^{*}$ \\ 
156 & 1 & 52/3 & -2 & $-22$ \\ 
160 & 2 & 28 & -3 & $-36\,^{\dag}$ \\ 
168 & 1 & 18 & -5/3 & $-71/3$ \\ 
172 & 2 & 21 & -2 & $-55/2\,^{*}$ \\ 
177 & 2 & 26 & -5/2 & $-34$ \\ 
180 & 1 & 40 & -4 & $-52\,^{\dag}$ \\ 
\hline
    \end{tabular}%
\end{minipage}\hspace{2.5cm}
\begin{minipage}[t]{0.4\linewidth}
    \begin{tabular}[t]{| l | c | c | c | c |}  \hline $D$ & $\#$ 
& $ \chi(X_{D}(\frakb))$ & $ \chi(\RED(\frakb))$ & $ \chi(G_{D}^{\epsilon\delta})$ \\ [0.5ex]  \hline 
184 & 2 & 74/3 & -7/3 & $-97/3\,^{*}$ \\ 
192 & 1 & 32 & -3 & $-42\,^{\dag}$ \\ 
193 & 4 & 98/3 & -10/3 & $-127/3$ \\ 
201 & 2 & 98/3 & -7/2 & $-42$ \\ 
204 & 1 & 26 & -8/3 & $-101/3$ \\ 
208 & 2 & 40 & -4 & $-52$ \\ 
216 & 1 & 36 & -10/3 & $-142/3$ \\ 
217 & 4 & 116/3 & -23/6 & $-151/3$ \\ 
220 & 2 & 92/3 & -10/3 & $-118/3$ \\ 
228 & 1 & 42 & -4 & $-55$ \\ 
232 & 2 & 33 & -7/2 & $-85/2$ \\ 
240 & 1 & 48 & -5 & $-62$ \\ 
241 & 4 & 142/3 & -14/3 & $-185/3$ \\ 
244 & 2 & 55 & -11/2 & $-143/2$ \\ 
249 & 2 & 46 & -9/2 & $-60$ \\ 
252 & 1 & 128/3 & -4 & $-56\,^{*}$ \\ 
264 & 1 & 112/3 & -4 & $-48$ \\ 
265 & 4 & 160/3 & -31/6 & $-209/3$ \\ 
268 & 2 & 41 & -4 & $-107/2$ \\ 
273 & 2 & 148/3 & -5 & $-64$ \\ 
276 & 1 & 60 & -6 & $-78$ \\ 
280 & 2 & 134/3 & -13/3 & $-175/3$ \\ 
288 & 1 & 80 & -8 & $-104\,^{\dag}$ \\ 
292 & 2 & 66 & -7 & $-85$ \\ 
297 & 2 & 72 & -22/3 & $-280/3$ \\ 
300 & 1 & 130/3 & -4 & $-57$ \\ 
304 & 2 & 76 & -15/2 & $-99$ \\ 
312 & 1 & 46 & -5 & $-59$ \\ 
313 & 4 & 200/3 & -41/6 & $-259/3$ \\ 
316 & 2 & 56 & -11/2 & $-73$ \\ 
321 & 2 & 66 & -13/2 & $-86$ \\ 
328 & 2 & 54 & -5 & $-71$ \\ 
336 & 1 & 80 & -8 & $-104$ \\
\hline \end{tabular}%
\end{minipage}\vspace{0.3cm}
\caption{Number of $\cO_{D}$-ideals $\frakb$ of norm 6 and volumes of each $X_{D}(\frakb)$, $\RED(\frakb)$ and $G_{D}(\frakb)$, for $D\le 385$. The cross and the asterisk indicate a Gothic or hexagons model, respectively.}
\end{table}
